\documentclass[10pt, 
twoside]{amsart}

\title{The Braids on your Blanket}
\date{\today}
\author{Michelle Cheng}
\author{Robert Laugwitz}
\address[R.~L.]{School of Mathematical Sciences,
University of Nottingham, University Park, Nottingham, NG7 2RD, UK}
\email{robert.laugwitz@nottingham.ac.uk}

\usepackage{mathabx}
\usepackage{import} 
\usepackage{amsmath}
\usepackage{amsfonts}
\usepackage{amsthm}
\usepackage{amssymb}
\usepackage[alphabetic, initials, nobysame]{amsrefs}
\usepackage[english]{babel}
\usepackage{url}
\usepackage{fancyhdr}
\usepackage{graphicx}
\usepackage[
colorlinks=true,
linkcolor=black, 
anchorcolor=black,
citecolor=black,
urlcolor=black, 
]{hyperref}
\usepackage[dvipsnames]{xcolor}

\usepackage{geometry}
\usepackage{caption, subcaption}

\newcommand{\mR}{\mathbb{R}}

\newcommand{\mZ}{\mathbb{Z}}

\newcommand{\cK}{\mathcal{K}}
\newcommand{\cL}{\mathcal{L}}

\newcommand{\rJ}{\mathrm{J}}
\newcommand{\rP}{\mathrm{P}}

\newcommand{\labelb}{\xrightarrow[\sim]{\text{Figure~\ref{figure-reidemeister}{\scshape (b)}}}}
\newcommand{\labela}{\xrightarrow[\sim]{\text{Figure~\ref{figure-reidemeister}{\scshape (a)}}}}

\newtheoremstyle{mystyle}
  {0.5cm}                   
  {0.5cm}                   
  {\normalfont}           
  {}                      
  {\itfont\bfseries}  
  {:}                     
  {0.3cm}              
  {\thmname{#1}}

\newtheoremstyle{defstyle}
  {0.5cm}                   
  {0.5cm}                   
  {\normalfont}           
  {}     
  {\normalfont\bfseries}  
  {:}                     
  {0.3cm}              
  {\thmname{#1}\thmnumber{ #2}\thmnote{ (#3)}}

\newtheorem{theorem}{Theorem}
\newtheorem{question}{Question}

\dedicatory{For Leo}

\begin{document}

\begin{abstract}
In this expositional essay, we introduce some elements of the study of groups by analysing the braid pattern on a knitted blanket. We determine that the blanket features pure braids with a minimal number of crossings. Moreover, we determine polynomial invariants associated to the links obtained by closing the braid patterns of the blanket. 
\end{abstract}
\maketitle



\begin{figure}[htbp]
\includegraphics[scale=0.5]{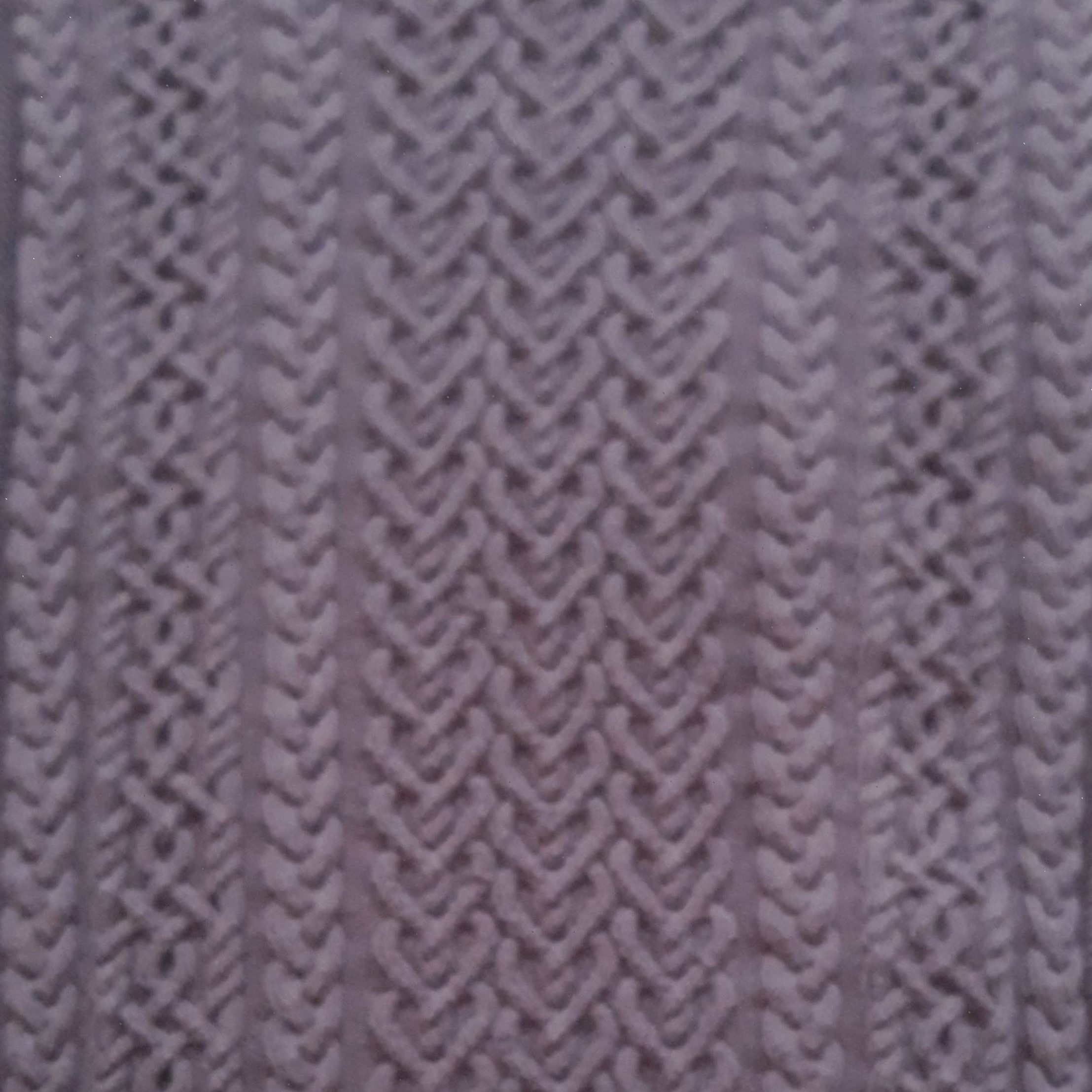}
\caption{The blanket}
\label{figure-blanket}
\end{figure}

Dear Leo, Michelle knitted this baby blanket\footnote{Based on the pattern \emph{Levi's Baby Blanket} found in the following blog \url{https://knittikins.wordpress.com/patterns/levis-baby-blanket/}.} for you which has an interesting design of crossing strands, called braids. The blanket is shown in Figure \ref{figure-blanket}. 
Braids like these have a long history of being used as a decorative element, for example, in Celtic art. In the Celtic tradition, braided patterns symbolize continuity and endless braids (which we will refer to as knots) are said to symbolize the eternity of life \cite{Bai}. Braids continue to be appreciated for their aesthetics in various cultures. 

In this text, we will use mathematics to describe the artistic pattern of the blanket. 
Mathematicians have developed precise formalism to describe symmetries and patterns beyond the study of numbers. 
Such a formalism has, for example, been developed for braids. Mathematicians have studied braids for at least one century and their work has found surprising applications in geometry and physics.
The pattern of the blanket appears symmetric to the eye and therefore feels aesthetically pleasing.  This impression can be confirmed by analyzing the mathematical properties of the braid patterns in the blanket, revealing interesting symmetries among the braids featured.

Mathematicians study braids through the more general concept of \emph{groups}. In mathematics, the word `group' is used to describe the general features of a collection of symmetries. This concept of groups appears ubiquitously in mathematics and physics. Groups first emerged from the work of a French mathematician, Evangeliste Galois, in the 19th century and revolutionized mathematicians' understanding of solving equations by making use of symmetries among solutions. The story of Galois is interesting by itself; it relates to the French revolution and provides a thrilling read. We will not discuss Galois contribution and the history of groups and how they help in solving equations here.\footnote{The book of Simon Singh \cite{Sin} gives an accessible account of the story of a mathematical problem that was solved after eluding mathematicians for over 358 years and relates to the work of great minds throughout the history of mathematics.}

Many believe that all important results in mathematics have already been discovered. However, the opposite is true. While some old questions (like Fermat's Last Theorem) have been solved, many problems remain unsolved and new questions continue to emerge. In fact, solutions to old problems often present new problems that nobody imagined thinking about  before. We illustrate this phenomenon in Section \ref{questions-math} by discussing such problems that emerged from the study of braids. 

\section{Groups}

\subsection{What are groups?}\label{section-groups}

Groups provide a helpful way to describe symmetries of mathematical structure that appear in different types of contexts. To define the mathematical concept of a group, key features that the collections symmetries of several structures have in common are abstracted to so-called \emph{axioms}. These axioms capture certain fundamental aspects of the nature of a collection symmetry. Mathematicians often approach defining an abstract concept this way: They identify a list of fundamental axioms shared by many structures and then refer to \emph{any} structure having these properties by a name---in this case, we call these structures \emph{groups}. Other examples of such abstract structures include number fields, functions, relations, differential operators, or probability distributions. By giving a universally accepted name to each of widely used these concepts, mathematics is a language that describes abstract and systematic structures that can be used to describe the world.

\begin{figure}
\centering
\begin{subfigure}[htb]{0.32\textwidth}
\centering
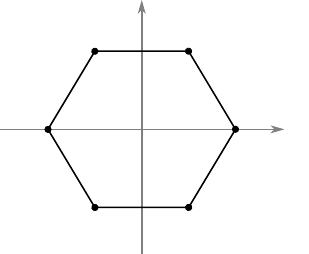
\caption{The rotation $x$}
\label{figure-hexagonx}
\end{subfigure}
\begin{subfigure}[htb]{0.32\textwidth}
\centering
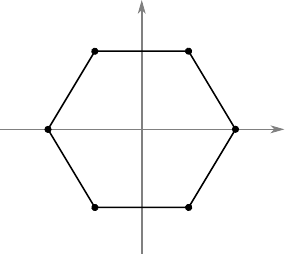
\caption{The reflection $y$}
\label{figure-hexagony}
\end{subfigure}
\begin{subfigure}[htb]{0.32\textwidth}
\centering
\begingroup%
  \makeatletter%
  \providecommand\color[2][]{%
    \errmessage{(Inkscape) Color is used for the text in Inkscape, but the package 'color.sty' is not loaded}%
    \renewcommand\color[2][]{}%
  }%
  \providecommand\transparent[1]{%
    \errmessage{(Inkscape) Transparency is used (non-zero) for the text in Inkscape, but the package 'transparent.sty' is not loaded}%
    \renewcommand\transparent[1]{}%
  }%
  \providecommand\rotatebox[2]{#2}%
  \newcommand*\fsize{\dimexpr\f@size pt\relax}%
  \newcommand*\lineheight[1]{\fontsize{\fsize}{#1\fsize}\selectfont}%
  \ifx\svgwidth\undefined%
    \setlength{\unitlength}{136.48035268bp}%
    \ifx\svgscale\undefined%
      \relax%
    \else%
      \setlength{\unitlength}{\unitlength * \real{\svgscale}}%
    \fi%
  \else%
    \setlength{\unitlength}{\svgwidth}%
  \fi%
  \global\let\svgwidth\undefined%
  \global\let\svgscale\undefined%
  \makeatother%
  \begin{picture}(1,0.89298088)%
    \lineheight{1}%
    \setlength\tabcolsep{0pt}%
    \put(0,0){\includegraphics[width=\unitlength,page=1]{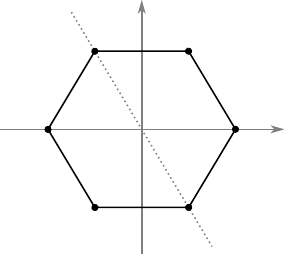}}%
    \put(0.69083703,0.76782733){\makebox(0,0)[lt]{\lineheight{1.25}\smash{\begin{tabular}[t]{l}$1$\end{tabular}}}}%
    \put(0.88322845,0.45736968){\makebox(0,0)[lt]{\lineheight{1.25}\smash{\begin{tabular}[t]{l}$2$\end{tabular}}}}%
    \put(0.69366886,0.08091454){\makebox(0,0)[lt]{\lineheight{1.25}\smash{\begin{tabular}[t]{l}$3$\end{tabular}}}}%
    \put(0.28252531,0.76725545){\makebox(0,0)[lt]{\lineheight{1.25}\smash{\begin{tabular}[t]{l}$6$\end{tabular}}}}%
    \put(0.09024561,0.45679764){\makebox(0,0)[lt]{\lineheight{1.25}\smash{\begin{tabular}[t]{l}$5$\end{tabular}}}}%
    \put(0.28535714,0.0803425){\makebox(0,0)[lt]{\lineheight{1.25}\smash{\begin{tabular}[t]{l}$4$\end{tabular}}}}%
    \put(0.52986959,0.27908038){\makebox(0,0)[lt]{\lineheight{1.25}\smash{\begin{tabular}[t]{l}${\color{Purple}xy}$\end{tabular}}}}%
    \put(0,0){\includegraphics[width=\unitlength,page=2]{hexagonxy.pdf}}%
    \put(0.28369712,0.55689978){\makebox(0,0)[lt]{\lineheight{1.25}\smash{\begin{tabular}[t]{l}${\color{Purple}xy}$\end{tabular}}}}%
  \end{picture}%
\endgroup%

\caption{The product $xy$}
\label{figure-product}
\end{subfigure}
\caption{Symmetries of a regular hexagon: $x$ is a rotation, $y$ is a reflection, and the  product $xy$ gives the reflection about the axis connecting $3$ and $6$}
\label{figure-hexagonbig}
\end{figure}

Certain groups can be described as collections of symmetries. Such a group is a collection of symmetries of a mathematical structure, for example, a geometric object, the set of solutions to an equation, or the possible states of a physical system, may all display symmetries. 
As a first example, we will explain the group of symmetries of a hexagon depicted in Figure \ref{figure-hexagonbig}. 

A symmetry of the hexagon is a transformation of the plane which sends any point on the lines describing the hexagon to another (or possibly the same) point on the hexagon. This means, any symmetry preserves the shape even though the different points on it might be swapped. By labelling the corners of the hexagon with numbers 1 to 6 as in Figure \ref{figure-hexagonbig}, we can describe all of its symmetries effectively. 

Examples of symmetries of the hexagon are \emph{rotations}. For instance, Figure~\ref{figure-hexagonbig}{\scshape (a)} shows the rotation that rotates corners 1 to 2, 2 to 3, 3 to 4, etc. For brevity, we can call this rotation, $x$. But there are other rotations such as the described by the following mapping of corners: 
\begin{equation}\label{eq-x4}
1\mapsto 5, \quad 2\mapsto 6, \quad 3\mapsto 1, \quad 4\mapsto 2, \quad 5\mapsto 3, \quad 6 \mapsto 4.
\end{equation}
If we apply the rotation $x$ four times, we obtain the above rotation described in \eqref{eq-x4}. 

Rotations are not the only symmetries of the hexagon, there are also \emph{reflections}. For example, there is a reflection that interchanges corners $1$ and $6$, $2$ and $5$, $3$ and $4$. This reflection is displayed in Figure~\ref{figure-hexagonbig}{\scshape (b)} where we call this reflection, $y$. It is the reflection about a vertical axis. Another reflection swaps corners $1$ and $5$, $2$ and $4$, while fixing corners $3$ and $6$---we may name this reflection, $z$. 
This reflection, $z$, is the reflection about the axis going through the corners $3$ and $6$.
In principle, we can reflect the hexagon at any axis that has an angle which is a multiple of $360^\circ/6=60^\circ$ degrees. 
The collection of all symmetries of the hexagon is an example of a group, it carries the name, the \emph{dihedral group} and is often denoted by $D_{6}$ since it captures symmetries of the hexagon which has six corners. One can demonstrate that there are precisely six distinct reflections and six distinct rotations for this hexagon. 

A group has a set of elements but also an operation that takes two elements as input and gives a single element as output, the \emph{product} of the two input elements. We usually use the symbol $x\cdot y$ to represent the product of the elements $x$ and $y$ of a group. In the case of the group of symmetries of the hexagon, the dihedral group, $D_6$, the product of $x$ and $y$ is given by composition of symmetries, $x\cdot y$. This means that we will first apply the symmetry $x$ followed by applying the symmetry $y$. The result could have been achieved through a single symmetry---their product, $x\cdot y$. 

Here is an example of a product in the dihedral group with $12$ elements, $D_{6}$. It is the product of, first, the rotation $x$ and, second, the reflection $y$, both of which we considered above. We can compute the product, which is denoted by $x\cdot y$---like a product of numbers---by tracing what happens to the corners of the hexagon. For example, $x$ sends corner $1$ to corner $2$, and then  $y$ sends corner $2$ to corner $5$. Hence, the product $x\cdot y$ has sent $1$ to $5$. Similarly, $x$ sends corner $2$ to $3$, and then $y$ sends corner $3$ to $4$. Therefore, the product $x\cdot y$ has sent corner $2$ to $4$. Continuing this way, the following Table \ref{table-product} describes the product $x\cdot y$ completely. Figure \ref{figure-hexagonbig}{\scshape (c)} shows us that this product is a reflection about the axis passing through $3$ and $6$, which are therefore fixed points of $x\cdot y$. We observe that the previously defined reflection $z$ coincides with the product $x\cdot y$ and may write $z=x\cdot y$. The order of $x$ and $y$ in the product $x\cdot y$ is crucial here as $y\cdot x$ does not give the same resulting symmetry. The last column in Table~\ref{table-product} shows the effect of the product $y\cdot x$ on the six corners. For instance, $y$ sends corner $1$ to $6$, and $y$ sends corner $6$ to $1$, so the symmety $y\cdot x$ sends corner $1$ to $1$, meaning, it fixes corner $1$.
\begin{table}[htb]
\[
\begin{array}{cccc}
x&y&x\cdot y&y\cdot x\\\hline
1\mapsto 2& 1\mapsto 6 & 1\mapsto 5&1\mapsto 1\\
2\mapsto 3& 2\mapsto 5 & 2\mapsto 4&2\mapsto 6\\
3\mapsto 4& 3\mapsto 4 & 3\mapsto 3&3\mapsto 5\\
4\mapsto 5& 4\mapsto 3 & 4\mapsto 2&4\mapsto 4\\
5\mapsto 6& 5\mapsto 2 & 5\mapsto 1&5\mapsto 3\\
6\mapsto 1& 6\mapsto 1 & 6\mapsto 6&6\mapsto 2
\end{array}
\]
\caption{Computing the products $x\cdot y$ and $y\cdot x$ in the dihedral group $D_{6}$}
\label{table-product}
\end{table}

Another example of a product computation was shown earlier in \eqref{eq-x4}. Here, we computed the product of $x$ with itself four times. This means, we computed $x^4=x\cdot x\cdot x\cdot x$ and the outcome of this computation is the rotation shown in \eqref{eq-x4}.

\subsection{Examples of groups}

Symmetry groups are a type of groups, but the concept of a group is much more universal. The range of numbers we usually compute with, called the \emph{real numbers}, form a group where the binary product operation is the addition of numbers. If $a$ and $b$ are numbers, then $a+b$ is also a number. Mathematicians call this group, $\mR$. Similarly, the whole numbers form a group called the \emph{integers} which is denoted by $\mZ$. Here, the product is also given by addition. There are many other groups. One of the easiest ones has exactly two elements, $0$ and $1$. The product operation, indicated by the symbol $+$ instead of $\cdot$ as it is a form of addition, is specified by the following list:
$$0+0=0, \qquad 0+1=1, \qquad 1+0=1, \qquad 1+1=0.$$
This might seem surprising, but as there are only two elements, we have no choice but to set $1+1=0$. This can be thought of as addition of whole numbers but only retaining the information whether a number is odd or even, also referred to as the \emph{parity} of the numbers. Then, even numbers are represented by $0$ and odd numbers are represented by $1$. This is justified by $0,2,4,\ldots$ being even numbers and $1,3,5$ being odd numbers. You can test the following rules of parity for the addition of numbers: 
$$\text{even}+\text{even}=\text{even}, \quad \text{even}+\text{odd}=\text{odd}, \quad \text{odd}+\text{even}=\text{odd},\quad \text{odd}+\text{odd}=\text{even}.$$
This two-element group is denoted by $\mZ_2$ among mathematicians. It is like the integers $\mZ$ but only has two elements.

A different way to describe the same group $\mZ_2$ with two elements is by denoting one element by $1$ and the other element by $-1$. The product operation works like the multiplication of integers and is indicated with the symbol $\cdot$. The following list describes all products for this two-element group:
$$ 1\cdot 1=1, \quad 1\cdot (-1)=-1, \quad (-1)\cdot 1=-1, \quad (-1)\cdot (-1)=1.$$
This is another incarnation of the same group $\mZ_2$ because we can match the element $1$ with $0$ and $-1$ with $1$ and the values of the operations $\cdot$ and $+$ correspond to one another under this matching. 

\subsection{The axioms of a group}

Axioms are used to provide a universal definition of what constitutes a group. These axioms are a short list of three rules that can be checked to determine if a given structure constitutes a group. If a statement about groups can be derived from the these fundamental axioms alone, it will hold for \emph{all} examples groups at once. This idea of axiomatic logic enables mathematicians to prove statements with absolute certainty.

The definition of a group involves three axioms. The first axiom is called \emph{associativity}. It states that for a product of three elements $x,y$, and $z$ of a group, we can ignore brackets. In formulas, this is expressed by the equality
$$(x\cdot y)\cdot z=x\cdot (y\cdot z).$$
It does not matter whether we form the product $x\cdot y$ and then form the product of $x\cdot y$ and $z$, which is $(x\cdot y)\cdot z$, or first form the product $y\cdot z$ and then the product of $x$ and $y\cdot z$, which is $x\cdot (y\cdot z)$. 

However, the order of elements appearing in a product does matter. In many groups, like the symmetries of the hexagon, $x\cdot y$ is not the same as $y\cdot x$. You can see this by comparing the last two columns of Table~\ref{table-product}. In some groups, such as the integers $\mZ$, the real number $\mR$, or the group $\mZ_2$ with two elements, we have $a+b=b+a$. This property is called \emph{commutativity}. Commutativity is not an axiom of a group. Groups satisfying commutativity are rather special, and we often use the addition symbol $+$ to denote a commutative product. 

The second axiom of a group states that any group contains a special element called the \emph{identity element}. This element is neutral with respect to the product operation in that forming the product with the identity element has no effect. In formulas, this means that for any element $x$ of a group,
$$x\cdot e= x, \qquad \text{and } \qquad e\cdot x=x.$$
 For example, when adding numbers, adding $0$ has no effect: $0+a=a=a+0$. Therefore, $0$ is the identity element of numbers with addition. For the dihedral group of symmetries of the hexagon discussed earlier, the neutral element for the composition of symmetries is the rotation by an angle of zero degrees. This trivial rotation fixes all corners. Hence, applying this symmetry before or after applying another symmetry has no effect. Therefore, this trivial symmetry is the identity element of the dihedral group. We usually denote the identity element of a group by $e$. For instance, for the group $\mZ$ of integers with addition, $e=0$.

The third fundamental axiom of a groups is called \emph{invertibility}. This axiom implies that group operations are reversible. Thinking about the symmetries of the hexagon again, we can undo each rotation by rotating back in the opposite direction. We can undo each reflection by applying it once again. This property of invertibility can be expressed, more generally, as follows. For a given element $x$ of a group, there exists an element called the \emph{inverse of $x$} which is denoted by $x^{-1}$. This inverse element satisfies the equations
$$x\cdot x^{-1}=e, \qquad \text{and}\qquad x^{-1}\cdot x=e,$$
where $e$ is the identity element of the group discussed before.
For example, the reflection $y$ that we discussed earlier is equal to its own inverse, i.e., $y^{-1}=y$. This is the case because applying a reflection twice returns the hexagon to its previous configuration of corners. In the notation of groups, $y \cdot  y=e$. 
The inverses for addition of integer numbers are the negatives of the given numbers. In fact, for any integer $a$, we have that
$$a+(-a)= 0, \qquad \text{and}\qquad (-a)+a=0.$$
Since $e=0$, this means $-a$ is the inverse of $a$ for the addition operation.

To summarize, a group is a collection of elements with a product operation that is associative, comes with an identity element, and inverse elements. There are many examples of groups and some of them we already encountered. Another fundamental example of a group is given by the set of all non-zero fractions. A fraction $\tfrac{a}{b}$ is described by two integer numbers, $a$ and $b$, and we require both of them to be non-zero. The product of fractions is given by 
$$\tfrac{a}{b}\cdot \tfrac{c}{d}=\tfrac{ac}{bd}.$$
One can check that this is an associative operation. In fact, commutativity also holds for this product:
$$\tfrac{a}{b}\cdot \tfrac{c}{d}=\tfrac{ac}{bd}=\tfrac{ca}{db}=\tfrac{c}{d}\cdot \tfrac{a}{b}.$$
The fraction $\tfrac{1}{1}$ is the identity element for this product since 
$$\tfrac{1}{1}\cdot \tfrac{a}{b}=\tfrac{1a}{1b}=\tfrac{a}{b}.$$
Inverses are given by reversing fractions, namely
$$\left(\tfrac{a}{b}\right)^{-1}=\tfrac{b}{a}.$$
You can verify this by using cancellation rules for fractions. Indeed,
$$\tfrac{a}{b}\cdot \tfrac{b}{a}=\tfrac{ab}{ba}=\tfrac{a}{a}=\tfrac{1}{1}.$$
What is interesting about the group of fractions is that we require non-zero fractions in order to have inverses of elements with respect to the multiplication operation. For instance, if $a$ is an integer, we can identify it with the fraction $\tfrac{a}{1}$. The inverse of $a$ is 
$$a^{-1}=\left(\tfrac{a}{1}\right)^{-1}=\tfrac{1}{a}.$$
This shows that $a$ cannot be zero as the fraction $\tfrac{1}{0}$ is ill-defined. This example also indicates that the set of non-zero integers is \emph{not} a group with respect to multiplication as there are no multiplicative inverses for integers other than $1$ and $-1$. For example, we cannot find an \emph{integer} $a$ such that 
$2\cdot a=1.$ We know that this equation is correct when $a=\tfrac{1}{2}$ which is \emph{not} a whole number and therefore not an integer.

Now we have encountered some of the most well-known groups. There are two more remarkable series of groups that we will consider---the \emph{braid groups} and the \emph{symmetric groups}. These groups also describe fundamental structures appearing in everyday life in concise and systematic ways.

\subsection{The symmetric groups}\label{section-Sn}

Consider the sequence of the first $n$ numbers $(1,2,\ldots, n)$. The symmetric group contains all possible ways to rearrange the order of these $n$ numbers. For example, if $n=5$, an element of the symmetric group on $5$ numbers is given by the sequence
$$(2, 3, 1, 5,4).$$
Here, we have reshuffled the order of the numbers to start with $2$, followed by $3$, etc. The symmetric group on $n$ numbers is usually denoted by $S_n$. Consequently, the above sequence $(2, 3, 1, 5,4)$ is an element of $S_5$. We can treat such a sequence like a function. For example, $(2, 3, 1, 5,4)$ corresponds to the function sending $1$ to $2$, $2$ to $3$, $3$ to $1$, $4$ to $5$, and finally $5$ to $4$. 

\begin{figure}
     \centering
     \begin{subfigure}[htb]{0.45\textwidth}
     \centering
		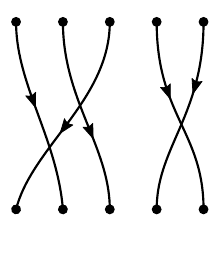 
		\centering
		\caption{}
         \label{figure-S5element-a}
     \end{subfigure}
     \begin{subfigure}[htb]{0.45\textwidth}
     \centering
     	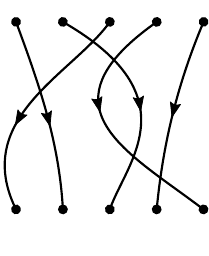
 		\caption{} 
         \label{figure-S5element-b}
     \end{subfigure}
\caption{The element $(2, 3, 1, 5,4)$ of the symmetric group $S_5$---displayed in two distinct ways}
   \label{figure-S5element}
\end{figure}

For the purpose of studying braids later, it helps to visualize the symmetric groups as pictures of crossing strings. Here, it does not matter which string crosses above and which one crosses below. These string pictures are read from top to bottom. The above sequence $(2, 3, 1, 5,4)$ corresponds to the picture in Figure~{\scshape\ref{figure-S5element}(a)}. To read such a diagram, start at one of the numbered strands at the top (say, the first one) and trace the corresponding string through the picture along the direction of the arrow. Then record the number at which this strand terminated (in this example, the number $2$). Then, we record this number in the first entry of the sequence. To verify the entire sequence, we trace all numbers from the top to the bottom along the direction of the arrows of the respective strands.  There is always some ambiguity about how to draw a certain sequence as a picture of crossing strings. For example, Figure {\scshape\ref{figure-S5element}(b)} is a valid picture for the same sequence $(2, 3, 1, 5,4)$.

The symmetric groups have a product structure which is given by composition of reorderings of the first $n$ numbers. To provide an example, we again look at the sequences $(2, 3, 1, 5,4)$ and view it as a function which, for example, sends $1$ to $2$.
If we are given a second sequence, say, $(3,5,2,1,4)$ then we can compute the product as the composition of the two functions. That means, we will apply the first sequence followed by the second one. In this example, first, $1$ is sent to $2$ by the first sequence, and then $2$ is sent to $5$ by the second sequence. Hence, the composition sends $1$ to $5$. Similarly, $2$ is sent to $3$ by the first sequence, and $3$ is sent to $2$ by the second sequence. Therefore, the composition sends $2$ to itself. Continuing like this we compute that the product is the sequence $(5,2,3,4,1)$. Written as a formula,
\begin{equation}\label{equation-S5ab}
(2, 3, 1, 5,4)(3,5,2,1,4)=(5,2,3,4,1).
\end{equation}

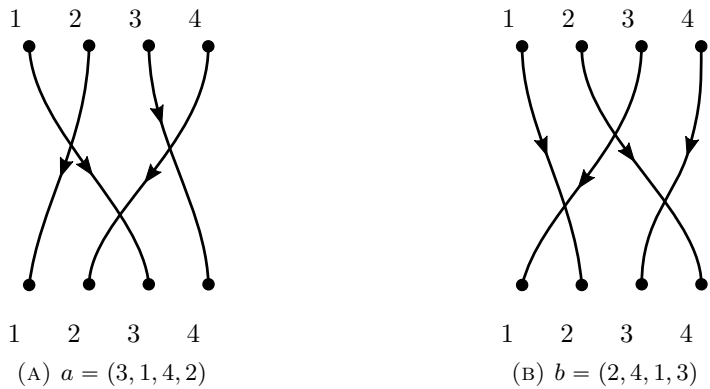
\begin{figure}
     \centering
     \begin{subfigure}[htb]{0.4\textwidth}
     \centering
\begingroup%
  \makeatletter%
  \providecommand\color[2][]{%
    \errmessage{(Inkscape) Color is used for the text in Inkscape, but the package 'color.sty' is not loaded}%
    \renewcommand\color[2][]{}%
  }%
  \providecommand\transparent[1]{%
    \errmessage{(Inkscape) Transparency is used (non-zero) for the text in Inkscape, but the package 'transparent.sty' is not loaded}%
    \renewcommand\transparent[1]{}%
  }%
  \providecommand\rotatebox[2]{#2}%
  \newcommand*\fsize{\dimexpr\f@size pt\relax}%
  \newcommand*\lineheight[1]{\fontsize{\fsize}{#1\fsize}\selectfont}%
  \ifx\svgwidth\undefined%
    \setlength{\unitlength}{77.62798681bp}%
    \ifx\svgscale\undefined%
      \relax%
    \else%
      \setlength{\unitlength}{\unitlength * \real{\svgscale}}%
    \fi%
  \else%
    \setlength{\unitlength}{\svgwidth}%
  \fi%
  \global\let\svgwidth\undefined%
  \global\let\svgscale\undefined%
  \makeatother%
  \begin{picture}(1,1.57875475)%
    \lineheight{1}%
    \setlength\tabcolsep{0pt}%
    \put(0,0){\includegraphics[width=\unitlength,page=1]{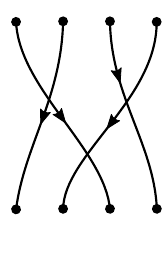}}%
    \put(0.00262383,1.53960572){\makebox(0,0)[lt]{\lineheight{1.25}\smash{\begin{tabular}[t]{l}$1$\end{tabular}}}}%
    \put(0.29246775,1.53960572){\makebox(0,0)[lt]{\lineheight{1.25}\smash{\begin{tabular}[t]{l}$2$\end{tabular}}}}%
    \put(0.58231187,1.53960572){\makebox(0,0)[lt]{\lineheight{1.25}\smash{\begin{tabular}[t]{l}$3$\end{tabular}}}}%
    \put(0.87215579,1.53960572){\makebox(0,0)[lt]{\lineheight{1.25}\smash{\begin{tabular}[t]{l}$4$\end{tabular}}}}%
    \put(-0.00427721,0.00757317){\makebox(0,0)[lt]{\lineheight{1.25}\smash{\begin{tabular}[t]{l}$1$\end{tabular}}}}%
    \put(0.28556671,0.00757317){\makebox(0,0)[lt]{\lineheight{1.25}\smash{\begin{tabular}[t]{l}$2$\end{tabular}}}}%
    \put(0.57541084,0.00757317){\makebox(0,0)[lt]{\lineheight{1.25}\smash{\begin{tabular}[t]{l}$3$\end{tabular}}}}%
    \put(0.86525504,0.00757317){\makebox(0,0)[lt]{\lineheight{1.25}\smash{\begin{tabular}[t]{l}$4$\end{tabular}}}}%
  \end{picture}%
\endgroup%
 
		\centering
		\caption{$a=(3,1,4,2)$}
         \label{figure-S4element-b}
     \end{subfigure}
     \begin{subfigure}[htb]{0.4\textwidth}
     \centering
\begingroup%
  \makeatletter%
  \providecommand\color[2][]{%
    \errmessage{(Inkscape) Color is used for the text in Inkscape, but the package 'color.sty' is not loaded}%
    \renewcommand\color[2][]{}%
  }%
  \providecommand\transparent[1]{%
    \errmessage{(Inkscape) Transparency is used (non-zero) for the text in Inkscape, but the package 'transparent.sty' is not loaded}%
    \renewcommand\transparent[1]{}%
  }%
  \providecommand\rotatebox[2]{#2}%
  \newcommand*\fsize{\dimexpr\f@size pt\relax}%
  \newcommand*\lineheight[1]{\fontsize{\fsize}{#1\fsize}\selectfont}%
  \ifx\svgwidth\undefined%
    \setlength{\unitlength}{77.62928265bp}%
    \ifx\svgscale\undefined%
      \relax%
    \else%
      \setlength{\unitlength}{\unitlength * \real{\svgscale}}%
    \fi%
  \else%
    \setlength{\unitlength}{\svgwidth}%
  \fi%
  \global\let\svgwidth\undefined%
  \global\let\svgscale\undefined%
  \makeatother%
  \begin{picture}(1,1.57872839)%
    \lineheight{1}%
    \setlength\tabcolsep{0pt}%
    \put(0,0){\includegraphics[width=\unitlength,page=1]{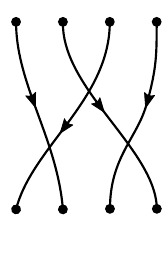}}%
    \put(0.00262378,1.53958002){\makebox(0,0)[lt]{\lineheight{1.25}\smash{\begin{tabular}[t]{l}$1$\end{tabular}}}}%
    \put(0.29246286,1.53958002){\makebox(0,0)[lt]{\lineheight{1.25}\smash{\begin{tabular}[t]{l}$2$\end{tabular}}}}%
    \put(0.58230215,1.53958002){\makebox(0,0)[lt]{\lineheight{1.25}\smash{\begin{tabular}[t]{l}$3$\end{tabular}}}}%
    \put(0.87214124,1.53958002){\makebox(0,0)[lt]{\lineheight{1.25}\smash{\begin{tabular}[t]{l}$4$\end{tabular}}}}%
    \put(-0.00427714,0.00757305){\makebox(0,0)[lt]{\lineheight{1.25}\smash{\begin{tabular}[t]{l}$1$\end{tabular}}}}%
    \put(0.28556195,0.00757305){\makebox(0,0)[lt]{\lineheight{1.25}\smash{\begin{tabular}[t]{l}$2$\end{tabular}}}}%
    \put(0.57540124,0.00757305){\makebox(0,0)[lt]{\lineheight{1.25}\smash{\begin{tabular}[t]{l}$3$\end{tabular}}}}%
    \put(0.8652406,0.00757305){\makebox(0,0)[lt]{\lineheight{1.25}\smash{\begin{tabular}[t]{l}$4$\end{tabular}}}}%
  \end{picture}%
\endgroup%

 		\caption{$b=(2,4,1,3)$} 
         \label{figure-S4element-a}
     \end{subfigure}
\caption{Two elements $a$ and $b$ of the symmetric group $S_4$}
   \label{figure-S4elements}
\end{figure}

\begin{figure}
     \centering
     	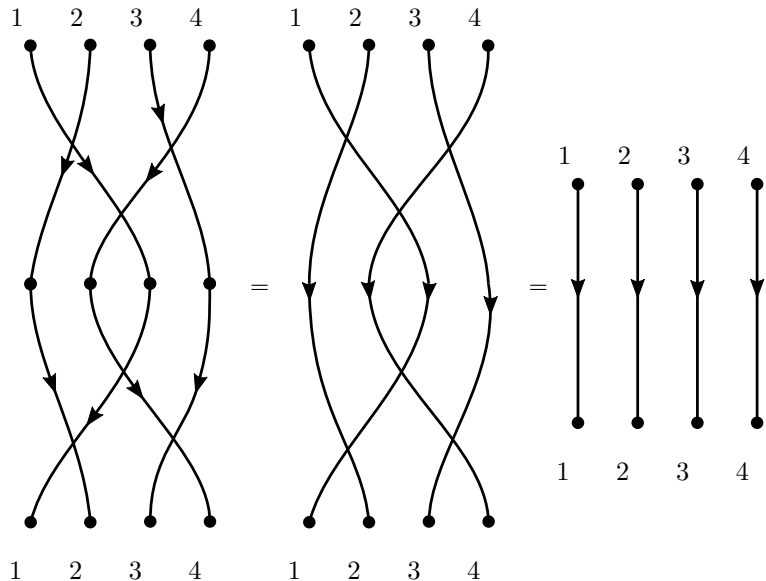
 		\caption{The product $ab$ of the elements $a=(3,1,4,2)$ and $b=(2,4,1,3)$ of the symmetric group $S_4$ (see Figure~\ref{figure-S4elements}) gives the identity element $(1,2,3,4)$} 
         \label{figure-S4element-ba}
\end{figure}

Another example of a product of two elements of symmetric groups is carried out using diagrams in Figure~\ref{figure-S4element-ba}. Here, we consider the elements $a=(3,1,4,2)$ and $b=(2,4,1,3)$ in $S_4$ depicted in Figure~\ref{figure-S4elements}. That is, we consider sequences of length four.
Visually, the product $ab$ is given by stacking $a$ on top of $b$ and removing the labels in the middle. Then we can simplify the resulting longer strings to strings of the original length. All that matters here is the information of where the input numbers come out in the string diagram. We see in Figure~\ref{figure-S4element-ba} how the product $ab$ is computed graphically. The product ends up being the identity element $e=(1,2,3,4)$ of $S_4$ which permutes none of the numbers. With the discussion at the end of the previous section, this means that $b=a^{-1}$ and $a=b^{-1}$, i.e., $a$ and $b$ are mutually inverse elements in $S_4$. In other words, applying $b$ after $a$ undoes the rearrangement of the numbers $1,2,3,4$ caused by $a$.

The groups $S_n$ appear quite different than groups of numbers (such as the integer or real numbers with addition, or non-zero rational numbers with multiplication) that we are familiar with from calculations. In some sense, the symmetric groups, $S_n$, are much more like the dihedral group $D_{6}$ of symmetries of the regular hexagon which we encountered in the previous section. In fact, the group $S_n$ can be interpreted as the group of all symmetries of a list of $n$ points. Another similarity between $S_n$ and $D_{6}$ is that these type of groups are not commutative and only have a finite number of elements. To see that, for example, $S_5$ is not commutative, we  compute
\begin{equation}\label{equation-S5ba}
(3,5,2,1,4)(2, 3, 1, 5,4)=(1,4,3,2,5).
\end{equation}
This product is different from the product $(2, 3, 1, 5,4)(3,5,2,1,4)$ computed in Equation \eqref{equation-S5ab}.

The symmetric groups are universal. Every group with finitely many elements is contained in some symmetric group, $S_n$, as a subgroup. This means that any group is a subset of $S_n$ which contains the identity element and is closed under taking products and inverses.

\section{The Braid Groups}

\subsection{What are the braid groups?}

\begin{figure}
     \centering
     \begin{subfigure}[htb]{0.4\textwidth}
     \centering
$$\vcenter{\hbox{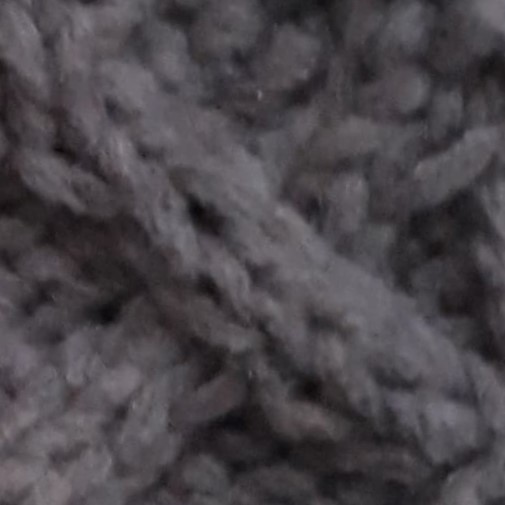}}\quad \longleftrightarrow \quad \vcenter{\hbox{\includegraphics[scale=0.1]{Graphics/psi.jpg}}}
$$
		\centering
		\caption{The left strand crossing above}
         \label{figure-Phi}
     \end{subfigure}
     \begin{subfigure}[htb]{0.4\textwidth}
     \centering
$$\vcenter{\hbox{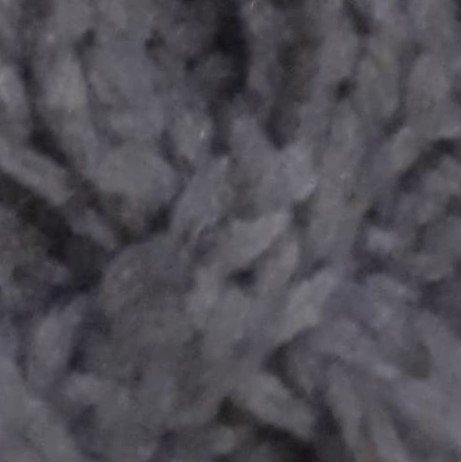}}\quad \longleftrightarrow \quad \vcenter{\hbox{\includegraphics[scale=0.11]{Graphics/psiinv.jpg}}}
$$ 		\caption{The right strand crossing above} 
         \label{figure-Phi-inv}
     \end{subfigure}
\caption{Crossing strands and how they appear in the blanket}
   \label{figure-crossings}
\end{figure}

In the previous section, we studied some examples of groups. The main groups of interest in this article, however, are the \emph{braid groups}. The braid groups can be seen as more complex symmetric groups as they track how the numbers were rearranged, rather than just the order of the numbers. If we look at a picture like Figure~\ref{figure-S5element}{\scshape (a)}, we see strings crossing, but in the crossings, there is no distinction which strand overlaps above and which one lays below at a crossing. But to study braids, the order of crossing strands becomes important. The pictures in Figure~\ref{figure-crossings} show this difference. 

\begin{figure}[htb]
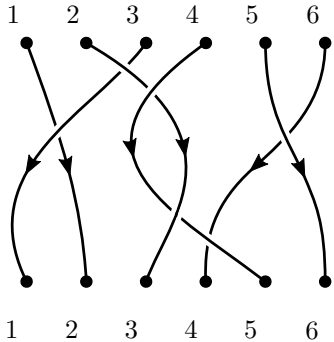
\caption{An example of a braid on $6$ strands}
\label{figure-B6}
\end{figure}

Similarly to the symmetric groups, braid groups form a series of groups indexed by a positive integer number $n$. The braid group on $n$ strands is denoted by $B_n$. Its elements are pictures of $n$ braided strands. For example, Figure \ref{figure-B6} shows an example of an element of $B_6$. In this example, six strings each connect one of the incoming labels $1,2,3,4,5,6$ to one of the outgoing labels (reading from top to bottom). The strings are embedded into three-dimensional space and cannot intersect. This way, the strings braid past one another.

\begin{figure}
     \centering
     \begin{subfigure}[htb]{0.4\textwidth}
     \centering
\begingroup%
  \makeatletter%
  \providecommand\color[2][]{%
    \errmessage{(Inkscape) Color is used for the text in Inkscape, but the package 'color.sty' is not loaded}%
    \renewcommand\color[2][]{}%
  }%
  \providecommand\transparent[1]{%
    \errmessage{(Inkscape) Transparency is used (non-zero) for the text in Inkscape, but the package 'transparent.sty' is not loaded}%
    \renewcommand\transparent[1]{}%
  }%
  \providecommand\rotatebox[2]{#2}%
  \newcommand*\fsize{\dimexpr\f@size pt\relax}%
  \newcommand*\lineheight[1]{\fontsize{\fsize}{#1\fsize}\selectfont}%
  \ifx\svgwidth\undefined%
    \setlength{\unitlength}{79.64151985bp}%
    \ifx\svgscale\undefined%
      \relax%
    \else%
      \setlength{\unitlength}{\unitlength * \real{\svgscale}}%
    \fi%
  \else%
    \setlength{\unitlength}{\svgwidth}%
  \fi%
  \global\let\svgwidth\undefined%
  \global\let\svgscale\undefined%
  \makeatother%
  \begin{picture}(1,1.35049798)%
    \lineheight{1}%
    \setlength\tabcolsep{0pt}%
    \put(0,0){\includegraphics[width=\unitlength,page=1]{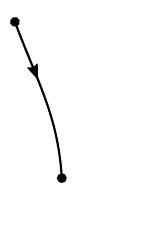}}%
    \put(-0.00416907,1.31233873){\makebox(0,0)[lt]{\lineheight{1.25}\smash{\begin{tabular}[t]{l}$1$\end{tabular}}}}%
    \put(0.27834689,1.31233873){\makebox(0,0)[lt]{\lineheight{1.25}\smash{\begin{tabular}[t]{l}$2$\end{tabular}}}}%
    \put(0.56086284,1.31233873){\makebox(0,0)[lt]{\lineheight{1.25}\smash{\begin{tabular}[t]{l}$3$\end{tabular}}}}%
    \put(0.84337866,1.31233873){\makebox(0,0)[lt]{\lineheight{1.25}\smash{\begin{tabular}[t]{l}$4$\end{tabular}}}}%
    \put(0.00356521,0.0073817){\makebox(0,0)[lt]{\lineheight{1.25}\smash{\begin{tabular}[t]{l}$1$\end{tabular}}}}%
    \put(0.27162033,0.0073817){\makebox(0,0)[lt]{\lineheight{1.25}\smash{\begin{tabular}[t]{l}$2$\end{tabular}}}}%
    \put(0.55413628,0.0073817){\makebox(0,0)[lt]{\lineheight{1.25}\smash{\begin{tabular}[t]{l}$3$\end{tabular}}}}%
    \put(0.83665237,0.0073817){\makebox(0,0)[lt]{\lineheight{1.25}\smash{\begin{tabular}[t]{l}$4$\end{tabular}}}}%
    \put(0,0){\includegraphics[width=\unitlength,page=2]{B4a.pdf}}%
  \end{picture}%
\endgroup%
 
		\centering
		\caption{The element $\sigma$ in $B_4$}
         \label{figure-B4element-b}
     \end{subfigure}
     \begin{subfigure}[htb]{0.4\textwidth}
     \centering
\begingroup%
  \makeatletter%
  \providecommand\color[2][]{%
    \errmessage{(Inkscape) Color is used for the text in Inkscape, but the package 'color.sty' is not loaded}%
    \renewcommand\color[2][]{}%
  }%
  \providecommand\transparent[1]{%
    \errmessage{(Inkscape) Transparency is used (non-zero) for the text in Inkscape, but the package 'transparent.sty' is not loaded}%
    \renewcommand\transparent[1]{}%
  }%
  \providecommand\rotatebox[2]{#2}%
  \newcommand*\fsize{\dimexpr\f@size pt\relax}%
  \newcommand*\lineheight[1]{\fontsize{\fsize}{#1\fsize}\selectfont}%
  \ifx\svgwidth\undefined%
    \setlength{\unitlength}{80.83056619bp}%
    \ifx\svgscale\undefined%
      \relax%
    \else%
      \setlength{\unitlength}{\unitlength * \real{\svgscale}}%
    \fi%
  \else%
    \setlength{\unitlength}{\svgwidth}%
  \fi%
  \global\let\svgwidth\undefined%
  \global\let\svgscale\undefined%
  \makeatother%
  \begin{picture}(1,0.97083213)%
    \lineheight{1}%
    \setlength\tabcolsep{0pt}%
    \put(-0.00410774,0.93323423){\makebox(0,0)[lt]{\lineheight{1.25}\smash{\begin{tabular}[t]{l}$1$\end{tabular}}}}%
    \put(0.2742523,0.93323423){\makebox(0,0)[lt]{\lineheight{1.25}\smash{\begin{tabular}[t]{l}$2$\end{tabular}}}}%
    \put(0.55261235,0.93323423){\makebox(0,0)[lt]{\lineheight{1.25}\smash{\begin{tabular}[t]{l}$3$\end{tabular}}}}%
    \put(0.83097239,0.93323423){\makebox(0,0)[lt]{\lineheight{1.25}\smash{\begin{tabular}[t]{l}$4$\end{tabular}}}}%
    \put(0,0){\includegraphics[width=\unitlength,page=1]{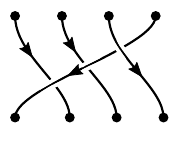}}%
    \put(0.04577095,0.00727312){\makebox(0,0)[lt]{\lineheight{1.25}\smash{\begin{tabular}[t]{l}$1$\end{tabular}}}}%
    \put(0.32413073,0.00727312){\makebox(0,0)[lt]{\lineheight{1.25}\smash{\begin{tabular}[t]{l}$2$\end{tabular}}}}%
    \put(0.60249077,0.00727312){\makebox(0,0)[lt]{\lineheight{1.25}\smash{\begin{tabular}[t]{l}$3$\end{tabular}}}}%
    \put(0.88085349,0.00727312){\makebox(0,0)[lt]{\lineheight{1.25}\smash{\begin{tabular}[t]{l}$4$\end{tabular}}}}%
  \end{picture}%
\endgroup%

 		\caption{The element $\tau$ in $B_4$} 
         \label{figure-B4element-a}
     \end{subfigure}
\caption{Two examples of braids on $4$ strands}
\label{figure-B4}
\end{figure}

\begin{figure}
     \centering
     	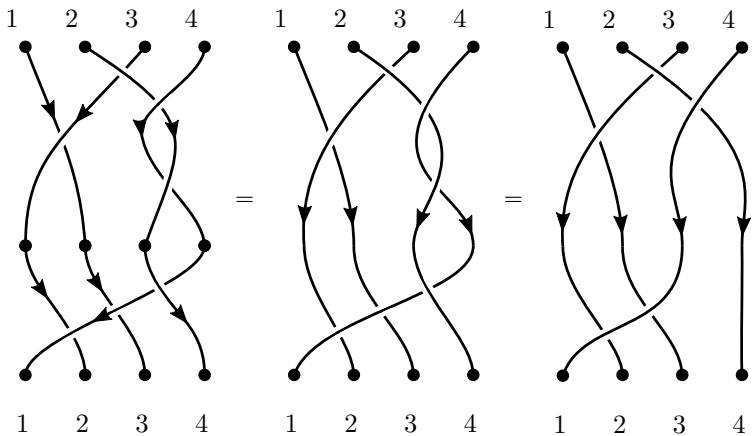
\caption{The product $\sigma\tau$ of the elements $\sigma$ and $\tau$ in the braid group $B_4$}
\label{figure-B4-ab}
\end{figure}

\begin{figure}
     \centering
     \begin{subfigure}[htb]{0.4\textwidth}
     \centering
\begingroup%
  \makeatletter%
  \providecommand\color[2][]{%
    \errmessage{(Inkscape) Color is used for the text in Inkscape, but the package 'color.sty' is not loaded}%
    \renewcommand\color[2][]{}%
  }%
  \providecommand\transparent[1]{%
    \errmessage{(Inkscape) Transparency is used (non-zero) for the text in Inkscape, but the package 'transparent.sty' is not loaded}%
    \renewcommand\transparent[1]{}%
  }%
  \providecommand\rotatebox[2]{#2}%
  \newcommand*\fsize{\dimexpr\f@size pt\relax}%
  \newcommand*\lineheight[1]{\fontsize{\fsize}{#1\fsize}\selectfont}%
  \ifx\svgwidth\undefined%
    \setlength{\unitlength}{82.41793305bp}%
    \ifx\svgscale\undefined%
      \relax%
    \else%
      \setlength{\unitlength}{\unitlength * \real{\svgscale}}%
    \fi%
  \else%
    \setlength{\unitlength}{\svgwidth}%
  \fi%
  \global\let\svgwidth\undefined%
  \global\let\svgscale\undefined%
  \makeatother%
  \begin{picture}(1,0.76931514)%
    \lineheight{1}%
    \setlength\tabcolsep{0pt}%
    \put(0,0){\includegraphics[width=\unitlength,page=1]{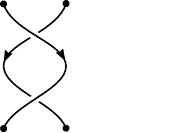}}%
    \put(0.47573434,0.36609309){\makebox(0,0)[lt]{\lineheight{1.25}\smash{\begin{tabular}[t]{l}$=$\end{tabular}}}}%
    \put(0,0){\includegraphics[width=\unitlength,page=2]{braidrel1.pdf}}%
  \end{picture}%
\endgroup%
 
		\centering
		\caption{}
         \label{figure-braidrelation1}
     \end{subfigure}
     \begin{subfigure}[htb]{0.4\textwidth}
     \centering
\begingroup%
  \makeatletter%
  \providecommand\color[2][]{%
    \errmessage{(Inkscape) Color is used for the text in Inkscape, but the package 'color.sty' is not loaded}%
    \renewcommand\color[2][]{}%
  }%
  \providecommand\transparent[1]{%
    \errmessage{(Inkscape) Transparency is used (non-zero) for the text in Inkscape, but the package 'transparent.sty' is not loaded}%
    \renewcommand\transparent[1]{}%
  }%
  \providecommand\rotatebox[2]{#2}%
  \newcommand*\fsize{\dimexpr\f@size pt\relax}%
  \newcommand*\lineheight[1]{\fontsize{\fsize}{#1\fsize}\selectfont}%
  \ifx\svgwidth\undefined%
    \setlength{\unitlength}{153.65887974bp}%
    \ifx\svgscale\undefined%
      \relax%
    \else%
      \setlength{\unitlength}{\unitlength * \real{\svgscale}}%
    \fi%
  \else%
    \setlength{\unitlength}{\svgwidth}%
  \fi%
  \global\let\svgwidth\undefined%
  \global\let\svgscale\undefined%
  \makeatother%
  \begin{picture}(1,0.60785837)%
    \lineheight{1}%
    \setlength\tabcolsep{0pt}%
    \put(0,0){\includegraphics[width=\unitlength,page=1]{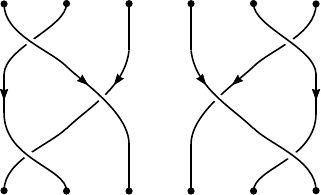}}%
    \put(0.45193899,0.27932067){\makebox(0,0)[lt]{\lineheight{1.25}\smash{\begin{tabular}[t]{l}$=$\end{tabular}}}}%
  \end{picture}%
\endgroup%

 		\caption{} 
         \label{figure-braidrelation2}
     \end{subfigure}

     \begin{subfigure}[htb]{0.8\textwidth}
     \centering
\begingroup%
  \makeatletter%
  \providecommand\color[2][]{%
    \errmessage{(Inkscape) Color is used for the text in Inkscape, but the package 'color.sty' is not loaded}%
    \renewcommand\color[2][]{}%
  }%
  \providecommand\transparent[1]{%
    \errmessage{(Inkscape) Transparency is used (non-zero) for the text in Inkscape, but the package 'transparent.sty' is not loaded}%
    \renewcommand\transparent[1]{}%
  }%
  \providecommand\rotatebox[2]{#2}%
  \newcommand*\fsize{\dimexpr\f@size pt\relax}%
  \newcommand*\lineheight[1]{\fontsize{\fsize}{#1\fsize}\selectfont}%
  \ifx\svgwidth\undefined%
    \setlength{\unitlength}{221.40333129bp}%
    \ifx\svgscale\undefined%
      \relax%
    \else%
      \setlength{\unitlength}{\unitlength * \real{\svgscale}}%
    \fi%
  \else%
    \setlength{\unitlength}{\svgwidth}%
  \fi%
  \global\let\svgwidth\undefined%
  \global\let\svgscale\undefined%
  \makeatother%
  \begin{picture}(1,0.28722747)%
    \lineheight{1}%
    \setlength\tabcolsep{0pt}%
    \put(0,0){\includegraphics[width=\unitlength,page=1]{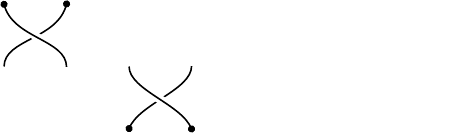}}%
    \put(0.47964152,0.13627501){\makebox(0,0)[lt]{\lineheight{1.25}\smash{\begin{tabular}[t]{l}$=$\end{tabular}}}}%
    \put(0,0){\includegraphics[width=\unitlength,page=2]{braidrel3.pdf}}%
  \end{picture}%
\endgroup%

 		\caption{} 
         \label{figure-braidrelation3}
     \end{subfigure}
\caption{The three fundamental relations among braids, see Equations \eqref{eq-braidrel1}--\eqref{eq-braidrel3}}
\label{figure-braidrelations}
\end{figure}

The product of two braids is computed similarly to the graphical computation of the product in the symmetric groups. The difference is that we need to be careful about the order of strands at each crossing. Take, for example, the two braids $\sigma$ and $\tau$ depicted in Figure~\ref{figure-B4}.
To compute the product $\sigma \tau$, we vertically stack the two braid pictures on top of one another in the plane, $\sigma$ on top of $\tau$. This is shown in the leftmost picture in Figure~\ref{figure-B4-ab}. Then we remove the dots in the middle in the second picture from the left in Figure~\ref{figure-B4-ab}. Finally, we simplify the picture, if possible. When simplifying, we are not allowed to cut the strings or change the positions of the endpoints. The rightmost picture in Figure~\ref{figure-B4-ab} shows a simplification of the product $\sigma \tau$.
To simplify $\sigma \tau$, we used the fundamental relation that crossing two strings with the same string on top twice is the same as two parallel strings. 
This relation is displayed in Figure~\ref{figure-braidrelations}{\scshape (a)}.

\begin{figure}
     \centering
     \begin{subfigure}[htb]{0.48\textwidth}
     \centering
$$ \sigma_i \,= \,\vcenter{\hbox{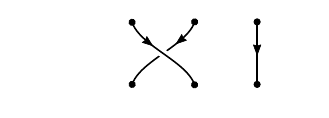}}$$
		\caption{}
         \label{figure-sigmai}
     \end{subfigure}
     \centering
     \begin{subfigure}[htb]{0.48\textwidth}
     \centering
$$ \sigma_i^{-1} \,= \,\vcenter{\hbox{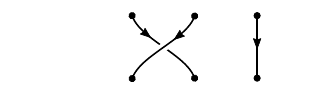}}$$
		\caption{}
         \label{figure-sigmaiinv}
     \end{subfigure}
\caption{The $i$-th generator $\sigma_i$ of the braid group $B_n$ and its inverse $\sigma_i^{-1}$}
\label{figure-sigmais}
\end{figure}

Pictures of braids as used in this text are called \emph{planar projections} because they project a three-dimensional object onto a two-dimensional plane. Such pictures are a way to illustrate manipulations of braids and to give an intuitive understanding. However, sometimes it is useful to have a more compact methods to denote elements in the braid groups. To do this, we introduce abbreviations for some fundamental braids that only braid two neighboring strings and nothing else. These fundamental braids are called \emph{generators} of the braid groups. Remember that the braid group $B_n$ consists of pictures with $n$ strands. We can choose $i$ to be any number from $1$ to $n-1$ and abbreviate (or denote) the braid that only braids the $i$-th strand \emph{over} its right neighbor (the $(i+1)$-th strand) by $\sigma_i$. In addition, we have its inverse $\sigma_i^{-1}$ which braids the $i$-th \emph{under} the $(i+1)$-th strand. Any braid picture can be written by stacking braids $\sigma_i$ and $\sigma_i^{-1}$, where $i$ ranges between $1$ and $n-1$. Figure~\ref{figure-sigmais} contains a picture of the generator $\sigma_i$ its inverse $\sigma_i^{-1}$.
 
Every braid in $B_4$ can be written as a product of (possibly multiple copies of) some of the braids $\sigma_1,\ldots, \sigma_{n-1}$ and $\sigma_1^{-1}, \ldots, \sigma_{n-1}^{-1}$. To demonstrate this idea, recall the two braids $\sigma$ and $\tau$ from Figure~\ref{figure-B4}. We find that 
\begin{equation}
\sigma =\sigma_2\sigma_1^{-1}\sigma_3^{-1}\sigma_3^{-1} \qquad \text{and }\qquad \tau = \sigma_3\sigma_2^{-1}\sigma_1^{-1}.
\end{equation}
An expression of a braid in terms of the generators $\sigma_i$ is by no means unique. For example, we can write $\sigma$ as such a product in multiple ways, for example,
\begin{equation}\label{eq-sigma-versions}
\sigma = \sigma_2\sigma_1^{-1}\sigma_3^{-1}\sigma_3^{-1}=\sigma_2\sigma_3^{-1}\sigma_1^{-1}\sigma_3^{-1}=\sigma_2\sigma_3^{-1}\sigma_3^{-1}\sigma_1^{-1}
\end{equation}
are all valid ways to write the braid group element $\sigma$. They correspond to slightly different pictures but all capture the essence of the same braid. 

There are, in fact, three fundamental types of relations that can be used to relate any two pictures (or, combinations of the elements $\sigma_i, \sigma_i^{-1}$) that display the same element of the braid group. These relations are found in Figure~\ref{figure-braidrelations}. 
The relations can be translated into the following formulas:
\begin{align}\label{eq-braidrel1}
\sigma_i \sigma_i^{-1}&= 1_n,&&\text{for all $i=1,\ldots, n-1$},& \text{Figure~\ref{figure-braidrelations}{\scshape (a)}},\\\label{eq-braidrel2}
\sigma_i \sigma_{i+1}\sigma_i&= \sigma_{i+1}\sigma_i\sigma_{i+1},&&\text{for all $i=1,\ldots, n-2$},& \text{Figure~\ref{figure-braidrelations}{\scshape (b)}},\\
\sigma_i \sigma_j&= \sigma_j\sigma_i,&&\text{for all $j<i-1$ or $j>i+1$,}& \text{Figure~\ref{figure-braidrelations}{\scshape (c)}},\label{eq-braidrel3}
\end{align}
where in the third relation \eqref{eq-braidrel3} we require that $|j-i|>1$, so that $i$ and $j$ are neither neigboring indices (which would be the case if $|j-i|=1$) nor equal (meaning $|j-i|=0$). In Equation \eqref{eq-braidrel1} the symbol $1_n$ denotes the identity element of $B_n$ which only consists of $n$ unbraided strands. Equation \eqref{eq-braidrel2} is an especially famous equation. It goes by the name of the \emph{3rd Reidemeister move} (see Figure \ref{figure-reidemeister}{\scshape (c)}) in knot theory and is called the \emph{Yang--Baxter equation}\footnote{This equation appeared first as a consistency equation in a multi-body problem on a line in quantum mechanics, and in statistical mechanics, see \cite{Jim}.} in physics. For the study of braids, it is astonishing that that successive application of these three types of relations is all that is required to relate any two pictures that represent the same element of the braid group. For example, we have seen how Equation \eqref{eq-braidrel1} was applied in the computation in Figure~\ref{figure-B4-ab} and applying Equation \eqref{eq-braidrel3} is used to relate the different ways to write $\sigma$ in Equation \eqref{eq-sigma-versions}.

\subsection{Passing from braids to permutations}\label{section-HomBtoS}

The braid groups can be seen as an enhancement of the symmetric groups that we discussed before. Given a braid, we can start with one of the top vertices of a braid picture and trace it along the connected string, in the direction of the arrow, to the bottom to see which number we end up with. Doing this for all top vertices, we associate an element of the symmetric group $S_n$ (a \emph{permutation}) to an element of the braid group $B_n$. 

For example, recall the braid $\sigma$ from Figure~\ref{figure-B4}{\scshape (a)}, namely
$$\sigma = \vcenter{\hbox{}}.
$$ We can start at the top vertex $1$ and, tracing its arrow through the picture, we end up at $2$. Similarly, starting at $2$, we end up at $3$, etc. Therefore, the braid $\sigma$ corresponds to the permutation $(2,3,1,4)$ of $S_4$ because there are four vertices at the top of the braid.

This procedure of turning braids into permutations loses information at each crossing, namely, the information of which strand crosses above and which one crosses below is not retained in the symmetric groups. For this reason, there are elements in the braid group that are non-equal but correspond to the same element in the symmetric group. For instance, the elements $\sigma_1$ and $\sigma_1^{-1}$ in $B_2$ are \emph{not} equal as we need to keep track of which strand crosses above and which crosses below. However, both of these elements correspond to the element $(2,1)$ of the symmetric group $S_2$, which swaps $1$ and $2$. The braid groups $B_n$ each contain infinitely many elements while there are only finitely many elements in $S_n$. Any permutation can be obtained from infinitely many distinct braids.

Passing from the braid group $B_n$ to the symmetric group $S_n$ is compatible with the corresponding products of the elements. This means that it does not matter whether we first multiply two elements $a,b$ in $B_n$ and translate the product $ab$ to $S_n$, or first translate both $a$ and $b$ to $S_n$ and then compute the product in the group $S_n$. Mathematicians call this kind of map a \emph{group homomorphism}. We will make use of this observation when analyzing the patterns of the blanket later. We can summarize this observation in formulas. Let $\sigma,\tau$ denote braids, that is, two elements of $B_n$. We write $\rP(\sigma)$ and $\rP(\tau)$ for the resulting permutations obtained from these elements. This way, $\rP$ defines a \emph{function} with inputs from $B_n$ and outputs in $S_n$. For this function, the equation 
\begin{equation}\label{eq-group-hom}
\rP(\sigma)\rP(\tau)=\rP(\sigma\tau)
\end{equation}
holds true as an equality of permutation---both sides are exactly the same reorderings of the set $\{1,2,\ldots, n\}$. You can check this property using the product computation in Figure \ref{figure-B4-ab}. First, we check that 
\begin{align*}
\rP(\sigma)=(2,3,1,4), \qquad \rP(\tau)=(2,3,4,1).
\end{align*}
We can now compute the product two ways, once in $B_4$ and once in $S_4$, and see that Equation \eqref{eq-group-hom} holds in this example. First, we compute the product of permutations as in Section \ref{section-Sn} which yields
\begin{align*}
\rP(\sigma)\rP(\tau)=(2,3,1,4)(2,3,4,1)=(3,4,2,1).
\end{align*}
Second, we read off the permutation associated to the product $\sigma\tau$ from Figure~\ref{figure-B4-ab}, which gives
\begin{align*}
\rP(\sigma\tau)=(3,4,2,1),
\end{align*}
the same result.

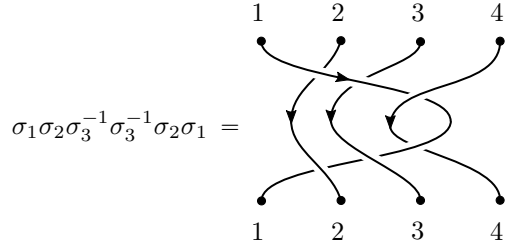
\begin{figure}[htb]
   $$\sigma_1\sigma_2\sigma_3^{-1}\sigma_3^{-1}\sigma_2\sigma_1\,=\,\vcenter{\hbox{
\begingroup%
  \makeatletter%
  \providecommand\color[2][]{%
    \errmessage{(Inkscape) Color is used for the text in Inkscape, but the package 'color.sty' is not loaded}%
    \renewcommand\color[2][]{}%
  }%
  \providecommand\transparent[1]{%
    \errmessage{(Inkscape) Transparency is used (non-zero) for the text in Inkscape, but the package 'transparent.sty' is not loaded}%
    \renewcommand\transparent[1]{}%
  }%
  \providecommand\rotatebox[2]{#2}%
  \newcommand*\fsize{\dimexpr\f@size pt\relax}%
  \newcommand*\lineheight[1]{\fontsize{\fsize}{#1\fsize}\selectfont}%
  \ifx\svgwidth\undefined%
    \setlength{\unitlength}{97.16830038bp}%
    \ifx\svgscale\undefined%
      \relax%
    \else%
      \setlength{\unitlength}{\unitlength * \real{\svgscale}}%
    \fi%
  \else%
    \setlength{\unitlength}{\svgwidth}%
  \fi%
  \global\let\svgwidth\undefined%
  \global\let\svgscale\undefined%
  \makeatother%
  \begin{picture}(1,0.88689511)%
    \lineheight{1}%
    \setlength\tabcolsep{0pt}%
    \put(0,0){\includegraphics[width=\unitlength,page=1]{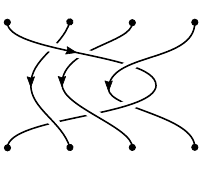}}%
    \put(-0.00243621,0.85366292){\makebox(0,0)[lt]{\lineheight{1.25}\smash{\begin{tabular}[t]{l}$1$\end{tabular}}}}%
    \put(0.30826246,0.85366292){\makebox(0,0)[lt]{\lineheight{1.25}\smash{\begin{tabular}[t]{l}$2$\end{tabular}}}}%
    \put(0.61805795,0.85561886){\makebox(0,0)[lt]{\lineheight{1.25}\smash{\begin{tabular}[t]{l}$3$\end{tabular}}}}%
    \put(0.92484453,0.85418907){\makebox(0,0)[lt]{\lineheight{1.25}\smash{\begin{tabular}[t]{l}$4$\end{tabular}}}}%
    \put(-0.00341707,0.00605023){\makebox(0,0)[lt]{\lineheight{1.25}\smash{\begin{tabular}[t]{l}$1$\end{tabular}}}}%
    \put(0.3072816,0.00605023){\makebox(0,0)[lt]{\lineheight{1.25}\smash{\begin{tabular}[t]{l}$2$\end{tabular}}}}%
    \put(0.61707709,0.00800616){\makebox(0,0)[lt]{\lineheight{1.25}\smash{\begin{tabular}[t]{l}$3$\end{tabular}}}}%
    \put(0.9238639,0.00657638){\makebox(0,0)[lt]{\lineheight{1.25}\smash{\begin{tabular}[t]{l}$4$\end{tabular}}}}%
  \end{picture}%
\endgroup%
}}$$
\caption{An example of a pure braid in $B_4$}
\label{figure-pure-braid}
\end{figure}

There are even braids $\beta$ whose associated permutation $\rP(\beta)$ is just the identity permutation $(1,2,\ldots, n)$ which does not permute any of the elements $1,2,\ldots, n$. Such braids are called \emph{pure braids}. Figure~\ref{figure-pure-braid} shows an example that illustrates that pure braids can already be rather complicated, although the associated permutation is trivial. To see that this braid is a pure braid, start with any input vertex and trace it through the picture. You will see that you end up at the same number that you started from.

The set of all pure braids contained in $B_4$ is denoted by $PB_4$. The subset $PB_4$ is closed under taking inverses and products and contains the identity element of $B_4$. Therefore, $PB_4$ is a \emph{subgroup} of $B_n$, called the \emph{pure braid group on $n$ strands}.

\begin{figure}
\centering 
$$
\vcenter{\hbox{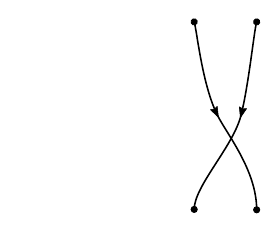}}\,=\,
\vcenter{\hbox{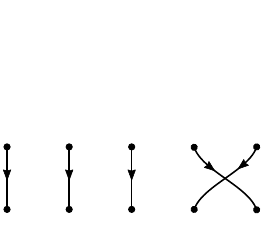}}\,=\,\tau_2\tau_1\tau_4
$$
\caption{The decomposition of the element from Equation~\eqref{eq-prod-trans} in $S_5$ as a product of transpositions}
\label{figure-prod-trans}
\end{figure}

In Equations \eqref{eq-braidrel1}--\eqref{eq-braidrel3}, we explained fundamental relations that are enough to transform any two different pictures of the same braid into one another. Because the symmetric group $S_n$ of permutations is closely related to the braid group $B_n$, the generators and relations for $S_n$ are very similar, but easier to work with in practice. We write $\tau_i$ for the elementary permutation that interchanges only the $i$-th and $(i+1)$-th element in the set $\{1,2,\ldots, n\}$. In our notation, this means 
\begin{align}\label{eq-taui}
\tau_i=(1,2,\ldots, i-1, i+1, i, i+2, \ldots, n).
\end{align}
For instance, $\tau_2=(1,3,2,4)$ as an element of $S_4$.
\begin{figure}[htb]
$$ \tau_i \,= \,\vcenter{\hbox{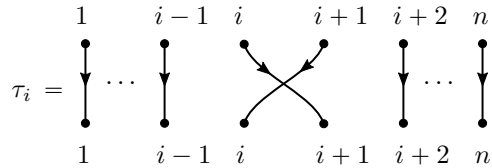}}$$
\caption{The $i$-th generator $\tau_i$ of $S_n$ which transposes $i$ and $i+1$.}
\label{figure-Sn-generator}
\end{figure}
The generator $\tau_i$ is depicted in Figure \ref{figure-Sn-generator} using the same graphical depiction of permutation we used before.
The notation is understood in the way that if $i-1$, $i+2$ do not exist, they are simply omitted. For example, $\tau_1=(2,1)$ is the only generator for the symmetric group $S_2$, while $S_3$ has the two generators
\begin{align}
\tau_1=(2,1,3), \qquad \tau_2=(1,3,2).
\end{align}
The element $\tau_5$ exists in all symmetric groups $S_n$ with $n\geq 6$. For example,
$$\tau_5=(1,2,3,4,6,5,7,8,9,10)$$
when viewed as a generator of $S_{10}$. The elements $\tau_i$ are called \emph{transpositions} or \emph{elementary transpositions}. Any element of $S_n$ can be written as an iterated product of these transpositions. For example
\begin{align}\label{eq-prod-trans}
(2,3,1,5,4)=(1,3,2,4,5)(2,1,3,4,5)(1,2,3,5,4)=\tau_2\tau_1\tau_4.
\end{align}
Here, we recall that products are read from left to right. So in the product in the middle, first, $3$ is sent to $2$, then $2$ is sent to $1$, and finally $1$ is sent to $1$.
Hence, tracing through the entire string of mappings, $3$ is sent to $1$. We illustrate this example in Figure \ref{figure-prod-trans}

The relation describing the symmetric group based on these generators are
\begin{align}\label{eq-permrel1}
\tau_i \tau_i&= 1_n,&&\text{for all $i$,}\\\label{eq-permrel2}
\tau_i \tau_{i+1}\tau_i&= \tau_{i+1}\tau_i\tau_{i+1},&&\text{for all $i<n-1$,}\\
\tau_i \tau_j&= \tau_j\tau_i, && \text{whenever $j\neq i$, $i-1$, and $i+1$}.\label{eq-permrel3}
\end{align}
These relations are almost the same as for the braid group. The only difference is that $\tau_i^{-1}=\tau_i$, which may be derived from Equation \eqref{eq-permrel1}. The geometric explanation for this relation is that the order in which the strands cross does not matter for permutation. One can observe that the generators $\tau_i$ for $S_n$ satisfy all of the relations that the $\sigma_i$ satisfy in $B_n$. This explains why the map $\rP\colon B_n\to S_n$ described earlier is a homomorphism, that is, it is compatible with products. In fact, it follows from the definitions that
$\rP(\sigma_i)=\tau_i$.

\subsection{Orders of elements in a group}\label{section-orders}

We will discuss one more concept from the theory of groups before starting to look at the braid patterns of the blanket more closely. This is the concept of the \emph{order} of an element. Let $g$ be an element inside of a group (for example, the symmetric groups, but the concept of order is defined for any group). If there exists a positive integer number $n$ such that $g^n=1$ and $n$ is the smallest positive number with this property, then $g$ is said to have \emph{order $n$}. If no such number $n$ exists, then $g$ has \emph{infinite order}. A basic result is that if a group has $m$ elements, then the order of any element $g$ in this group divides $m$.\footnote{This result, named after 18th century mathematician Joseph-Louis Lagrange, can be found in standard textbooks on the theory of groups such as  \cite{Rot}.} This result strongly limits the possible orders of elements in a given group and implies that if the group only have finitely many elements, then there is no element of infinite order. 

The symmetric groups only have finitely many elements. In fact, we can verify that the symmetric group $S_n$ has 
$$n!=n(n-1)\cdot\ldots \cdot 1$$
elements. This number grows extremely quickly. For example,
\begin{align}
4!&=4\cdot 3\cdot 2\cdot 1=24, & 5!&=5\cdot 4!=120,& 6!&=6\cdot 5!=720, &
7!&=7\cdot 6!=5040, &\ldots 
\end{align}
In fact, we can have elements in $S_n$ of any prime order. An example of the order of an element we may study the element $\rP(\sigma)=(2,3,1,4)$ in the symmetric group $S_4$. We already know from the start that the order of this element can be $2,3,4,6,8,12$, or even $24$. To find the order, we compute successive products of $\rP(\sigma)$ with itself. That is, we compute powers of this element. We see that 
\begin{align*}
\rP(\sigma)^2&=(3,1,2,4), &\rP(\sigma)^3&=(1,2,3,4).
\end{align*}
We do not need to compute further since the third power equals the identity. This tells us that $\rP(\sigma)$ has order $3$.\footnote{Group theory offers more effective methods to compute the order of elements in symmetric groups by decomposing such elements into products of \emph{cycles}.  A cycle is a sequence of numbers obtained by repeated application of a permutation that returns to the initial element. For instance, the element $(3,1,2,4)$ contains only one cycle of order $3$, the cycle $1\mapsto 3\mapsto 2\mapsto 1$. The order of an element of a symmetric group is the least common multiple of all cycle lengths of a given element. For more details, see for example \cite{Rot}*{Chapter~2}.
}

There are infinitely many elements in the braid group $B_n$ and even in the pure braid group $PB_n$, for any number $n$ of strands. For example, we can always multiply the element $\sigma_1^2$ with itself and one gets increasingly tangled up strands $\sigma_1^2$, $\sigma_1^4$, $\sigma_1^6, \ldots$ which will never be the same braid. A harder fact to prove about the braid groups is that the identity element is the only element of finite order. The identity always has order $1$. All other element of $B_n$ have infinite order.

\section{The Blanket}

We have now provided all the mathematical definitions to discuss the patterns on the blanket in more detail. We will first look for symmetries and repeating patterns in the braids of the blanket. 

\subsection{The repeating braid patterns of the blanket}
\label{section-reppatterns}

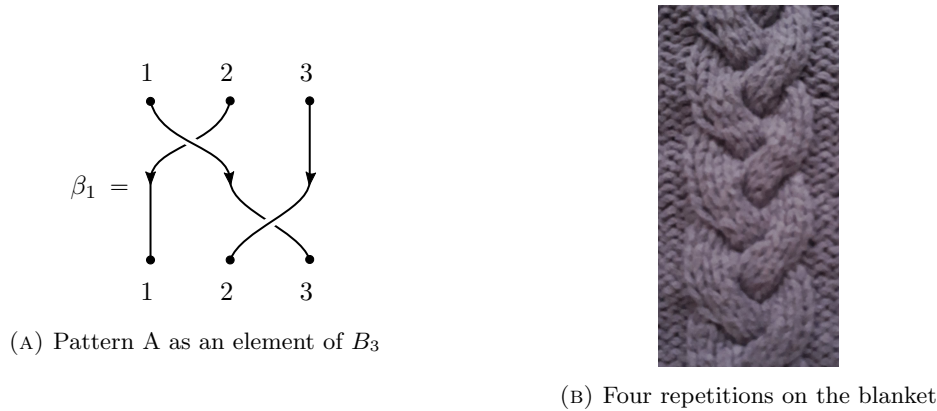
\begin{figure}
\centering
\begin{subfigure}[htb]{0.4\textwidth}
\centering
   $$\beta_1\,=\,\vcenter{\hbox{
\begingroup%
  \makeatletter%
  \providecommand\color[2][]{%
    \errmessage{(Inkscape) Color is used for the text in Inkscape, but the package 'color.sty' is not loaded}%
    \renewcommand\color[2][]{}%
  }%
  \providecommand\transparent[1]{%
    \errmessage{(Inkscape) Transparency is used (non-zero) for the text in Inkscape, but the package 'transparent.sty' is not loaded}%
    \renewcommand\transparent[1]{}%
  }%
  \providecommand\rotatebox[2]{#2}%
  \newcommand*\fsize{\dimexpr\f@size pt\relax}%
  \newcommand*\lineheight[1]{\fontsize{\fsize}{#1\fsize}\selectfont}%
  \ifx\svgwidth\undefined%
    \setlength{\unitlength}{67.00685038bp}%
    \ifx\svgscale\undefined%
      \relax%
    \else%
      \setlength{\unitlength}{\unitlength * \real{\svgscale}}%
    \fi%
  \else%
    \setlength{\unitlength}{\svgwidth}%
  \fi%
  \global\let\svgwidth\undefined%
  \global\let\svgscale\undefined%
  \makeatother%
  \begin{picture}(1,1.28591371)%
    \lineheight{1}%
    \setlength\tabcolsep{0pt}%
    \put(0,0){\includegraphics[width=\unitlength,page=1]{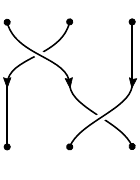}}%
    \put(-0.00441554,1.24055924){\makebox(0,0)[lt]{\lineheight{1.25}\smash{\begin{tabular}[t]{l}$1$\end{tabular}}}}%
    \put(0.44329934,1.24055924){\makebox(0,0)[lt]{\lineheight{1.25}\smash{\begin{tabular}[t]{l}$2$\end{tabular}}}}%
    \put(0.89101519,1.24055924){\makebox(0,0)[lt]{\lineheight{1.25}\smash{\begin{tabular}[t]{l}$3$\end{tabular}}}}%
    \put(-0.00495518,0.00877358){\makebox(0,0)[lt]{\lineheight{1.25}\smash{\begin{tabular}[t]{l}$1$\end{tabular}}}}%
    \put(0.4427597,0.00877358){\makebox(0,0)[lt]{\lineheight{1.25}\smash{\begin{tabular}[t]{l}$2$\end{tabular}}}}%
    \put(0.89047555,0.00877358){\makebox(0,0)[lt]{\lineheight{1.25}\smash{\begin{tabular}[t]{l}$3$\end{tabular}}}}%
  \end{picture}%
\endgroup%
}}$$
\caption{Pattern A as an element of $B_3$}
\label{figure-patternA-B3}
\end{subfigure}
\begin{subfigure}[htb]{0.5\textwidth}
\centering
\includegraphics[scale=0.25]{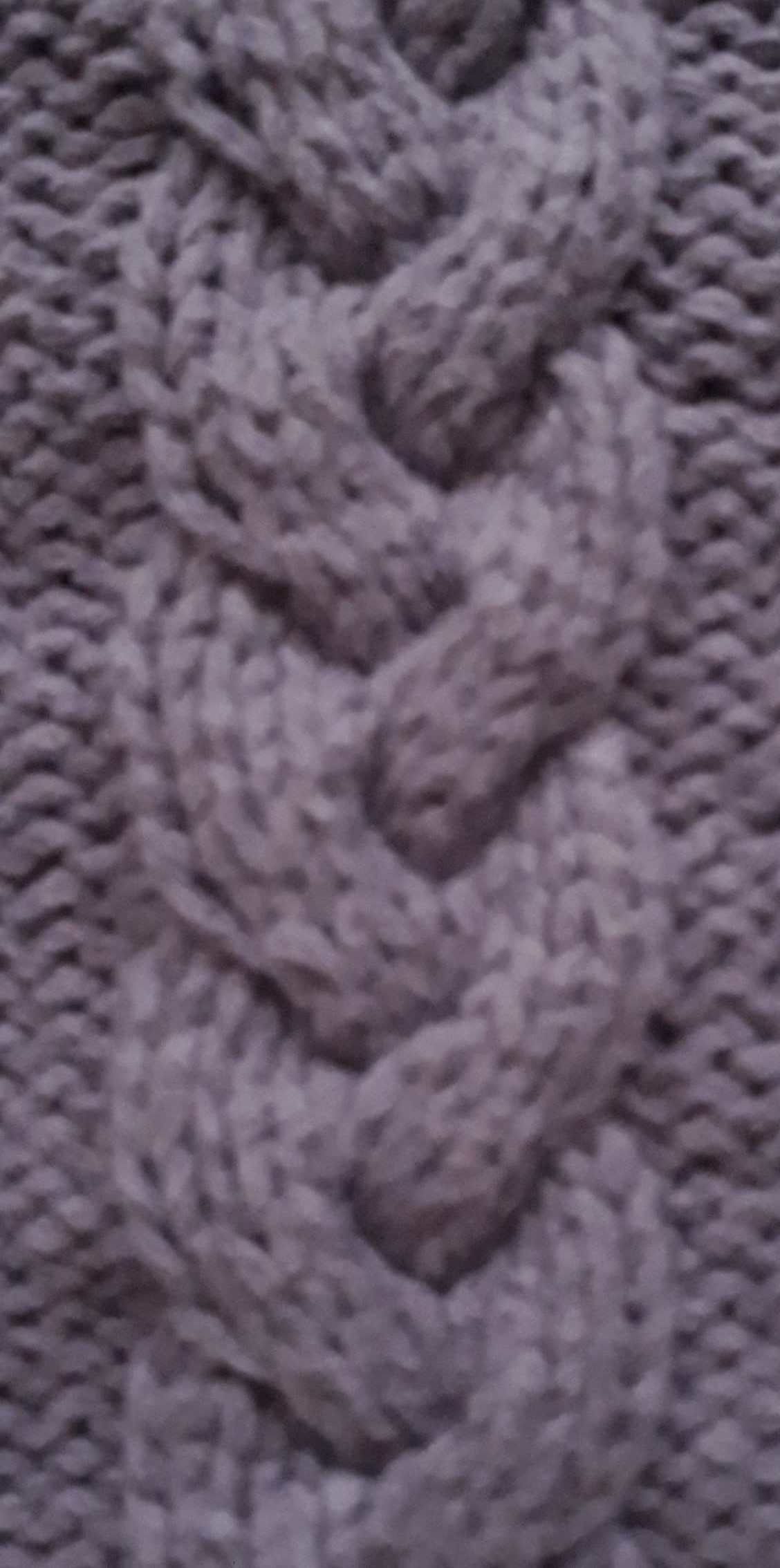}
\caption{Four repetitions on the blanket}
\label{figure-patternA-blanket}
\end{subfigure}
\caption{The repeating Pattern A of the blanket}
\label{figure-patternA}
\end{figure}

There are three distinct patterns that can be observed on the blanket. The first pattern, Pattern A, is the smallest. It involves only three strands and is a very recognizable braid, reminiscent of patterns of braided hair. Pattern A is shown in Figure~\ref{figure-patternA}, once as a picture of four repetitions of the pattern on the blanket and once as a schematic drawing of a single repetition, identified as the element 
\begin{align}\label{pattern1}
\beta_1 = \sigma_1\sigma_2^{-1}
\end{align}
of the braid group $B_3$. 
Therefore, the picture in Figure~\ref{figure-patternA}{\scshape (b)} corresponds to 
\begin{align}
\beta_1^4=\sigma_2^{-1}\sigma_1\sigma_2^{-1}\sigma_1\sigma_2^{-1}\sigma_1\sigma_2^{-1}\sigma_1.
\end{align}

\begin{figure}
\centering
\begin{subfigure}[htb]{0.5\textwidth}
\centering
   $$\beta_2\,=\,\vcenter{\hbox{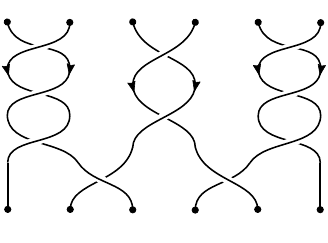}}$$
\caption{Pattern B as an element of $B_6$}
\label{figure-patternB-B6}
\end{subfigure}
\begin{subfigure}[htb]{0.4\textwidth}
\centering
\includegraphics[scale=0.25]{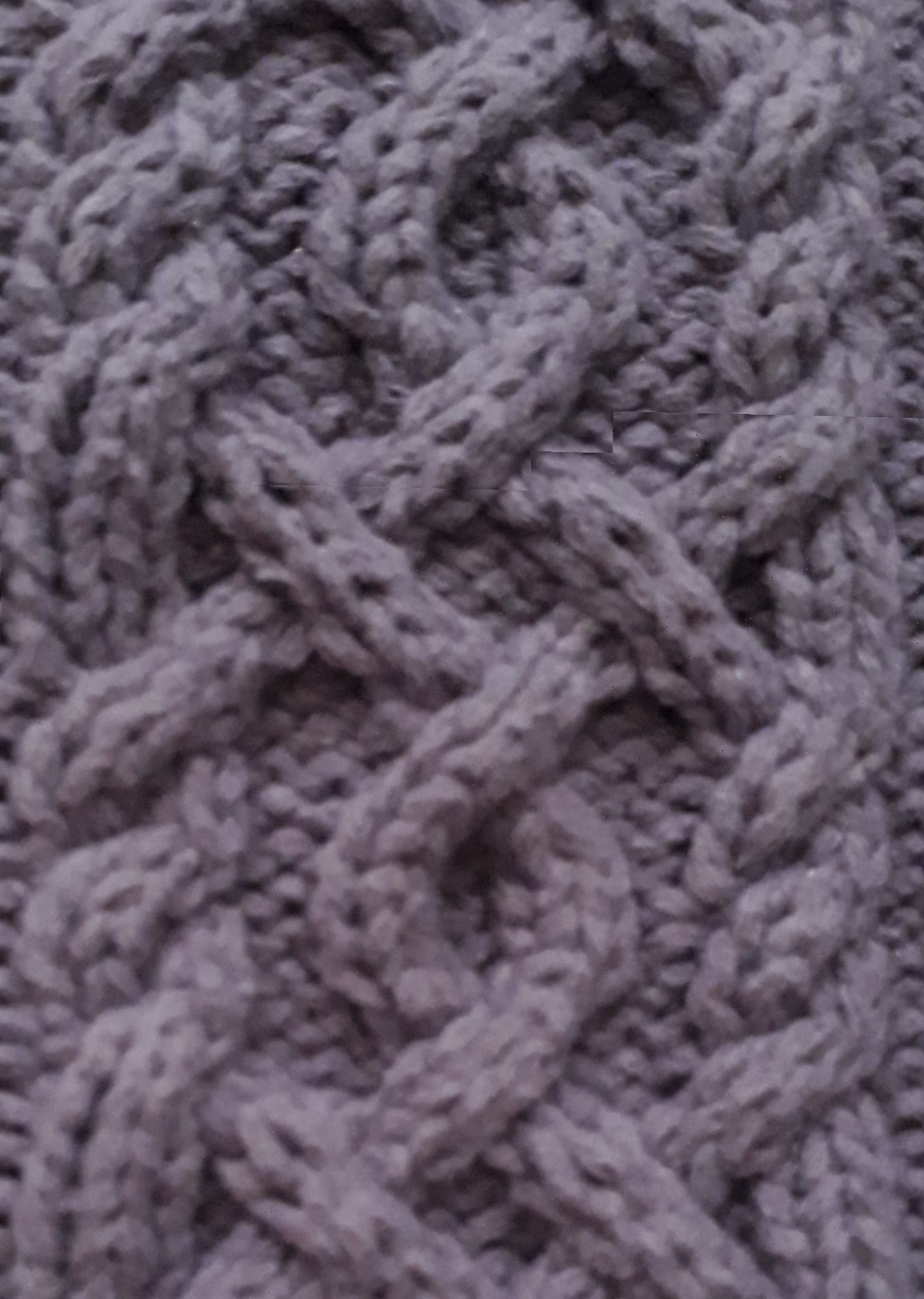}
\caption{Two repetitions on the blanket}
\label{figure-patternB-blanket}
\end{subfigure}
\caption{The repeating Pattern B of the blanket}
\label{figure-patternB}
\end{figure}

The second repeating braid pattern, Pattern B, is depicted in Figure~\ref{figure-patternB}. This pattern corresponds to an element of $B_6$ as 6 braids are used. In Figure~{\scshape\ref{figure-patternB}(a)}, we see a schematic depiction of a single repetition of this pattern. It corresponds to the braid group element 
\begin{align}\label{pattern2}
\beta_2=\sigma_{1}^{-3}\sigma_{3}^{-2}\sigma_{5}^{-3}\sigma_2\sigma_4
\end{align}
of $B_6$. We have summarized products into powers, so that, for example
\begin{align}
\sigma_1^{-3}=\sigma_1^{-1}\sigma_1^{-1}\sigma_1^{-1}.
\end{align}
The right hand side, Figure~{\scshape\ref{figure-patternB}(b)},  shows two repetitions of this pattern on the blanket.

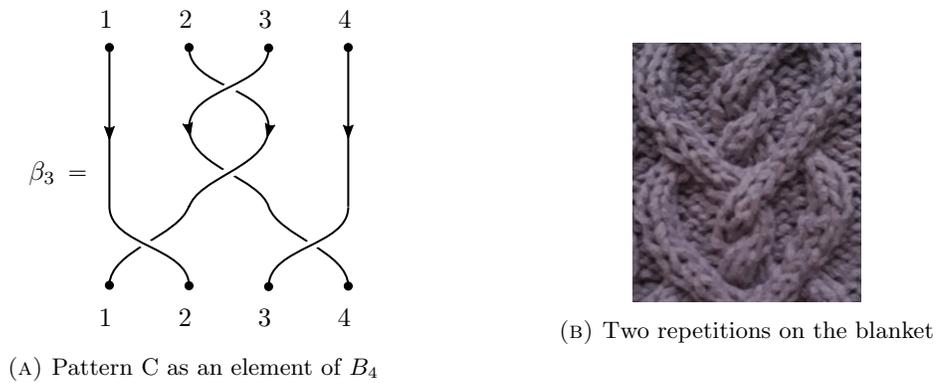
\begin{figure}
\centering
\begin{subfigure}[htb]{0.4\textwidth}
\centering
   $$\beta_3\,=\,\vcenter{\hbox{
\begingroup%
  \makeatletter%
  \providecommand\color[2][]{%
    \errmessage{(Inkscape) Color is used for the text in Inkscape, but the package 'color.sty' is not loaded}%
    \renewcommand\color[2][]{}%
  }%
  \providecommand\transparent[1]{%
    \errmessage{(Inkscape) Transparency is used (non-zero) for the text in Inkscape, but the package 'transparent.sty' is not loaded}%
    \renewcommand\transparent[1]{}%
  }%
  \providecommand\rotatebox[2]{#2}%
  \newcommand*\fsize{\dimexpr\f@size pt\relax}%
  \newcommand*\lineheight[1]{\fontsize{\fsize}{#1\fsize}\selectfont}%
  \ifx\svgwidth\undefined%
    \setlength{\unitlength}{97.23856531bp}%
    \ifx\svgscale\undefined%
      \relax%
    \else%
      \setlength{\unitlength}{\unitlength * \real{\svgscale}}%
    \fi%
  \else%
    \setlength{\unitlength}{\svgwidth}%
  \fi%
  \global\let\svgwidth\undefined%
  \global\let\svgscale\undefined%
  \makeatother%
  \begin{picture}(1,1.19424788)%
    \lineheight{1}%
    \setlength\tabcolsep{0pt}%
    \put(0,0){\includegraphics[width=\unitlength,page=1]{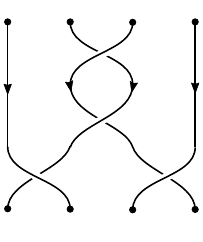}}%
    \put(-0.00065986,1.16299423){\makebox(0,0)[lt]{\lineheight{1.25}\smash{\begin{tabular}[t]{l}$1$\end{tabular}}}}%
    \put(0.30785978,1.16299423){\makebox(0,0)[lt]{\lineheight{1.25}\smash{\begin{tabular}[t]{l}$2$\end{tabular}}}}%
    \put(0.61637898,1.16299423){\makebox(0,0)[lt]{\lineheight{1.25}\smash{\begin{tabular}[t]{l}$3$\end{tabular}}}}%
    \put(0.92489884,1.16299423){\makebox(0,0)[lt]{\lineheight{1.25}\smash{\begin{tabular}[t]{l}$4$\end{tabular}}}}%
    \put(-0.0034146,0.00604585){\makebox(0,0)[lt]{\lineheight{1.25}\smash{\begin{tabular}[t]{l}$1$\end{tabular}}}}%
    \put(0.30510526,0.00604585){\makebox(0,0)[lt]{\lineheight{1.25}\smash{\begin{tabular}[t]{l}$2$\end{tabular}}}}%
    \put(0.61362446,0.00604585){\makebox(0,0)[lt]{\lineheight{1.25}\smash{\begin{tabular}[t]{l}$3$\end{tabular}}}}%
    \put(0.92214432,0.00604585){\makebox(0,0)[lt]{\lineheight{1.25}\smash{\begin{tabular}[t]{l}$4$\end{tabular}}}}%
  \end{picture}%
\endgroup%
}}$$
\caption{Pattern C as an element of $B_4$}
\label{figure-patternC-B4}
\end{subfigure}
\begin{subfigure}[htb]{0.5\textwidth}
\centering
\includegraphics[scale=0.28]{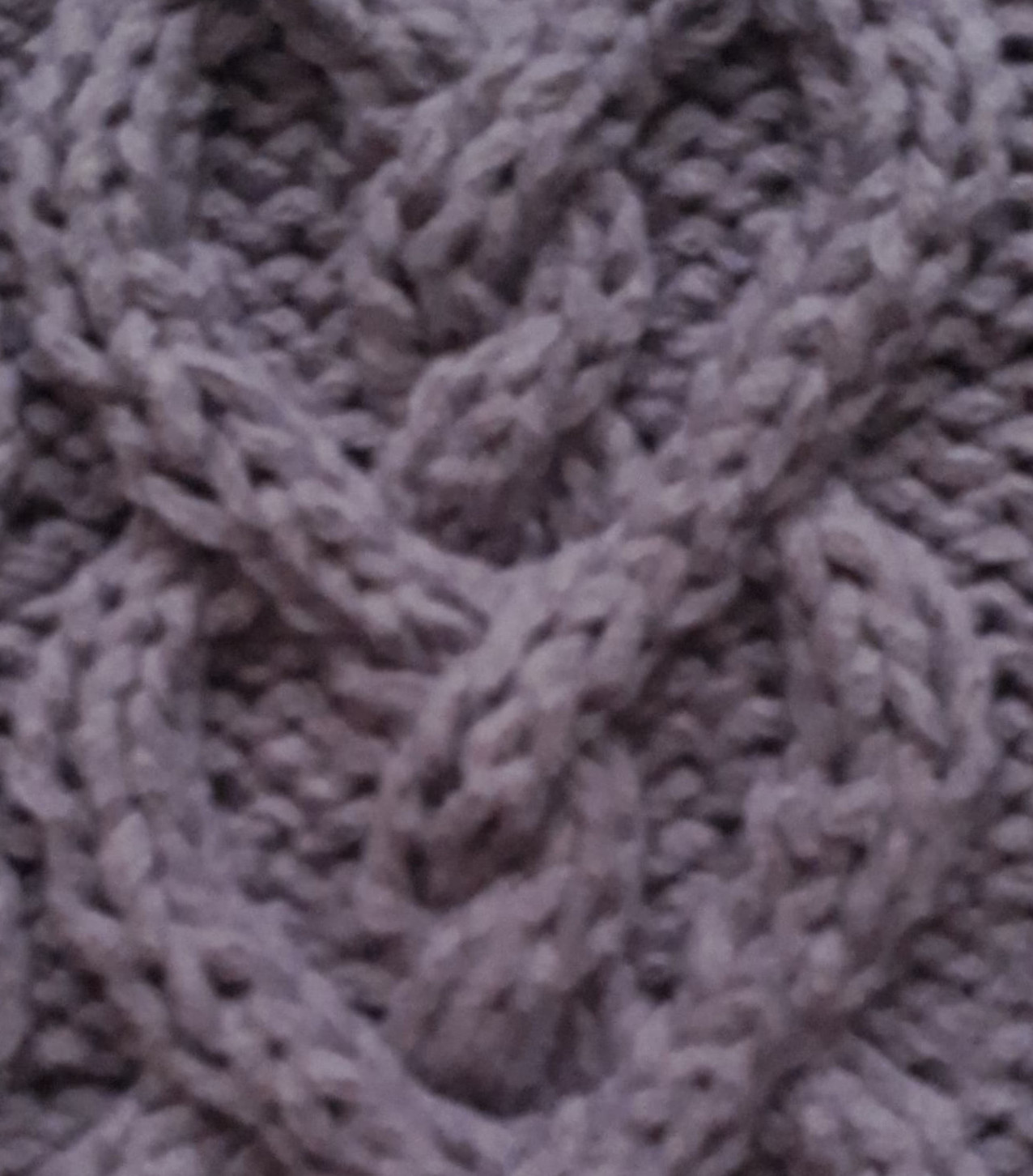}
\caption{Two repetitions on the blanket}
\label{figure-patternC-blanket}
\end{subfigure}
\caption{The repeating Pattern C of the blanket}
\label{figure-patternC}
\end{figure}

There is a third repeating braid pattern found on the blanket. This pattern is displayed in Figure~\ref{figure-patternC}.
Figure~\ref{figure-patternC}{\scshape (a)} schematically displays a single repetition of this pattern as the element 
\begin{align}\label{pattern3}
\beta_3=\sigma_2^{-2}\sigma_1\sigma_3^{-1}
\end{align}
of the braid group $B_4$ while Figure~{\scshape\ref{figure-patternC}(b)} shows two repetitions of this pattern on the blanket.

The three repeating patterns found on the blanket are elements $\beta_1$ of $B_3$, $\beta_2$ of $B_4$, and $\beta_3$ of $B_6$. Using the discussion at the end of Section~\ref{section-Sn}, we find associated permutations $\rP(\beta_1)$ of $3$, $\rP(\beta_2)$ of $4$, respectively, $\rP(\beta_3)$ of $6$ numbers. These permutations trace where the input labels and up traced through the braids. One computes that these permutations are given by 
\begin{align}\label{Pofbetas}
\rP(\beta_1)&=(3,1,2), & \rP(\beta_2)&=(3,1,2,5,6,4), & \rP(\beta_3)&=(2,1,4,3).
\end{align}
In Section~\ref{section-orders}, we discussed the concept of the \emph{order} of an element $g$ as the minimal $n$ such that $g^n=1$. Once checks that $\rP(\beta_1)$ and $\rP(\beta_2)$ both have order $3$ while $\rP(\beta_3)$ has order $2$. 

\subsection{The braid structure of the entire blanket}

We will now study the overall structure of the blanket as a braid. For this we observe the repetition of the basic patterns A, B, and C discussed in the previous section. Looking at Figure~\ref{figure-blanket}, we observe that the basic patterns A, B, and C repeat horizontally using the following scheme:
\begin{align}\label{eq-overallpattern}
\begin{array}{ccccccccc}
A&B&A&C&C&C&A&B&A
\end{array}.
\end{align}
This shows that the repetition pattern is symmetric with respect to $180$ degree rotation about the vertical axis. However, the entire blanket is not rotation symmetric at the vertical middle axis because Patterns A and C are not symmetric with respect to such a rotation. Algebraically, a rotation by $180$ about the vertical axis correspond to exchanging $\sigma_i$ with $\sigma_{n-i}$, for all $i$ between $1$ and $n-1$. The braid $\beta_2$ is symmetric under this operations because
$$\beta_2=\sigma_{1}^{-3}\sigma_{3}^{-2}\sigma_{5}^{-3}\sigma_2\sigma_4= \sigma_{5}^{-3}\sigma_{3}^{-2}\sigma_{1}^{-3}\sigma_4\sigma_2,$$
using the relations from Equation~\eqref{eq-braidrel3}. However, $\beta_1$ is not symmetric under this operation as 
$$\beta_1=\sigma_1\sigma_2^{-1}\neq \sigma_2\sigma_1^{-1}.$$
The braid $\beta_3$ also does not display such a rotation symmetry\footnote{However, using the symmetry could be achieved by using the braid $\sigma_2^{-2}\sigma_1\sigma_3$ instead of $\beta_3$.} since
$$\beta_3=\sigma_2^{-2}\sigma_1\sigma_3^{-1}\neq \sigma_2^{-2}\sigma_3\sigma_1^{-1}.$$

The braids use a total of $36=4\cdot 3+ 2\cdot 6+ 3\cdot 4$ strands, so it is represented by an element of $B_{36}$. An interesting symmetry here is that all types A, B, and C of the repeating patterns use the same number of $12$ strands.

\begin{figure}[p]
 \resizebox{0.2333\width}{0.2333\height}{
 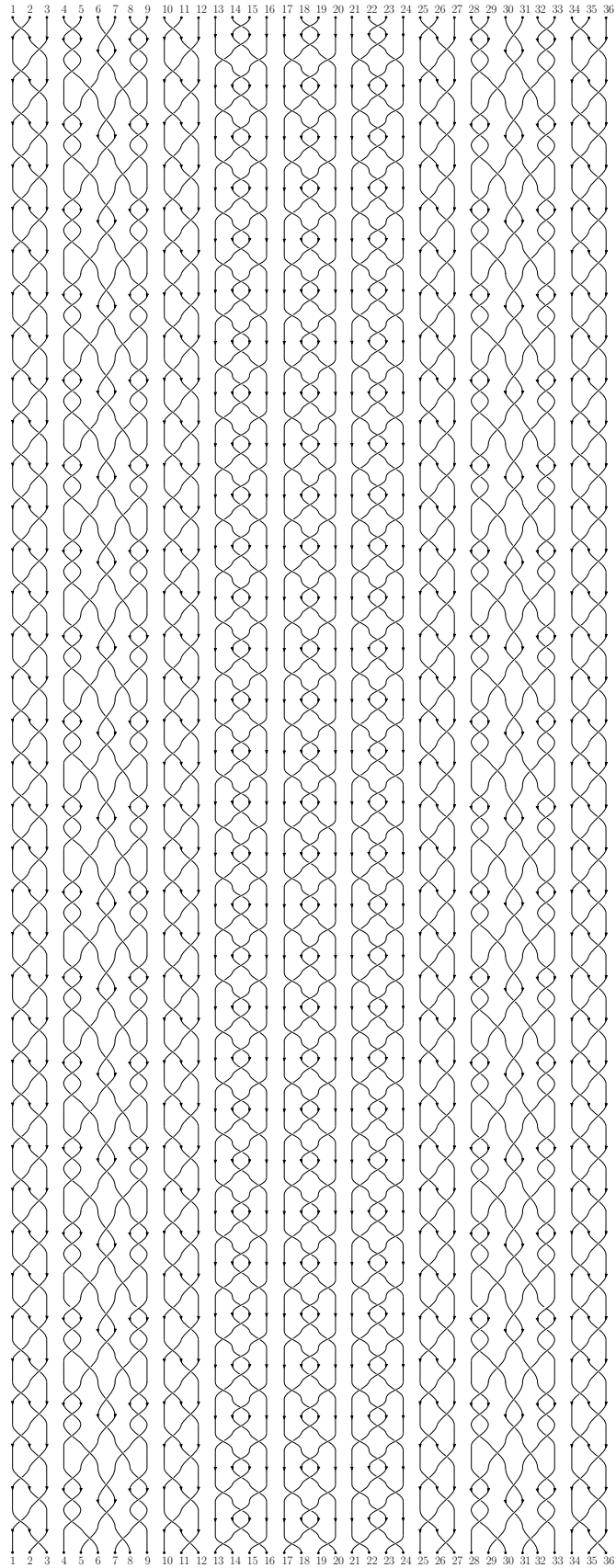}
\caption{The braid $\beta$ represented by the entire blanket}
\label{figure-entire-blanket}
\end{figure}

A careful study of the blanket shows that Pattern A repeats 36 times horizontally, Pattern B repeats 18 times, and Pattern C repeats 30 times. From this we can report how many of the generators $\sigma_i$ or their inverses are used to write the braids on the blanket as an element of $B_{36}$---the number of crossings. The entire braid pattern on the blanket is displayed in Figure~\ref{figure-entire-blanket}. Further data about the number of strands occupied by the pattern and the number of crossings used is summarized in Table \ref{table-stats}. The question whether the number of crossings used in the blanket's design is minimal to capture these braids will be addressed in Section \ref{questions-math}.

\begin{table}[htb]
\begin{tabular}{c|cccccc}
& Strands/copy & Copies & Total strands& Repetitions & Crossings/rep. & Total crossings\\\hline
Pattern A & 3 & 4 & 12 & 36 & 2 & 288\\ 
Pattern B & 6 & 2 & 12 & 18 & 10 & 360\\
Pattern C & 4 & 3 & 12 & 30 & 4 & 360\\\hline
Blanket &&&36&&&1,008
\end{tabular}
\caption{The repeating pattern on the blanket}
\label{table-stats}
\end{table}

We note that if there was one more horizontal repetition of Pattern A, then each Pattern would feature $360$ crossings. The way the blanket is designed, only patterns B and C feature $360$ crossings each, and Pattern A has $288$ crossings. However, $5$ horizontal copies of Pattern A, $3$ of Pattern B, and $4$ of Pattern C cannot be placed in an arrangement that is rotation symmetric about the horizontal middle axis.

We can now identify the overall braid as a product of the braids associated to the repeating patterns, Patterns A, B, and C. For example, Pattern A repeats $36$ times as indicated in Table \ref{table-stats}.  Thus, a copy of repetitions of Pattern A in vertical direction corresponds to $\beta_1^{36}$, the $36$-fold product of $\beta_1$ with itself. 

In the horizontal directions, braids may be placed next to each other. For example, on the left edge of the blanket, we have a copy of Pattern A next to copy of Pattern B on the right. Mathematically, this means that the $36$-fold power of $\beta_1$ is placed on the left of the $18$-fold power of $\beta_1$. We separate horizonally neighboring braids by commas. This, on the left we start with $(\beta_1^{36}, \beta_2^{18})$. We continue following the scheme of pattern repetition observed in Equation \eqref{eq-overallpattern} and find that the blanket corresponds to the braid $\beta$ defined as 
\begin{equation}
\label{eq-blanketbraid}
\beta=(\beta_1^{36}, \beta_2^{18}, \beta_1^{36}, \beta_3^{30},\beta_3^{30},\beta_3^{30},\beta_1^{36}, \beta_2^{18}, \beta_1^{36}).
\end{equation}
A schematic picture of the braid corresponding to the entire blanket is given in Figure \ref{figure-entire-blanket}.

\subsection{Mathematical questions}\label{questions-math}

The blanket is, first and foremost, a piece of art and craft. However, the artwork translates into the abstract structure of a braid, which can be analyzed using a branch of algebra called group theory. Therefore, we are able to use mathematics to extract more precise information from the blanket.  Mathematical questions concerning the braids on this blanket may emerge and may  be answered using the current mathematical understanding of the group theory behind braids. Even if we could not answer these questions using the current understanding of braids by mathematicians, we may be able can post questions that are worth exploring through research. 

\subsubsection*{Pure braids}

One question that we can ask about the blanket is whether it represents a pure braid, as defined in Section~\ref{section-HomBtoS}. 

\begin{question}\label{question-pure}
Does the braid pattern of  the blanket represent a pure braid? 
\end{question}

\begin{proof}[Answer to Question \ref{question-pure}] To answer this question, we recall the definition of a pure braid. A pure braid is a braid $\beta$ (of $n$ strands) such that the associated permutation $\rP(\beta)$ is the identity element of the symmetric group $S_n$. We computed the permutations associated to each of the repeating patterns $\beta_1,\beta_2$, and $\beta_3$ at the end of Section~\ref{section-reppatterns}. We also know how many times each pattern repeats in a vertical direction and in horizontal direction using the information from Table \ref{table-stats}. This way, we identified the braid $\beta$ representing the entire blanket in Equation \ref{eq-blanketbraid}, a braid in $B_{36}$. Thus, if 
$\rP(\beta)=1_{36}$, the identity of $S_{36}$, then we know that $\beta$ is a pure braid. 

To compute $\rP(\beta)$, we use the rule that the function $\rP$ which takes elements of $B_{n}$ (i.e., braids) as inputs and produces elements of $S_{36}$ (i.e. permutations) as outputs is a homomorphism of groups. This implies that 
\begin{align}\label{eq-powersofbraids}
\rP(\beta^n)=\rP(\beta)^n,
\end{align}
see the discussion at the end of Section \ref{section-HomBtoS}. Using this property we conclude that 
\begin{gather}
\rP(\beta_1^{36})=\rP(\beta_1)^{36}=\big(\rP(\beta_1)^3\big)^{12}= 1_{3}^{12}=1_3.
\end{gather}
Here, we use the earlier observation that the order of $\rP(\beta_1)=(3,1,2)$ is $3$ which divides $36$. Hence, the $36$-th power of $\beta_1$ is the identity. Similarly, we calculate that 
\begin{gather}
\rP(\beta_2^{18})=\rP(\beta_2)^{18}=\big(\rP(\beta_2)^3\big)^{6}= 1_{6}^{6}=1_6,\\
\rP(\beta_3^{30})=\rP(\beta_3)^{30}=\big(\rP(\beta_3)^{3}\big)^{10}= 1_{4}^{10}=1_4.
\end{gather}
Thus, the three braids $\beta_1,\beta_2$, and $\beta_3$ appearing as repeating patterns  on the blanket are all pure braids.

We further observe that the function $\rP$ is compatible with putting braids $\beta_1, \beta_2$ next to one another in the sense that 
\begin{equation}\label{eq-parallelbraids}
\rP(\beta_1,\beta_2)=(\rP(\beta_1), \rP(\beta_2)).
\end{equation}
Here, on the right, we represent permutations that permute a disjoint set of strands separated by commas. For example, 
\begin{equation}
\left((3,1,2),(6,4,5,8,9,7)\right)=(3,1,2,6,4,5,8,9,7).
\end{equation}
On the left, we have the permutation $(3,1,2)$ permuting the first three numbers separated by a comma from the permutation $(6,4,5,8,9,7)$ permuting the numbers from $4$ to $9$. The right hand side is a single permutation permuting the first $9$ numbers. 

By Equation \ref{eq-blanketbraid}, the braid $\beta$ representing the entire blanket consists of parallel braids of the form $\beta_1^{36}$, $\beta_2^{18}$, and $\beta_3^{30}$. As we observed above, these are all pure braids. Thus, by the observations of the last paragraph, $\beta$ itself is a pure braid. Indeed, applying the map $\rP$ to $\beta$ gives
\begin{align}
\rP(\beta)& = \rP\big(\beta_1^{36}, \beta_2^{18}, \beta_1^{36}, \beta_3^{30},\beta_3^{30},\beta_3^{30},\beta_1^{36}, \beta_2^{18}, \beta_1^{36}\big)\\
&=\big(\rP(\beta_1^{36}), \rP(\beta_2^{18}), \rP(\beta_1^{36}), \rP(\beta_3^{30}),\rP(\beta_3^{30}),\rP(\beta_3^{30}),\rP(\beta_1^{36}), \rP(\beta_2^{18}), \rP(\beta_1^{36})\big)\\
&=(1_3,1_6,1_3,1_4,1_4,1_4,1_3,1_6,1_3)= 1_{36}.
\end{align}
The last equation uses that parallel identities on subsets of the strands combine into the identity on all strands. Thus, the blanket represents a pure braid, and Question \ref{question-pure} is answered in the affirmative.
\end{proof}

Question~\ref{question-pure} was not a very difficult question. We were able to answer this question using elementary group theory found in textbooks. Alternatively, to confirm that the blanket represents a pure braid, we could have considered Figure~\ref{figure-blanket} and traced every one of the $36$ strands from top to bottom to check that they each start and end at the same number. 

\subsubsection*{Crossing numbers}

Looking at the total number of $1,008$ crossing of the blanket, see Table \ref{table-stats}, raises the question of whether the braids of the blankets are realized with a \emph{minimal} number of crossings. In other words, is there some way to (abstractly) rearrange the braids using the three types of relations from Figure \ref{figure-braidrelations}, see also Equations \eqref{eq-braidrel1}--\eqref{eq-braidrel3}, such that equivalent braids with a smaller number of crossings appear?
\begin{question}\label{question-crossing}
Is the number of $1,008$ crossing in the blanket minimal in representing the braids?
\end{question}

\begin{proof}[Answer to Question \ref{question-crossing}]
It turns out that the numbers of crossings are minimal for all braids $\beta_1,\beta_2,\beta_3$ (of Equations \eqref{pattern1}, \eqref{pattern2}, and \eqref{pattern3}) and all of their products, and hence in the braid $\beta$ symbolizing the entire blanket (see Equation \eqref{eq-blanketbraid}). But to explain why this is the case, powerful theorems are needed.

For the braid $\beta_1$ symbolizing Pattern A, it is easier to check by hand that the number of crossings used is minimal. Recall that in Equation \eqref{Pofbetas}, we say that the permutation associated to $\beta_1$ is $\rP(\beta_1)=(3,1,2)$. In the symmetric group of permutations, the element $(3,1,2)$ can be written as 
$$\rP(\beta_1)=(3,1,2)=(2,1,3)(1,3,2)=\tau_1\tau_2,$$
with the notation for transpositions $\tau_1,\tau_2$ from Equation \eqref{eq-taui}. The element $\rP(\beta_1)$ can be written using at two elementary transpositions ($\tau_1$ and $\tau_1$ once each). It is not possible to write $(3,1,2)$ as a product of just one transposition as it has order $3$ and all transpositions have order $2$. Hence, we need at least two generators (that is, two crossings) to display $\rP(\beta_1)$ in $S_3$. Because all relations among the $\sigma_i$ in $B_n$ are satisfied by the $\tau_i$ in $S_n$, we need at least as many braid group generators $\sigma_i$ to define a braid $\beta$ than we need $\tau_i$ to define the associated permutation $\rP(\beta)$. So, whichever way we try to write $\beta_1$ as a product of elements $\sigma_i$, we will need at least two such generators. The number of generators used in a product equals the number of crossings. So we need at least two crossings to draw $\beta_1$. So the number of crossings used in Pattern A is minimal. 

The same reasoning that we used to argue that $\beta_1$ requires at least two crossings to be drawn cannot be used to determine the minimal number of crossings in $\beta_2$ and $\beta_3$. The reason for this is that in the defining picture for $\beta_2$ (see Figure \ref{figure-patternB}{\scshape (a)}) we use 10 crossings. However, decomposing $\rP(\beta_2)$ gives
\begin{align*}
\rP(\beta_2)=(3,1,2,5,6,4)=(2,1,3,4,5,6)(1,3,2,4,5,6)(1,2,3,4,6,5)(1,2,3,5,4,6)=\tau_1\tau_2\tau_5\tau_4,
\end{align*} 
and only uses four transposition. For this reason, it is possible that there might be a way to display the braid $\beta_2$ with less than 10 crossings. A similar observation can be made for $\beta_3$. As we will see below, we require more clever mathematical methods to be sure about the minimal number of crossings needed to denote the braids $\beta_2$ and $\beta_3$.

\subsubsection*{Positive and homogeneous braids}
A braid is called \emph{positive} if it can be displayed in a way only containing crossings where the left stand crosses over the right strand. In the language of group theory, a braid is positive if it can be written as a product of the generators $\sigma_1,\ldots, \sigma_{n-1}$ without using any inverse elements. As an example, consider the braid $\beta$ in Figure \ref{figure-positive}. This braid reminds us of Pattern B, but it is not quite the same. For instance, comparing to  Figure \ref{figure-patternB}, we see that the crossing of the first two strands are reversed compared to those in Pattern B.

\begin{figure}[htb]
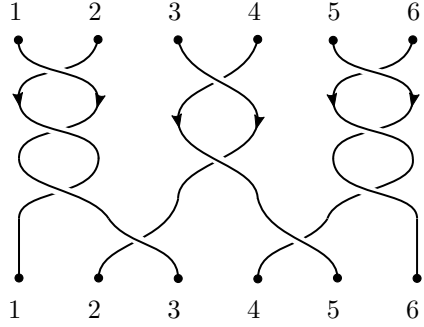
\caption{The positive braid $\sigma_1^3\sigma_3^2\sigma_5^3\sigma_2\sigma_4$}
\label{figure-positive}
\end{figure}

Each positive braid happens to also contain a minimal number of crossings. This is because the only relations between positive braids are given by the 2nd and 3rd types of braid relations Equations \eqref{eq-braidrel2} and \eqref{eq-braidrel3}, depicted in Figure \ref{figure-braidrelations}. These type of relations never reduce the number of generators $\sigma_i$ appearing in a given expression of a positive braid in the generators. Thus, the number of such generators must be minimal. 

None of the braids, $\beta_1,\beta_2,\beta_3$, in the blanket appear to be positive. This can be seen by considering the expressions for $\beta_1$, $\beta_2$, and $\beta_3$ in Equations \eqref{pattern1}, \eqref{pattern2}, and \eqref{pattern3}. In  these braids, we see inverses of the generators $\sigma_i$ appearing. For this reason, we cannot yet determined whether braids $\beta_2$ and $\beta_3$ have a minimal number of crossings. 

Hence, in order to determine whether braids such as $\beta_2$ and $\beta_3$, that are not positive, use a minimal number of crossings, more work by mathematicians is required. By generalizing the concept of a positive braid to larger classes of braids, called \emph{alternating} and \emph{homogeneous} braids, we have the necessary tools at our disposal so that we can be sure that the number of crossings in the blanket's design is minimal. 

The first class of braids that we encounter are \emph{alternating braids}. A braid is called \emph{alternating} if it can be expressed as a product of the elements, $\sigma_1,\sigma_2^{-1},\sigma_3, \sigma_4^{-1}, \ldots$ only. This means that the inverse $\sigma_{2k}^{-1}$ may appear for an even index (that is, an index $2k$ for an integer $k$), but $\sigma_{2k}$ cannot appear. For odd indices (that is, indices $2k+1$ for $k$ an integer), $\sigma_{2k+1}$ may appear, but its inverse $\sigma_{2k+1}^{-1}$ cannot appear. The braid $\beta_1$ from Equation \eqref{pattern1} is an alternating braid. Its definition features only $\sigma_1$ and $\sigma_2^{-1}$. The braid $\beta_2$ is not itself alternating, but its inverse
$$\beta_2^{-1}=\sigma_4^{-1}\sigma_2^{-1}\sigma_{5}\sigma_{3} \sigma_{1}
$$ 
is alternating. 
The result that any alternating braid necessarily has a minimal number of crossings was proved by V. Turaev in \cite{Tur}. 

\begin{theorem}[V.~Turaev, 1988]\label{theorem-Turaev}
Any alternating braid has a minimal crossing number.
\end{theorem}

Therefore, we can use Turaev's theorem to conclude that the braids $\beta_1$ and $\beta_2^{-1}$ have a minimal number of crossings. We already knew this for the braid $\beta_1$ because of the earlier analysis. Moreover, $\beta_2$ itself has a minimal number of crossings. Similarly, all powers of $\beta_1$ and $\beta_2$ are alternating braids and have minimal crossing numbers.

The braid $\beta_3$, however, is neither positive nor alternating, but we observe that its definition in Equation \eqref{pattern3} only features $\sigma_1, \sigma_2^{2}$, and $\sigma_3^{-1}$.  Such a braid is called \emph{$(1,-1,-1)$-homogeneous}. The concept of such a \emph{homogeneous braid} generalizes the idea of both positive and alternating braids. For example, a braid in $B_4$ is positive if and only if it is \emph{$(1,1,1)$-homogeneous}, as it only features $\sigma_1,\sigma_2,\sigma_3$, and $\sigma_4$ but none of their inverses. The string $(1,1,1)$ indices which power of which generator is used. Similarly, an alternating braid such as $\beta_2^{-1}$ in $B_6$ is a \emph{$(1,-1,1,-1,1)$-homogeneous} braids, and $\beta_2$ itself is \emph{$(-1,1,-1,1,-1)$-homogeneous}. As we can see in this blanket, the set of homogeneous braid is larger than the set of positive or alternating braids. It was unknown until recently whether any homogeneous braid has a minimal crossing number, or if there might be braids that are homogeneous but with a crossing number that could be reduced further.

The following theorem is a recent result that was only proved four years ago, see \cite{AM}. 
 
\begin{theorem}[I. Alekseev and G. Mamedov, 2019]\label{theorem-AM}
Any homogeneous braid has a minimal crossing number.
\end{theorem}

Due to this theorem, we know that the braid $\beta_3$ has a minimal crossing number. Therefore, all braids $\beta_1,\beta_2$, and $\beta_3$ have minimal crossing numbers. Further, all powers of these braids have minimal crossing numbers and placing such braids next to one another also produces braids with minimal crossing numbers. Hence, we conclude that the braid $\beta$ from Equation \eqref{eq-blanketbraid} displaying the entire blanket has a minimal crossing number as it is a homogeneous braid itself.
\end{proof}

Question~\ref{question-crossing} turned out to be more complex than Question~\ref{question-pure}. In order to answer Question~\ref{question-crossing} we used recent results from research mathematicians obtained in the last decades, showing how mathematics is an ever-changing field and that there is still much to be discovered.

\section{Braids and knots}\label{epilogue}

The beauty of mathematics reveals itself when connections are created between different kinds of mathematical structures. In order to understand why the theorems of Turaev and Alekseev--Mamedov hold true, we have to understand their \emph{proofs}. Proofs are formal verifications of the validity of a certain mathematical claim---in these cases, the statements of the theorems. For Theorems \ref{theorem-Turaev} \& \ref{theorem-AM}, the results are proved using techniques from the study of \emph{knots}. In mathematics, the concept of a knot can be formalized using a subject called \emph{topology}, where spaces are studied up to continuous deformation. This means that in topology, unlike geometry, distances can be distorted, and a space is only studied as if it was made from soft rubber line, sheets, or blocks.
In this section, we will briefly introduce some key ideas from knot theory and illustrate how studying knots can give us information about the complexity of the braid patterns of the blanket.

\subsection{Knot invariants}\label{section-invariants}

\begin{figure}
\begin{subfigure}[htb]{0.3\textwidth}
\centering
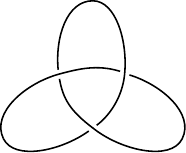
\caption{The knot $\cK_1$}
\end{subfigure}
\begin{subfigure}[htb]{0.3\textwidth}
\centering
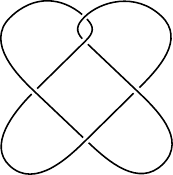
\caption{The knot $\cK_2$}
\end{subfigure}
\begin{subfigure}[htb]{0.3\textwidth}
\centering
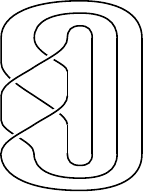
\caption{The knot $\cK_3$}
\end{subfigure}
\caption{Which of these three knots are equivalent?}
\label{figure-knots}
\end{figure}

\begin{figure}
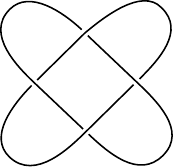
\caption{An example of a link with two components}
\label{figure-link}
\end{figure}

\emph{Knot theory} may be an even older subject of mathematical investigation than braids.\footnote{The first systematic account of the braid group was given by Emil Artin in the 1920s \cite{Art} while knots or links were already studied by Gauss in the 1830s \cite{RN}.} A central goal of knot theory is to list all knots up to so-called \emph{knotting operations}. These operations twist the string of a knot without cutting. In this context, a \emph{knot} is a closed curve in three-dimensional Euclidean space. Three examples of knots are presented in Figure~\ref{figure-knots}. A \emph{link} is a collection of several knots which might be intertwined. An example of such a link, consisting of two smaller, intertwined knots, called its \emph{components}, is displayed in Figure~\ref{figure-link}. A primary objective in knot theory is to classify all knots. The earliest example of such a knot tabulation was given by Peter Guthrie Tait in 1885 by listing all inequivalent knots with up to 10 crossings \cite{MTH}. Providing such a knot table means to write down part of the infinite list of all possible knots ordered by increasing complexity. 

\begin{figure}
\begin{subfigure}[htb]{0.45\textwidth}
\centering
\begingroup%
  \makeatletter%
  \providecommand\color[2][]{%
    \errmessage{(Inkscape) Color is used for the text in Inkscape, but the package 'color.sty' is not loaded}%
    \renewcommand\color[2][]{}%
  }%
  \providecommand\transparent[1]{%
    \errmessage{(Inkscape) Transparency is used (non-zero) for the text in Inkscape, but the package 'transparent.sty' is not loaded}%
    \renewcommand\transparent[1]{}%
  }%
  \providecommand\rotatebox[2]{#2}%
  \newcommand*\fsize{\dimexpr\f@size pt\relax}%
  \newcommand*\lineheight[1]{\fontsize{\fsize}{#1\fsize}\selectfont}%
  \ifx\svgwidth\undefined%
    \setlength{\unitlength}{184.50209766bp}%
    \ifx\svgscale\undefined%
      \relax%
    \else%
      \setlength{\unitlength}{\unitlength * \real{\svgscale}}%
    \fi%
  \else%
    \setlength{\unitlength}{\svgwidth}%
  \fi%
  \global\let\svgwidth\undefined%
  \global\let\svgscale\undefined%
  \makeatother%
  \begin{picture}(1,0.4911852)%
    \lineheight{1}%
    \setlength\tabcolsep{0pt}%
    \put(0,0){\includegraphics[width=\unitlength,page=1]{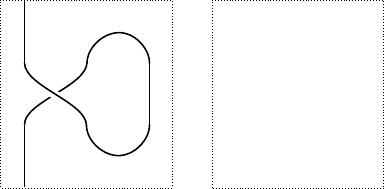}}%
    \put(0.47028201,0.22484261){\makebox(0,0)[lt]{\lineheight{1.25}\smash{\begin{tabular}[t]{l}$=$\end{tabular}}}}%
    \put(0,0){\includegraphics[width=\unitlength,page=2]{reide1.pdf}}%
  \end{picture}%
\endgroup%

\caption{The 1st Reidemeister move}
\end{subfigure}
\begin{subfigure}[htb]{0.45\textwidth}
\centering
\begingroup%
  \makeatletter%
  \providecommand\color[2][]{%
    \errmessage{(Inkscape) Color is used for the text in Inkscape, but the package 'color.sty' is not loaded}%
    \renewcommand\color[2][]{}%
  }%
  \providecommand\transparent[1]{%
    \errmessage{(Inkscape) Transparency is used (non-zero) for the text in Inkscape, but the package 'transparent.sty' is not loaded}%
    \renewcommand\transparent[1]{}%
  }%
  \providecommand\rotatebox[2]{#2}%
  \newcommand*\fsize{\dimexpr\f@size pt\relax}%
  \newcommand*\lineheight[1]{\fontsize{\fsize}{#1\fsize}\selectfont}%
  \ifx\svgwidth\undefined%
    \setlength{\unitlength}{123.94984537bp}%
    \ifx\svgscale\undefined%
      \relax%
    \else%
      \setlength{\unitlength}{\unitlength * \real{\svgscale}}%
    \fi%
  \else%
    \setlength{\unitlength}{\svgwidth}%
  \fi%
  \global\let\svgwidth\undefined%
  \global\let\svgscale\undefined%
  \makeatother%
  \begin{picture}(1,0.4893684)%
    \lineheight{1}%
    \setlength\tabcolsep{0pt}%
    \put(0,0){\includegraphics[width=\unitlength,page=1]{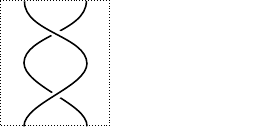}}%
    \put(0.47990234,0.2136643){\makebox(0,0)[lt]{\lineheight{1.25}\smash{\begin{tabular}[t]{l}$=$\end{tabular}}}}%
    \put(0,0){\includegraphics[width=\unitlength,page=2]{reide2.pdf}}%
  \end{picture}%
\endgroup%

\caption{The 2nd Reidemeister move}
\end{subfigure}
\begin{subfigure}[htb]{0.45\textwidth}
\centering
\begingroup%
  \makeatletter%
  \providecommand\color[2][]{%
    \errmessage{(Inkscape) Color is used for the text in Inkscape, but the package 'color.sty' is not loaded}%
    \renewcommand\color[2][]{}%
  }%
  \providecommand\transparent[1]{%
    \errmessage{(Inkscape) Transparency is used (non-zero) for the text in Inkscape, but the package 'transparent.sty' is not loaded}%
    \renewcommand\transparent[1]{}%
  }%
  \providecommand\rotatebox[2]{#2}%
  \newcommand*\fsize{\dimexpr\f@size pt\relax}%
  \newcommand*\lineheight[1]{\fontsize{\fsize}{#1\fsize}\selectfont}%
  \ifx\svgwidth\undefined%
    \setlength{\unitlength}{183.93355561bp}%
    \ifx\svgscale\undefined%
      \relax%
    \else%
      \setlength{\unitlength}{\unitlength * \real{\svgscale}}%
    \fi%
  \else%
    \setlength{\unitlength}{\svgwidth}%
  \fi%
  \global\let\svgwidth\undefined%
  \global\let\svgscale\undefined%
  \makeatother%
  \begin{picture}(1,0.492703)%
    \lineheight{1}%
    \setlength\tabcolsep{0pt}%
    \put(0,0){\includegraphics[width=\unitlength,page=1]{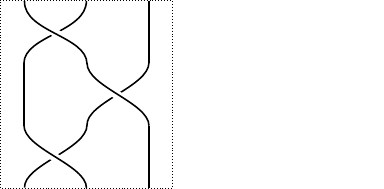}}%
    \put(0.4701907,0.22553714){\makebox(0,0)[lt]{\lineheight{1.25}\smash{\begin{tabular}[t]{l}$=$\end{tabular}}}}%
    \put(0,0){\includegraphics[width=\unitlength,page=2]{reide3.pdf}}%
  \end{picture}%
\endgroup%

\caption{The 3rd Reidemeister move}
\end{subfigure}
\caption{The three fundamental (un)knotting operations called Reidemeister moves}
\label{figure-reidemeister}
\end{figure}

Like we observed in braids, some pictures of knots might look different but describe \emph{equivalent} knots. In this case, one of the knots can be transformed into the other by means of knotting operations that do not cut the string. To give a precise mathematical definition of a knotting operation, we use the concept of \emph{ambient isotopy}. Such an ambient isotopy is, roughly speaking, a globally defined smooth deformation of the ambient space. A fundamental breakthrough in knot theory was achieved by Kurt Reidemeister in the 1920s \cite{Rei}. Reidemeister proved Theorem~\ref{theorem-Reide} below which states that two knots are equivalent through ambient isotopy precisely if they may be deformed into one another using a finite list of the three elementary moves---now called \emph{Reidemeister moves}---displayed in Figure \ref{figure-reidemeister}. In these pictures, the dotted box represents a section inside of a picture of a knot and the moves indicate that the two boxes linked with equality maybe interchanged to produce equivalent knots.

\begin{theorem}[K.~Reidemeister, 1927]\label{theorem-Reide}
Two knots are equivalent if and only if they can be related using a finite sequence of the three types of operations from Figure \ref{figure-reidemeister}.
\end{theorem}

\begin{figure}
\begingroup%
  \makeatletter%
  \providecommand\color[2][]{%
    \errmessage{(Inkscape) Color is used for the text in Inkscape, but the package 'color.sty' is not loaded}%
    \renewcommand\color[2][]{}%
  }%
  \providecommand\transparent[1]{%
    \errmessage{(Inkscape) Transparency is used (non-zero) for the text in Inkscape, but the package 'transparent.sty' is not loaded}%
    \renewcommand\transparent[1]{}%
  }%
  \providecommand\rotatebox[2]{#2}%
  \newcommand*\fsize{\dimexpr\f@size pt\relax}%
  \newcommand*\lineheight[1]{\fontsize{\fsize}{#1\fsize}\selectfont}%
  \ifx\svgwidth\undefined%
    \setlength{\unitlength}{327.28021408bp}%
    \ifx\svgscale\undefined%
      \relax%
    \else%
      \setlength{\unitlength}{\unitlength * \real{\svgscale}}%
    \fi%
  \else%
    \setlength{\unitlength}{\svgwidth}%
  \fi%
  \global\let\svgwidth\undefined%
  \global\let\svgscale\undefined%
  \makeatother%
  \begin{picture}(1,0.22314372)%
    \lineheight{1}%
    \setlength\tabcolsep{0pt}%
    \put(0,0){\includegraphics[width=\unitlength,page=1]{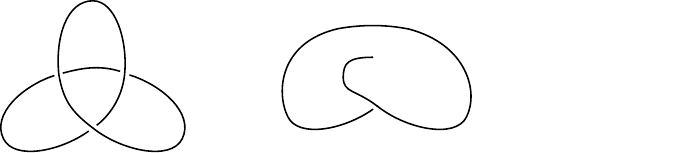}}%
    \put(0.25129364,0.10741542){\makebox(0,0)[lt]{\lineheight{1.25}\smash{\begin{tabular}[t]{l}$\labelb$\end{tabular}}}}%
    \put(0.70656173,0.10760626){\makebox(0,0)[lt]{\lineheight{1.25}\smash{\begin{tabular}[t]{l}$\labela$\end{tabular}}}}%
    \put(0,0){\includegraphics[width=\unitlength,page=2]{equiv-unknot.pdf}}%
  \end{picture}%
\endgroup%

\caption{A knot equivalent to the unknot}
\label{figure-unknotting}
\end{figure}

An example of using the Reidemeister moves from Theorem~\ref{theorem-Reide} to simplify a knot is shown in Figure~\ref{figure-unknotting}. The first knotting operations is of type {\scshape (b)} which uncrosses two neighboring strands. The second knotting operation removes a loop, which is of type {\scshape (c)}. At the end of the manipulation, we are left with the so-called \emph{unknot}, the easiest knot which does not contain any crossing strings and is simply a circle. 

If we are presented with two different knots, sometimes it is impractical to check all the possible ways to deform the ambient space in order to see if the two knots are equivalent. For example, try to confirm which knots Figure \ref{figure-knots} are equivalent and which are not. 
In order to conclusively determine that it is impossible to turn one knot into another through knotting operations, mathematicians employ \emph{invariants}. 
 An invariant is a method to associate a numeric quantity to each knot picture. This quantity could be a simple integer number, or something more complicated, like a polynomial. 
 What defines an invariant is the property stating that \emph{if} two knots are equivalent, \emph{then} the invariant takes the same value for these knots. Using such an invariant, we can affirm that if two knots are associated with different values through the invariant, then they cannot be equivalent.\footnote{This conclusion uses the logical \emph{contrapositive} that the statement ``If A then B'' holds true if and only if the statement ``If not B, then not A'' holds true.} A basic knot invariant is, for example, the number of components in a link which is the number of closed curve inside the link. A knot is a link with just one component and cannot be equivalent to a proper link with several components. Another invariant, which is harder to compute, is the \emph{crossing number}, which is the minimal number of crossings needed to display a knot, using the same idea as the crossing number of a \emph{braid} discussed earlier. The knots $\cK_1$ and $\cK_2$ in Figure \ref{figure-knots} are displayed with a minimal number of crossings---with three and five crossings respectively. In the Alexander--Briggs knot table \cite{AB}, which lists all distinct knots with up to nine crossings, these knots are labelled as the knots $3_1$ and $5_2$ respectively. The knot $\cK_3$ is displayed with four crossings but can be transformed into depictions that only use three crossings. Thus, its crossing number also equals three.
 
Some powerful invariants of knots were found in the 1980 and these invariants of knots were used to prove the theorems about braids that we discussed earlier in this article. These invariants are called the \emph{Jones polynomial} and the \emph{Alexander--Conway polynomial} \cites{Ale,Con,Jon}. For the knots displayed in Figure \ref{figure-knots}, we include the values of the Jones and Alexander--Conway polynomials in Table \ref{table-invariants}.

The beauty of the invariants of Jones and Alexander--Conway is that they can be computed \emph{recursively}, meaning that we can compute the invariant of a more complicated knot or link by knowing the invariants or slightly simpler knot or link (with one or two less crossings) by virtue of a universal formula. In the case of the Alexander--Conway polynomial $\nabla(\cL)$ of a link $\cL$, this formula\footnote{The original \emph{Alexander polynomial} was defined using other methods in the 1920s \cite{Ale}. The recursive formula was discovered by Conway in the late 1960s \cite{Con}.} is given by 
\begin{equation}\label{eq-AC}
\nabla(\cL_+)-\nabla(\cL_-)=z\nabla(\cL_0). 
\end{equation}
Similarly, for the Jones polynomial $\rJ(\cL)$, the formula\footnote{Often, $q=\sqrt{t}$, or $q=1/\sqrt{t}$, is used as the polynomial variable.} is given by 
\begin{equation}\label{eq-Jones}
q^{-2}\rJ(\cL_+)-q^2 \rJ(\cL_-)=(q-q^{-1})\rJ(\cL_0).
\end{equation}
In addition to these recursive formulas, we need to specify the value of the invariants at the simplest possible knot, the unknot. The requirement is simply that 
\begin{equation}\label{eq-normalize}
\nabla(\bigcirc)=1, \qquad \text{and}\qquad \rJ(\bigcirc)=1.
\end{equation}
\begin{figure}[htb]
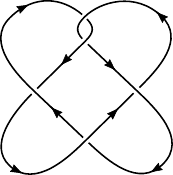
\caption{A choice of an orientation on the knot $\cK_2$ from Figure~\ref{figure-knots}{\scshape (b)}}
\label{figure-oriented}
\end{figure}
To understand the formulas in Equations \ref{eq-AC} and \eqref{eq-Jones}, we have to first give an \emph{orientation} to our link which specifies a direction to travel around each component of the link. An example of a knot with an orientation is shown in Figure~\ref{figure-oriented}. It is not important, which orientation we choose, but we have to continue using one orientation in the recursive computation of the invariants. Next, we fix a single crossing of the knot.  If this crossing is positive, which means it is of the form $\vcenter{\hbox{
\begingroup%
  \makeatletter%
  \providecommand\color[2][]{%
    \errmessage{(Inkscape) Color is used for the text in Inkscape, but the package 'color.sty' is not loaded}%
    \renewcommand\color[2][]{}%
  }%
  \providecommand\transparent[1]{%
    \errmessage{(Inkscape) Transparency is used (non-zero) for the text in Inkscape, but the package 'transparent.sty' is not loaded}%
    \renewcommand\transparent[1]{}%
  }%
  \providecommand\rotatebox[2]{#2}%
  \newcommand*\fsize{\dimexpr\f@size pt\relax}%
  \newcommand*\lineheight[1]{\fontsize{\fsize}{#1\fsize}\selectfont}%
  \ifx\svgwidth\undefined%
    \setlength{\unitlength}{15.75118201bp}%
    \ifx\svgscale\undefined%
      \relax%
    \else%
      \setlength{\unitlength}{\unitlength * \real{\svgscale}}%
    \fi%
  \else%
    \setlength{\unitlength}{\svgwidth}%
  \fi%
  \global\let\svgwidth\undefined%
  \global\let\svgscale\undefined%
  \makeatother%
  \begin{picture}(1,1.00000054)%
    \lineheight{1}%
    \setlength\tabcolsep{0pt}%
    \put(0,0){\includegraphics[width=\unitlength,page=1]{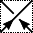}}%
  \end{picture}%
\endgroup%
}}$, we call the link that has this crossing $\cL_+$. Replacing the fixed crossing by the negative crossing $\vcenter{\hbox{
\begingroup%
  \makeatletter%
  \providecommand\color[2][]{%
    \errmessage{(Inkscape) Color is used for the text in Inkscape, but the package 'color.sty' is not loaded}%
    \renewcommand\color[2][]{}%
  }%
  \providecommand\transparent[1]{%
    \errmessage{(Inkscape) Transparency is used (non-zero) for the text in Inkscape, but the package 'transparent.sty' is not loaded}%
    \renewcommand\transparent[1]{}%
  }%
  \providecommand\rotatebox[2]{#2}%
  \newcommand*\fsize{\dimexpr\f@size pt\relax}%
  \newcommand*\lineheight[1]{\fontsize{\fsize}{#1\fsize}\selectfont}%
  \ifx\svgwidth\undefined%
    \setlength{\unitlength}{15.75118201bp}%
    \ifx\svgscale\undefined%
      \relax%
    \else%
      \setlength{\unitlength}{\unitlength * \real{\svgscale}}%
    \fi%
  \else%
    \setlength{\unitlength}{\svgwidth}%
  \fi%
  \global\let\svgwidth\undefined%
  \global\let\svgscale\undefined%
  \makeatother%
  \begin{picture}(1,1.00000054)%
    \lineheight{1}%
    \setlength\tabcolsep{0pt}%
    \put(0,0){\includegraphics[width=\unitlength,page=1]{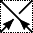}}%
  \end{picture}%
\endgroup%
}}$, gives the knot $\cL_-$. Finally, replacing the original fixed crossing by two parallel strands $\vcenter{\hbox{
\begingroup%
  \makeatletter%
  \providecommand\color[2][]{%
    \errmessage{(Inkscape) Color is used for the text in Inkscape, but the package 'color.sty' is not loaded}%
    \renewcommand\color[2][]{}%
  }%
  \providecommand\transparent[1]{%
    \errmessage{(Inkscape) Transparency is used (non-zero) for the text in Inkscape, but the package 'transparent.sty' is not loaded}%
    \renewcommand\transparent[1]{}%
  }%
  \providecommand\rotatebox[2]{#2}%
  \newcommand*\fsize{\dimexpr\f@size pt\relax}%
  \newcommand*\lineheight[1]{\fontsize{\fsize}{#1\fsize}\selectfont}%
  \ifx\svgwidth\undefined%
    \setlength{\unitlength}{15.75118489bp}%
    \ifx\svgscale\undefined%
      \relax%
    \else%
      \setlength{\unitlength}{\unitlength * \real{\svgscale}}%
    \fi%
  \else%
    \setlength{\unitlength}{\svgwidth}%
  \fi%
  \global\let\svgwidth\undefined%
  \global\let\svgscale\undefined%
  \makeatother%
  \begin{picture}(1,0.99999838)%
    \lineheight{1}%
    \setlength\tabcolsep{0pt}%
    \put(0,0){\includegraphics[width=\unitlength,page=1]{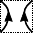}}%
  \end{picture}%
\endgroup%
}}$ gives the knot $\cL_0$. Note that the dotted box refers to a small region inside of a larger knot or link. The following formulas apply Equations~\eqref{eq-AC} and \eqref{eq-Jones} to the knot $\cK_1$ from Figure~\ref{figure-knots}{\scshape (a)}:
\begin{align}\label{eq-AC-appl}
\nabla\left(\resizebox{0.6\width}{!}{$\vcenter{\hbox{
\begingroup%
  \makeatletter%
  \providecommand\color[2][]{%
    \errmessage{(Inkscape) Color is used for the text in Inkscape, but the package 'color.sty' is not loaded}%
    \renewcommand\color[2][]{}%
  }%
  \providecommand\transparent[1]{%
    \errmessage{(Inkscape) Transparency is used (non-zero) for the text in Inkscape, but the package 'transparent.sty' is not loaded}%
    \renewcommand\transparent[1]{}%
  }%
  \providecommand\rotatebox[2]{#2}%
  \newcommand*\fsize{\dimexpr\f@size pt\relax}%
  \newcommand*\lineheight[1]{\fontsize{\fsize}{#1\fsize}\selectfont}%
  \ifx\svgwidth\undefined%
    \setlength{\unitlength}{89.50341702bp}%
    \ifx\svgscale\undefined%
      \relax%
    \else%
      \setlength{\unitlength}{\unitlength * \real{\svgscale}}%
    \fi%
  \else%
    \setlength{\unitlength}{\svgwidth}%
  \fi%
  \global\let\svgwidth\undefined%
  \global\let\svgscale\undefined%
  \makeatother%
  \begin{picture}(1,0.82023259)%
    \lineheight{1}%
    \setlength\tabcolsep{0pt}%
    \put(0,0){\includegraphics[width=\unitlength,page=1]{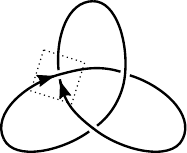}}%
  \end{picture}%
\endgroup%
}}$}\right)-\nabla\left(\resizebox{0.6\width}{!}{$\vcenter{\hbox{
\begingroup%
  \makeatletter%
  \providecommand\color[2][]{%
    \errmessage{(Inkscape) Color is used for the text in Inkscape, but the package 'color.sty' is not loaded}%
    \renewcommand\color[2][]{}%
  }%
  \providecommand\transparent[1]{%
    \errmessage{(Inkscape) Transparency is used (non-zero) for the text in Inkscape, but the package 'transparent.sty' is not loaded}%
    \renewcommand\transparent[1]{}%
  }%
  \providecommand\rotatebox[2]{#2}%
  \newcommand*\fsize{\dimexpr\f@size pt\relax}%
  \newcommand*\lineheight[1]{\fontsize{\fsize}{#1\fsize}\selectfont}%
  \ifx\svgwidth\undefined%
    \setlength{\unitlength}{89.50403517bp}%
    \ifx\svgscale\undefined%
      \relax%
    \else%
      \setlength{\unitlength}{\unitlength * \real{\svgscale}}%
    \fi%
  \else%
    \setlength{\unitlength}{\svgwidth}%
  \fi%
  \global\let\svgwidth\undefined%
  \global\let\svgscale\undefined%
  \makeatother%
  \begin{picture}(1,0.82022692)%
    \lineheight{1}%
    \setlength\tabcolsep{0pt}%
    \put(0,0){\includegraphics[width=\unitlength,page=1]{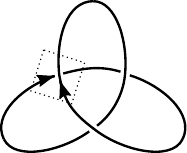}}%
  \end{picture}%
\endgroup%
}}$}\right)&=z \nabla\left(\resizebox{0.6\width}{!}{$\vcenter{\hbox{
\begingroup%
  \makeatletter%
  \providecommand\color[2][]{%
    \errmessage{(Inkscape) Color is used for the text in Inkscape, but the package 'color.sty' is not loaded}%
    \renewcommand\color[2][]{}%
  }%
  \providecommand\transparent[1]{%
    \errmessage{(Inkscape) Transparency is used (non-zero) for the text in Inkscape, but the package 'transparent.sty' is not loaded}%
    \renewcommand\transparent[1]{}%
  }%
  \providecommand\rotatebox[2]{#2}%
  \newcommand*\fsize{\dimexpr\f@size pt\relax}%
  \newcommand*\lineheight[1]{\fontsize{\fsize}{#1\fsize}\selectfont}%
  \ifx\svgwidth\undefined%
    \setlength{\unitlength}{89.56208507bp}%
    \ifx\svgscale\undefined%
      \relax%
    \else%
      \setlength{\unitlength}{\unitlength * \real{\svgscale}}%
    \fi%
  \else%
    \setlength{\unitlength}{\svgwidth}%
  \fi%
  \global\let\svgwidth\undefined%
  \global\let\svgscale\undefined%
  \makeatother%
  \begin{picture}(1,0.81969529)%
    \lineheight{1}%
    \setlength\tabcolsep{0pt}%
    \put(0,0){\includegraphics[width=\unitlength,page=1]{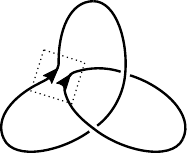}}%
  \end{picture}%
\endgroup%
}}$}\right)
\\ 
q^{-2}\rJ\left(\resizebox{0.6\width}{!}{$\vcenter{\hbox{}}$}\right)-q^2\rJ\left(\resizebox{0.6\width}{!}{$\vcenter{\hbox{}}$}\right)&=(q-q^{-1})\rJ\left(\resizebox{0.6\width}{!}{$\vcenter{\hbox{}}$}\right)
\label{eq-Jones-appl}
\end{align}
There are theoretical reasons why such relatively simple recursive formula provides some of the best known invariants of knots.  These reasons delve into other areas of mathematics such as analysis and algebra. However, the computation of these invariants is not difficult. Equation~\eqref{eq-AC-appl} implies that 
$$\nabla(\cK_1)= 1+z \nabla(\cL_H),$$
where 
\begin{equation}\label{eq-HL}
\cL_H=\vcenter{\hbox{
\begingroup%
  \makeatletter%
  \providecommand\color[2][]{%
    \errmessage{(Inkscape) Color is used for the text in Inkscape, but the package 'color.sty' is not loaded}%
    \renewcommand\color[2][]{}%
  }%
  \providecommand\transparent[1]{%
    \errmessage{(Inkscape) Transparency is used (non-zero) for the text in Inkscape, but the package 'transparent.sty' is not loaded}%
    \renewcommand\transparent[1]{}%
  }%
  \providecommand\rotatebox[2]{#2}%
  \newcommand*\fsize{\dimexpr\f@size pt\relax}%
  \newcommand*\lineheight[1]{\fontsize{\fsize}{#1\fsize}\selectfont}%
  \ifx\svgwidth\undefined%
    \setlength{\unitlength}{54.51488541bp}%
    \ifx\svgscale\undefined%
      \relax%
    \else%
      \setlength{\unitlength}{\unitlength * \real{\svgscale}}%
    \fi%
  \else%
    \setlength{\unitlength}{\svgwidth}%
  \fi%
  \global\let\svgwidth\undefined%
  \global\let\svgscale\undefined%
  \makeatother%
  \begin{picture}(1,0.66850104)%
    \lineheight{1}%
    \setlength\tabcolsep{0pt}%
    \put(0,0){\includegraphics[width=\unitlength,page=1]{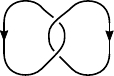}}%
  \end{picture}%
\endgroup%
}}
\end{equation}
is called the \emph{Hopf link}. We have used that the knot $(\cK_1)_-$ in this example is equivalent to the unknot. This can be seen by first using the 2nd Reidemeister move to remove the two crossings at the top, and then using the 1st Reidemeister move to remove the extra loop---see Figure~\ref{figure-unknotting}. Using that $\nabla(\cL_H)=z$, see Equation~\eqref{eq-AC-HL} in the Appendix, we find that 
\begin{equation}\label{eq-AC-K1}
\nabla(\cK_1)=1+z^2.
\end{equation}
Similarly, we derive that 
$$\rJ(\cK_1)=q^4\cdot 1 +(q-q^{-1})\rJ(\cK_H)=q^4+q^2(q-q^{-1})(-q^5-q),$$
using that $\rJ(\cK_H)=-q^5-q$, see Equation~\eqref{eq-Jones-HL}.
From this we compute that 
\begin{align}\label{eq-Jones-K1}
\rJ(\cK_1)=-q^8+q^6+q^2.
\end{align}
We include the detailed computation of the knot invariants associated to the Hopf link and the $3$-twist knot $\cK_2$ in Appendix \ref{appendix}.

\begin{table}[hb]
\[
\begin{array}{c|cc}
\text{Knot $\cK_i$}&\text{Alexander--Conway Polynomial $\nabla(\cK_i)$}&\text{Jones Polynomial $\rJ(\cK_i)$}\\ \hline
\cK_1& 1+z^2 & -q^8+q^6+q^2 \\
\cK_2& 1+2z^2 & -q^{12}+q^{10}-q^8+2q^6-q^4+q^2\\
\cK_3& 1+z^2 & -q^8+q^6+q^2   \end{array}
\]
\caption{The knot invariants associated to the knots $\cK_1,\cK_2,\cK_3$ from Figure \ref{figure-knots}}
\label{table-invariants}
\end{table}

Invariants can be seen as a measure for the complexity of knots and links. A good example for this is the observation that the minimal number of crossings that is required to display a knot or link  is an invariant of knots. For the three examples of knots $\cK_1,\cK_2,\cK_3$ from Figure~\ref{figure-knots}, the knots $\cK_1$ and $\cK_2$ are displayed with a minimal number of crossings. The knot $\cK_1$ requires $3$ crossings, while $\cK_2$ requires $5$ crossings. The knot $\cK_3$ is equivalent to $\cK_1$ and, hence, necessarily has the same minimal crossing numbers, even though the chosen presentation uses $4$ crossings. For other invariants, like the Alexander--Conway and Jones polynomial, there less of a clear interpretation of the information the invariants provide but there appears to be a tendency that more complicated polynomials are associated to more complex knots.

\subsection{Turning braids into knots and links}\label{section-braidstoknots}

\begin{figure}
\begin{center}
\begingroup%
  \makeatletter%
  \providecommand\color[2][]{%
    \errmessage{(Inkscape) Color is used for the text in Inkscape, but the package 'color.sty' is not loaded}%
    \renewcommand\color[2][]{}%
  }%
  \providecommand\transparent[1]{%
    \errmessage{(Inkscape) Transparency is used (non-zero) for the text in Inkscape, but the package 'transparent.sty' is not loaded}%
    \renewcommand\transparent[1]{}%
  }%
  \providecommand\rotatebox[2]{#2}%
  \newcommand*\fsize{\dimexpr\f@size pt\relax}%
  \newcommand*\lineheight[1]{\fontsize{\fsize}{#1\fsize}\selectfont}%
  \ifx\svgwidth\undefined%
    \setlength{\unitlength}{243.18578977bp}%
    \ifx\svgscale\undefined%
      \relax%
    \else%
      \setlength{\unitlength}{\unitlength * \real{\svgscale}}%
    \fi%
  \else%
    \setlength{\unitlength}{\svgwidth}%
  \fi%
  \global\let\svgwidth\undefined%
  \global\let\svgscale\undefined%
  \makeatother%
  \begin{picture}(1,0.49040646)%
    \lineheight{1}%
    \setlength\tabcolsep{0pt}%
    \put(0,0){\includegraphics[width=\unitlength,page=1]{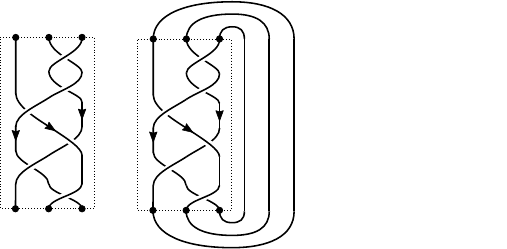}}%
    \put(0.19502451,0.24782922){\makebox(0,0)[lt]{\lineheight{1.25}\smash{\begin{tabular}[t]{l}$\underset{\sim}{\longrightarrow}$\end{tabular}}}}%
    \put(0.61137277,0.24782922){\makebox(0,0)[lt]{\lineheight{1.25}\smash{\begin{tabular}[t]{l}$\underset{\sim}{\longrightarrow}$\end{tabular}}}}%
    \put(0,0){\includegraphics[width=\unitlength,page=2]{braidclosure.pdf}}%
  \end{picture}%
\endgroup%

\end{center}
\caption{The braid closure of the braid $\sigma_2^{-2}\sigma_1^{-1}\sigma_2\sigma_1^{-1}\sigma_2^{-1}$}
\label{figure-braidclosure}
\end{figure}

We will now explore the close connection that exists between braids and links. Given a braid, we can produce a knot or link by connecting each top end of a string to the corresponding end at the bottom. Figure \ref{figure-braidclosure} illustrates this process of closing a braid with an example.  
We can try to check that the knot obtained in this figure is the same as  the knot $\cK_2$ from Figure \ref{figure-knots}{\scshape (b)}. A braid might have several components and is therefore not necessarily a knot but a proper link. The link obtained from the braid $\beta$ is called the \emph{braid closure} of $\beta$. It can be shown that every knot and every link can be obtained as the closure of some braid. 

We can verify by drawing pictures that the knots $\cK_1,\cK_2$, and $\cK_3$ in Figure \ref{figure-knots} can be obtained by closing the braids
$$\gamma_1=\sigma_2^{-1}\sigma_1^{-1}\sigma_2^{-1}\sigma_1^{-1}, \qquad  \gamma_2=\sigma_2^{-2}\sigma_1^{-1}\sigma_2\sigma_1^{-1}\sigma_2^{-1}, \qquad \text{and} \qquad \gamma_3=\sigma_2^{-1}\sigma_1^{-1}\sigma_2^{-1}\sigma_1^{-1}.$$ 
At this point, we gave away that the knots, $\cK_1$ and $\cK_3$, from Figure \ref{figure-knots} are, indeed, equivalent.\footnote{The knots $\cK_1$ and $\cK_3$ are different depictions of the \emph{trefoil knot}. We show that these two knots are equivalent in Figure~\ref{figure-K1equivK3} in Appendix~\ref{appendix}. The knot $\cK_2$ is not equivalent to the other knots---it is known as the \emph{3-twist knot}.} Can we transform one into the other using knotting operations? We observe that braid closures are not unique for a given braid. For example, the braid closure of 
$\sigma_2^{-3}\sigma_1^{-1}\sigma_2\sigma_1^{-1}$ also produces the knot $\cK_2$.
Theorem \ref{theorem-markov} by A.~Markov \cite{Mar} solves the question of when two braid closures are equivalent as knots. Describing knots and links as braid closures is useful, for example, in order to 
implement knots and links in a computer. 

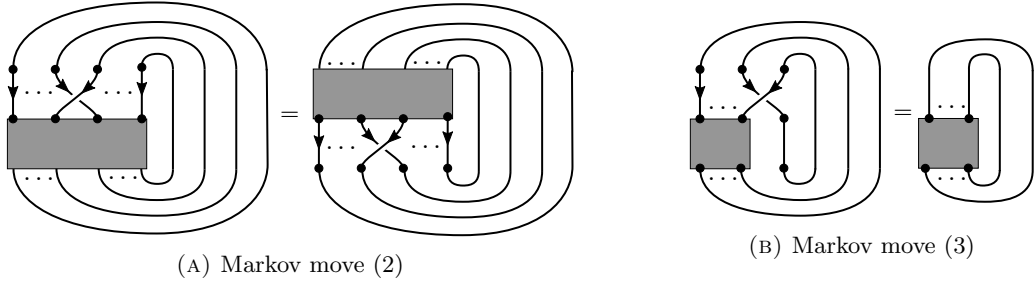
\begin{figure}
\begin{subfigure}[htb]{0.50\textwidth}
\centering
\begingroup%
  \makeatletter%
  \providecommand\color[2][]{%
    \errmessage{(Inkscape) Color is used for the text in Inkscape, but the package 'color.sty' is not loaded}%
    \renewcommand\color[2][]{}%
  }%
  \providecommand\transparent[1]{%
    \errmessage{(Inkscape) Transparency is used (non-zero) for the text in Inkscape, but the package 'transparent.sty' is not loaded}%
    \renewcommand\transparent[1]{}%
  }%
  \providecommand\rotatebox[2]{#2}%
  \newcommand*\fsize{\dimexpr\f@size pt\relax}%
  \newcommand*\lineheight[1]{\fontsize{\fsize}{#1\fsize}\selectfont}%
  \ifx\svgwidth\undefined%
    \setlength{\unitlength}{213.49214221bp}%
    \ifx\svgscale\undefined%
      \relax%
    \else%
      \setlength{\unitlength}{\unitlength * \real{\svgscale}}%
    \fi%
  \else%
    \setlength{\unitlength}{\svgwidth}%
  \fi%
  \global\let\svgwidth\undefined%
  \global\let\svgscale\undefined%
  \makeatother%
  \begin{picture}(1,0.41456731)%
    \lineheight{1}%
    \setlength\tabcolsep{0pt}%
    \put(0,0){\includegraphics[width=\unitlength,page=1]{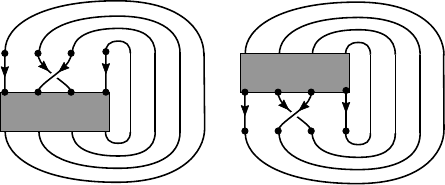}}%
    \put(0.48085728,0.20435224){\makebox(0,0)[lt]{\lineheight{1.25}\smash{\begin{tabular}[t]{l}$=$\end{tabular}}}}%
    \put(0.02905913,0.25064831){\makebox(0,0)[lt]{\lineheight{1.25}\smash{\begin{tabular}[t]{l}$\ldots$\end{tabular}}}}%
    \put(0.16957915,0.25064831){\makebox(0,0)[lt]{\lineheight{1.25}\smash{\begin{tabular}[t]{l}$\ldots$\end{tabular}}}}%
    \put(0.02905913,0.09958879){\makebox(0,0)[lt]{\lineheight{1.25}\smash{\begin{tabular}[t]{l}$\ldots$\end{tabular}}}}%
    \put(0.17823641,0.10046685){\makebox(0,0)[lt]{\lineheight{1.25}\smash{\begin{tabular}[t]{l}$\ldots$\end{tabular}}}}%
    \put(0.56291244,0.15404067){\makebox(0,0)[lt]{\lineheight{1.25}\smash{\begin{tabular}[t]{l}$\ldots$\end{tabular}}}}%
    \put(0.71158717,0.15529516){\makebox(0,0)[lt]{\lineheight{1.25}\smash{\begin{tabular}[t]{l}$\ldots$\end{tabular}}}}%
    \put(0.56353989,0.30208896){\makebox(0,0)[lt]{\lineheight{1.25}\smash{\begin{tabular}[t]{l}$\ldots$\end{tabular}}}}%
    \put(0.7147236,0.30334345){\makebox(0,0)[lt]{\lineheight{1.25}\smash{\begin{tabular}[t]{l}$\ldots$\end{tabular}}}}%
  \end{picture}%
\endgroup%

\caption{Markov move (2)}
\end{subfigure}
\begin{subfigure}[htb]{0.43\textwidth}
\centering
\begingroup%
  \makeatletter%
  \providecommand\color[2][]{%
    \errmessage{(Inkscape) Color is used for the text in Inkscape, but the package 'color.sty' is not loaded}%
    \renewcommand\color[2][]{}%
  }%
  \providecommand\transparent[1]{%
    \errmessage{(Inkscape) Transparency is used (non-zero) for the text in Inkscape, but the package 'transparent.sty' is not loaded}%
    \renewcommand\transparent[1]{}%
  }%
  \providecommand\rotatebox[2]{#2}%
  \newcommand*\fsize{\dimexpr\f@size pt\relax}%
  \newcommand*\lineheight[1]{\fontsize{\fsize}{#1\fsize}\selectfont}%
  \ifx\svgwidth\undefined%
    \setlength{\unitlength}{129.56531386bp}%
    \ifx\svgscale\undefined%
      \relax%
    \else%
      \setlength{\unitlength}{\unitlength * \real{\svgscale}}%
    \fi%
  \else%
    \setlength{\unitlength}{\svgwidth}%
  \fi%
  \global\let\svgwidth\undefined%
  \global\let\svgscale\undefined%
  \makeatother%
  \begin{picture}(1,0.56946817)%
    \lineheight{1}%
    \setlength\tabcolsep{0pt}%
    \put(0,0){\includegraphics[width=\unitlength,page=1]{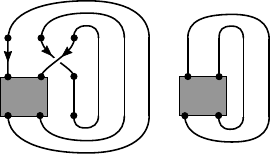}}%
    \put(0.59045903,0.28571123){\makebox(0,0)[lt]{\lineheight{1.25}\smash{\begin{tabular}[t]{l}$=$\end{tabular}}}}%
    \put(0.05194397,0.31465428){\makebox(0,0)[lt]{\lineheight{1.25}\smash{\begin{tabular}[t]{l}$\ldots$\end{tabular}}}}%
    \put(0.05194397,0.10626497){\makebox(0,0)[lt]{\lineheight{1.25}\smash{\begin{tabular}[t]{l}$\ldots$\end{tabular}}}}%
    \put(0.71704324,0.31524216){\makebox(0,0)[lt]{\lineheight{1.25}\smash{\begin{tabular}[t]{l}$\ldots$\end{tabular}}}}%
    \put(0.71704324,0.10685285){\makebox(0,0)[lt]{\lineheight{1.25}\smash{\begin{tabular}[t]{l}$\ldots$\end{tabular}}}}%
  \end{picture}%
\endgroup%

\caption{Markov move (3)}
\end{subfigure}
\caption{Two of Markov's fundamental operations on braid closures}
\label{figure-markov}
\end{figure}

\begin{theorem}[A.~Markov, 1935]\label{theorem-markov}
The braid closures of two braids are equivalent links if and only if the braids are related by a finite list of the following three types of operations:
\begin{enumerate}
\item Changing a braid to an equivalent braid;
\item Changing a braid $\beta$ of the form $\beta=\gamma\sigma_{i}$ to $\beta'=\sigma_{i}\gamma$;
\item Exchanging a braid $\beta$ in $B_n$ with the braid
$(\beta \times 1)\sigma_n$ in $B_{n+1}$.
\end{enumerate}
\end{theorem}
The first operation (1) shows that the braid closure is a well-defined operation. If two braids are equivalent then their closures must be equivalent. The second operation (2) describes moving a crossing from the top of the braid to the bottom by moving the crossing along the unbraided right side of the link. This operation is illustrated in Figure \ref{figure-markov}{\scshape (a)} while the third operation (3) refers to adding (or removing) a loop as in Figure \ref{figure-markov}{\scshape (b)}. In these pictures, an arbitrary braid may be placed into the gray boxes.

\subsection{Invariants of the braid patterns}

We will now apply ideas from knot theory to the analysis of the blanket. Recall that each repeating pattern, Pattern A--C, of the blanket can be described by braid group elements $\beta_1$, $\beta_2$, and $\beta_3$, which were defined in Equations \eqref{pattern1}, \eqref{pattern2}, and \eqref{pattern3}.
We discussed in Section \ref{section-braidstoknots} how the closure of a braid is a link or knot. This means, we can compute the Alexander--Conway polynomial $\nabla$ and Jones polynomial $\rJ$ associated with the links $\cL_{1}$, $\cL_{2}$, and $\cL_{3}$ obtained as the braid closures of $\beta_1$, $\beta_2$, and $\beta_3$, respectively.

\begin{question}\label{question-invariants}
What are the values of Alexander--Conway and Jones polynomials for the braid closures of Pattern A, B, and C of the blanket?
\end{question}

\begin{proof}[Answer to Question~\ref{question-invariants}]
We start by looking at Pattern A. The associated braid $\beta_1$ was identified in Equation \eqref{pattern1}. We see in Figure \ref{figure-braidclosure1} that the closure of this braid can be simplified using two Markov moves of type (3) in Theorem \ref{theorem-markov}.
Therefore, the braid closure of $\beta_1$ is the \emph{unknot} and, by definition, its knot polynomials are equal to $1$. This allows us to conclude that 
\begin{align}
\nabla(\cL_1)=1, && \rJ(\cL_1)=1.
\end{align}
\begin{figure}
\begin{center}
\begingroup%
  \makeatletter%
  \providecommand\color[2][]{%
    \errmessage{(Inkscape) Color is used for the text in Inkscape, but the package 'color.sty' is not loaded}%
    \renewcommand\color[2][]{}%
  }%
  \providecommand\transparent[1]{%
    \errmessage{(Inkscape) Transparency is used (non-zero) for the text in Inkscape, but the package 'transparent.sty' is not loaded}%
    \renewcommand\transparent[1]{}%
  }%
  \providecommand\rotatebox[2]{#2}%
  \newcommand*\fsize{\dimexpr\f@size pt\relax}%
  \newcommand*\lineheight[1]{\fontsize{\fsize}{#1\fsize}\selectfont}%
  \ifx\svgwidth\undefined%
    \setlength{\unitlength}{252.04160304bp}%
    \ifx\svgscale\undefined%
      \relax%
    \else%
      \setlength{\unitlength}{\unitlength * \real{\svgscale}}%
    \fi%
  \else%
    \setlength{\unitlength}{\svgwidth}%
  \fi%
  \global\let\svgwidth\undefined%
  \global\let\svgscale\undefined%
  \makeatother%
  \begin{picture}(1,0.27685313)%
    \lineheight{1}%
    \setlength\tabcolsep{0pt}%
    \put(0,0){\includegraphics[width=\unitlength,page=1]{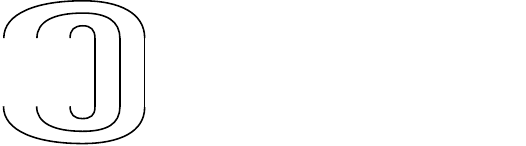}}%
    \put(0.30527246,0.13396037){\makebox(0,0)[lt]{\lineheight{1.25}\smash{\begin{tabular}[t]{l}$\underset{\sim}{\longrightarrow}$\end{tabular}}}}%
    \put(0,0){\includegraphics[width=\unitlength,page=2]{braidclosure1.pdf}}%
    \put(0.70944333,0.13199576){\makebox(0,0)[lt]{\lineheight{1.25}\smash{\begin{tabular}[t]{l}$\underset{\sim}{\longrightarrow}$\end{tabular}}}}%
    \put(0,0){\includegraphics[width=\unitlength,page=3]{braidclosure1.pdf}}%
  \end{picture}%
\endgroup%

\end{center}
\caption{The braid closure of the braid $\beta_1$ symbolizing Pattern A}
\label{figure-braidclosure1}
\end{figure}

\begin{figure}
\begin{center}
\begingroup%
  \makeatletter%
  \providecommand\color[2][]{%
    \errmessage{(Inkscape) Color is used for the text in Inkscape, but the package 'color.sty' is not loaded}%
    \renewcommand\color[2][]{}%
  }%
  \providecommand\transparent[1]{%
    \errmessage{(Inkscape) Transparency is used (non-zero) for the text in Inkscape, but the package 'transparent.sty' is not loaded}%
    \renewcommand\transparent[1]{}%
  }%
  \providecommand\rotatebox[2]{#2}%
  \newcommand*\fsize{\dimexpr\f@size pt\relax}%
  \newcommand*\lineheight[1]{\fontsize{\fsize}{#1\fsize}\selectfont}%
  \ifx\svgwidth\undefined%
    \setlength{\unitlength}{276.78400903bp}%
    \ifx\svgscale\undefined%
      \relax%
    \else%
      \setlength{\unitlength}{\unitlength * \real{\svgscale}}%
    \fi%
  \else%
    \setlength{\unitlength}{\svgwidth}%
  \fi%
  \global\let\svgwidth\undefined%
  \global\let\svgscale\undefined%
  \makeatother%
  \begin{picture}(1,0.33709251)%
    \lineheight{1}%
    \setlength\tabcolsep{0pt}%
    \put(0,0){\includegraphics[width=\unitlength,page=1]{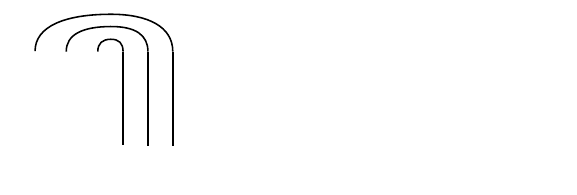}}%
    \put(0.36737587,0.181549){\makebox(0,0)[lt]{\lineheight{1.25}\smash{\begin{tabular}[t]{l}$\underset{\sim}{\longrightarrow}$\end{tabular}}}}%
    \put(0.73541691,0.17976001){\makebox(0,0)[lt]{\lineheight{1.25}\smash{\begin{tabular}[t]{l}$\underset{\sim}{\longrightarrow}$\end{tabular}}}}%
    \put(0,0){\includegraphics[width=\unitlength,page=2]{braidclosure2.pdf}}%
  \end{picture}%
\endgroup%

\end{center}
\caption{The braid closure of the braid $\beta_3$ symbolizing Pattern C}
\label{figure-braidclosure3}
\end{figure}

Next, we consider Pattern C and the braid $\beta_3$ from Equation \eqref{pattern3} extracted from this pattern. Figure~\ref{figure-braidclosure3} shows that we can simplify its braid closure to the Hopf link $\cL_H$ which we have previously encountered in Equation \eqref{eq-HL}. This derivation uses two Markov moves of type (3). The second time, before applying this move, we have to slide the interior part of the link past the outer arch. As this is a topological operation, or, a knotting operation, it can be performed using Reidemeister moves thanks to Theorem~\ref{theorem-Reide}. The Hopf link is a proper link with two components, its knot polynomials are given by 
\begin{align}
\nabla(\cL_3)=\nabla(\cL_H)=z, && \rJ(\cL_3)=\rJ(\cL_H)=-q^5-q.
\end{align}
Details of this computation can be found in Appendix \ref{appendix}.

Pattern~B has a significantly more complicated braid closure. This is due to the fact that the braid $\beta_2$, from Equation \eqref{pattern2}, features higher powers of the final generator, $\sigma_5$, namely,
$$\beta_2=\sigma_{1}^{-3}\sigma_{3}^{-2}\sigma_{5}^{-3}\sigma_2\sigma_4$$
features $\sigma_5^{-3}$. The same is true for the first generator, $\sigma_1$. This makes it impossible to apply a Markov move of type (3) to remove a loop from the braid closure, a move that we made significant use of when identifying the braid closures for $\beta_1$ and $\beta_3$. 
We compute the Alexander--Conway and Jones polynomials of the braid closure $\cL_2$ of $\beta_2$ and display the results in Table~\ref{table-pattern-invariants}.\footnote{This computation used the \emph{KNOT program} of K. Kodama available on this website: \url{http://www.math.kobe-u.ac.jp/~kodama/knot.html}.} As a derivation of these invariants by hand, under use of the recursive formulas in Equations \eqref{eq-AC} and \eqref{eq-Jones} is possible but relatively long, it is not included in this article. Table~\ref{table-pattern-invariants} summarizes the the values of the invariants for the links $\cL_1,\cL_2,$ and $\cL_3$.

\begin{table}[htb]
\[
\begin{array}{c||c|c|c}
\text{Knot $\cK_i$}&\text{Components}&\text{$\nabla(\cL_i)$}&\text{$\rJ(\cL_i)$}\\ \hline
\cL_1& 1 &  1 &  1  \\
\cL_2& 1 &  z^{5}+2z^{3}+z &  -q^{21}+2q^{19}-2q^{17}+4q^{15}-3q^{13}+2q^{11}-3q^{9}-q^{5}   \\
\cL_3& 2 &  z & -q^5-q     \end{array}
\]
\caption{The knot invariants associated to the links $\cL_1,\cL_2,\cL_3$ associated to the blanket's three distinct braid patterns}
\label{table-pattern-invariants}
\end{table}

By inspection, we can tell that the braid closure $\cL_2$ of $\beta_2$ is a link with two components. One way to verify this is by drawing a planar picture of the braid closure. Another method, using group theory, considers the permutation 
$$\rP(\beta_2)=(3,1,2,5,6,4)$$
from Equation \eqref{Pofbetas}. This permutation can be factored as the product 
$$\rP(\beta_2)=(3,1,2,4,5,6)\cdot (1,2,3,5,6,4).$$
The first permutation in this product, $(3,1,2,4,5,6)$ reorders $1,2,3$ but fixes $4,5,6$. However, the second permutation in this product, $(1,2,3,5,6,4)$ fixes $1,2,3$ but reorders $4,5,6$. This way, $\rP(\beta_2)$ separates the numbers from $1$ to $6$ into two disjoint sets: $\{1,2,3\}$ and $\{4,5,6\}$. Repeated application of the partition $\rP(\beta_2)$ connects all points in these subsets, but will never connect a number in the set $\{1,2,3\}$ to a number in the set $\{4,5,6\}$. For instance, starting with $1$, applying $\rP(\beta_2)$ gives us $3$, then $3$ is mapped to $2$, and $2$ is mapped back to $1$. This shows that the two sets $\{1,2,3\}$ and $\{4,5,6\}$ are precisely the subsets of vertices that will constitute components of the braid closure.  
\end{proof}

\begin{figure}
\centering
\begingroup%
  \makeatletter%
  \providecommand\color[2][]{%
    \errmessage{(Inkscape) Color is used for the text in Inkscape, but the package 'color.sty' is not loaded}%
    \renewcommand\color[2][]{}%
  }%
  \providecommand\transparent[1]{%
    \errmessage{(Inkscape) Transparency is used (non-zero) for the text in Inkscape, but the package 'transparent.sty' is not loaded}%
    \renewcommand\transparent[1]{}%
  }%
  \providecommand\rotatebox[2]{#2}%
  \newcommand*\fsize{\dimexpr\f@size pt\relax}%
  \newcommand*\lineheight[1]{\fontsize{\fsize}{#1\fsize}\selectfont}%
  \ifx\svgwidth\undefined%
    \setlength{\unitlength}{307.83696525bp}%
    \ifx\svgscale\undefined%
      \relax%
    \else%
      \setlength{\unitlength}{\unitlength * \real{\svgscale}}%
    \fi%
  \else%
    \setlength{\unitlength}{\svgwidth}%
  \fi%
  \global\let\svgwidth\undefined%
  \global\let\svgscale\undefined%
  \makeatother%
  \begin{picture}(1,0.26895009)%
    \lineheight{1}%
    \setlength\tabcolsep{0pt}%
    \put(0,0){\includegraphics[width=\unitlength,page=1]{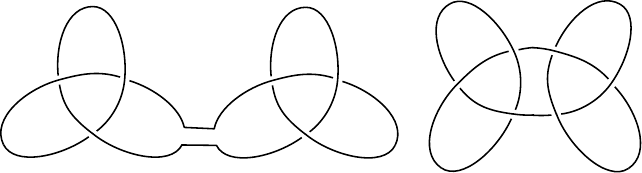}}%
    \put(0.60640709,0.12595149){\color[rgb]{0,0,0}\makebox(0,0)[lt]{\lineheight{1.25}\smash{\begin{tabular}[t]{l}$\underset{\sim}{\longrightarrow}$\end{tabular}}}}%
  \end{picture}%
\endgroup%

\caption{The connected sum of two trefoil knots $\cK_1$}
\label{figure-knotsum}
\end{figure}

We can go one step further and identify the link that appears as the braid closure $\cL_2$ of $\beta_2$. In order to do this, we introduce the operation of \emph{connected sum} of knots. This operation takes two knots, cuts their strings, and ties the open ends of the two knots together. This process results in a single new knot. Figure~\ref{figure-knotsum} shows an example of such a connected sum operation. Here, we computed the connected sum of two identical copies of the trefoil knot $\cK_1$ from Figure~\ref{figure-knots}. The resulting knot is called the \emph{Granny knot}. Given two knots $\cK_1$ and $\cK_2$, we denote their connected sum by $\cK_1\# \cK_2$. For example, the knot in Figure~\ref{figure-knotsum} is denoted by $\cK_1\# \cK_1$ using this notation. 

Taking the connected sum of any knot with the unknot $\bigcirc$ returns the knot itself. In formulas,
$$\cK \# \bigcirc=\cK \qquad \text{and }\qquad  \bigcirc\# \cK.$$
This way, the unknot functions as the identity element for the connected sum operation, similar  to the identity element of a group.
However, the operation $\#$ does \emph{not} define a group product in the sense of Section \ref{section-groups} because there are no inverses. That is, for a given knot $\cK$ we cannot find a knot $\cL$ such that $\cK\#\cL=\bigcirc$ unless both knots are already equivalent to the unknot. The reason for this is that the connected sum adds the minimal crossing numbers of the knots. If $\operatorname{Cr}(\cK)$ denotes the minimal crossing number of the knot $\cK$, then 
$$\operatorname{Cr}(\cK\#\cL)=\operatorname{Cr}(\cK)+\operatorname{Cr}(\cL).$$

A knot that cannot be displayed as a connected sum of two knots (that are not the unknot) is called a \emph{prime knot}. This naming follows the analogy with prime numbers, which cannot be factored as a product of two other positive integers that are not equal to $1$. For instance, the trefoil knot $\cK_1$ is a prime knot. In fact, it is the second easiest knot possible. The granny knot from Figure~\ref{figure-knotsum} is not prime as it is a connected sum of two knots, neither of which are the unknot. 

The connected sum operation is not completely unambiguous for links with several components. It includes a choice of which components of the two link we will cut. However, the Alexander--Conway and Jones polynomials follow the pattern that their value evaluated for connected sum is simply the product of the polynomials. In formulas, this means that for any two links $\cK$ and $\cL$, 
\begin{align}\label{eq-suminvariants}
\nabla(\cK\# \cL)=\nabla(\cK)\cdot \nabla(\cL)&&\text{and}&&
\rJ(\cK\# \cL)=\rJ(\cK)\cdot \rJ(\cL).
\end{align}
These formulas are valid independently of which components of the links were cut in order to form the connected sum.

\begin{figure}
\begin{center}
\resizebox{\textwidth}{!}{
\begingroup%
  \makeatletter%
  \providecommand\color[2][]{%
    \errmessage{(Inkscape) Color is used for the text in Inkscape, but the package 'color.sty' is not loaded}%
    \renewcommand\color[2][]{}%
  }%
  \providecommand\transparent[1]{%
    \errmessage{(Inkscape) Transparency is used (non-zero) for the text in Inkscape, but the package 'transparent.sty' is not loaded}%
    \renewcommand\transparent[1]{}%
  }%
  \providecommand\rotatebox[2]{#2}%
  \newcommand*\fsize{\dimexpr\f@size pt\relax}%
  \newcommand*\lineheight[1]{\fontsize{\fsize}{#1\fsize}\selectfont}%
  \ifx\svgwidth\undefined%
    \setlength{\unitlength}{477.37410405bp}%
    \ifx\svgscale\undefined%
      \relax%
    \else%
      \setlength{\unitlength}{\unitlength * \real{\svgscale}}%
    \fi%
  \else%
    \setlength{\unitlength}{\svgwidth}%
  \fi%
  \global\let\svgwidth\undefined%
  \global\let\svgscale\undefined%
  \makeatother%
  \begin{picture}(1,0.27851566)%
    \lineheight{1}%
    \setlength\tabcolsep{0pt}%
    \put(0,0){\includegraphics[width=\unitlength,page=1]{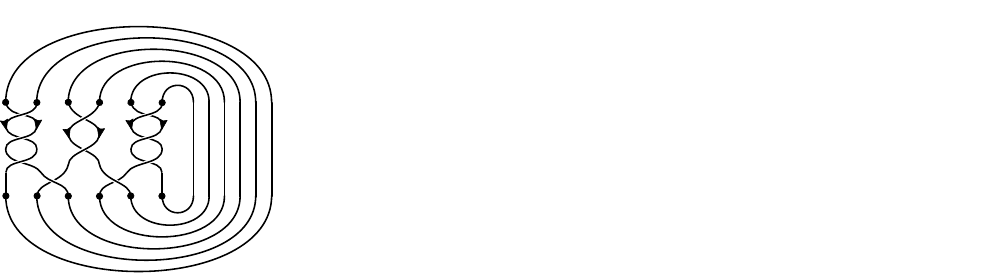}}%
    \put(0.29448693,0.12861431){\makebox(0,0)[lt]{\lineheight{1.25}\smash{\begin{tabular}[t]{l}$\underset{\sim}{\longrightarrow}$\end{tabular}}}}%
    \put(0,0){\includegraphics[width=\unitlength,page=2]{braidclosure3.pdf}}%
    \put(0.65044035,0.12876812){\makebox(0,0)[lt]{\lineheight{1.25}\smash{\begin{tabular}[t]{l}$\underset{\sim}{\longrightarrow}$\end{tabular}}}}%
    \put(0,0){\includegraphics[width=\unitlength,page=3]{braidclosure3.pdf}}%
  \end{picture}%
\endgroup%
}
\end{center}
\caption{The braid closure of the braid $\beta_2$ symbolizing Pattern B}
\label{figure-braidclosure2}
\end{figure}

Applying a few knotting operation shows that $\cL_2$ decomposes as the following connected sum:
\begin{equation}
\cL_2= \cK_1 \# \cL_H \# \cK_1.
 \end{equation}
 Here, $\cK_1$ is the trefoil knot from Figure~\ref{figure-knots}{\scshape (a)}, and $\cL_H$ is, again, the Hopf link. The knotting operations performed in Figure~\ref{figure-braidclosure2} show how this decomposition of $\cL_2$ as a connected sum can be verified. In the figure, we start by closing the braid $\beta_2$. Then, after removing the dots, which are not relevant to the link, we use a few knotting operations to rearrange the knot in a form that makes us recognize it as a connected sum as claimed.
  Now, since the values of the Alexander--Convey polynomial are
 $$\nabla(\cK_1)=z^2+1\qquad \text{and}\qquad \nabla(\cL_H)=z,$$
 by Equations \eqref{eq-AC-K1} and \eqref{eq-Jones-K1}, we can verify that Equation \eqref{eq-suminvariants} is correct in this example. Indeed,
 $$ \nabla(\cL_2)=z^{5}+2z^{3}+z= (z^2+1)z(z^2+1)=\nabla(\cK_1)\nabla(\cL_H)\nabla(\cK_1).$$
Similarly, for the Jones polynomial, we have 
\begin{align*}
\rJ(\cL_2)&=-q^{21}+2q^{19}-2q^{17}+4q^{15}-3q^{13}+2q^{11}-3q^{9}-q^{5}\\
&=(-q^{8}+q^{6}+q^{2})\cdot (-q-q^{5})\cdot (-q^{8}+q^{6}+q^{2})\\
&=\rJ(\cK_1)\rJ(\cL_H)\rJ(\cK_1).
\end{align*}
Finally, we observe that the minimal crossing number of $\cL_2$ equals eight, based on the following computation:
$$\operatorname{Cr}(\cL_2)=\operatorname{Cr}(\cK_1)+\operatorname{Cr}(\cL_H)+\operatorname{Cr}(\cK_1)=3+2+3=8.$$

\section{Conclusion}

In this article, we used mathematics to describe the braid pattern on a knitted blanket. This involved a theory called \emph{group theory}. This way, we were able to answer some questions about the braids, observe some inherit symmetries, and provide measures of the complexity of the structure of the patterns. We were able to demonstrate, by answering Question~\ref{question-pure} that the strands of the braid return to their original position when traced through the entire blanket. This means that the braids featured in the blanket are so-called pure braids. Next, we were able to see that the number of crossings used in the blanket is minimal among all possible equivalent ways to represent these braids by answering Question~\ref{question-crossing}. To answer the second question, we applied results of recent research in mathematics by  Turaev and Alekseev--Mamedov .

We explained the connection between braids and knots in Section~\ref{epilogue}. Studying the knots and links obtained by closing the braid patterns of the blanket lead to Question~\ref{question-invariants} of what values the invariants of Alexander--Conway and Jones assign to these links. While answering this questions, we observed that invariants can serve as a measure of the complexity of a link. For instance, the minimal crossing number of a knot directly indicates its complexity, or, how intertwined the string is. Other invariants like the Jones polynomial tend to associate more complex polynomials to more intricate knots. We were able to identify all three links obtained from the braid patterns among known lists of knots and using the operation of connected sums of links.

Studying the braid pattern of this  blanket shows something remarkable about mathematics. It shows us that mathematics is more than the study of numbers. It is the study of all kinds of structures that we can use to explain the world around us. For example, braids are a fundamental structure that can  describe how strings are intertwined. The study of braids has a long history, but even to this day, new knowledge about braids can be gained. The theorems of Turaev and Alekseev--Mamedov are remarkably recent\footnote{In fact, at the time of writing this text, Alekseev--Mamedov's article had not yet been published in a mathematical journal where their work will be reviewed by other experts in the theory of braids. The process of reviewing an article may take several years. However, nowadays most mathematics articles are made available on the internet as a \emph{pre-print} shortly after completion, so other mathematicians can start studying and applying the author's work.} compared to the age of the mathematical study of braids, which goes back over a hundred years \cite{Mag}. Using a new perspective to look at an old structure may broaden our understanding of fundamental structures such as braids. This shows that mathematics is an evolving area of research. Solving old questions engenders new questions for mathematicians to explore. The question that Alekseev--Mamedov answered was asked by mathematicians only recently after P.~R. Cromwell thought of the concept of a homogeneous braid as a generalization of alternating and positive braids in the late 1980s \cite{Cro}. 

The power of a general mathematical theorem such as Theorem \ref{theorem-AM} allows us to know with certainty that something is true, even if the collection of mathematical objects it refers to can be astronomically large. In theory, there are braids of arbitrary size, so that no person or computer could ever list \emph{all} of them. By only working with fundamental axioms that all of these objects share, mathematicians can derive strong conclusions about all of them at the same time. Sometimes, new mathematical results require connections between different kind of objects such as the the results of Turaev and Alekseev--Mamedov which built on the connection between braids and links, and used the invariants of Jones and Alexander--Conway to prove the results on minimal crossing numbers of braids which we applied to answer Question~\ref{question-crossing}.

Beyond what was explored in this article, invariants of knots have found applications in various other areas of mathematics and even DNA folding \cite{DNA}. 
The Jones polynomial displays the full beauty of mathematics as it relates to several seemingly unrelated areas of research: Mathematical frameworks for quantum physics, through connections to the geometry of three-dimensional spaces and transitions of such spaces over time, and through statistical mechanics; Algebra, through the study of symmetries of so-called \emph{quantum groups} which have emerged from mathematical physics in the late 1980s; Analysis, through the study of so-called subfactors of von Neumann algebras. 
These connections continue to inspire a wealth of research in mathematics and physics to this date. More information about these connections can, for example, be found in the research monograph \cite{JonBook} by Vaughan Jones who discovered the Jones Polynomial and received the \emph{Fields medal}, an award regarded as the Nobel prize for mathematicians, for his work.

\appendix
\section{Computation of knot invariants}\label{appendix}

In this appendix, we include detailed computations of the knot invariants listed in Table~\ref{table-invariants}. Recall that for a single circle, both invariants $\nabla(\bigcirc)$ and $\rJ(\bigcirc)$ take the value $1$. From this, we can derive the value of these invariants on two copies of the circle. For this, we consider the following three links:
\begin{align*}
\cL_+=\vcenter{\hbox{
\begingroup%
  \makeatletter%
  \providecommand\color[2][]{%
    \errmessage{(Inkscape) Color is used for the text in Inkscape, but the package 'color.sty' is not loaded}%
    \renewcommand\color[2][]{}%
  }%
  \providecommand\transparent[1]{%
    \errmessage{(Inkscape) Transparency is used (non-zero) for the text in Inkscape, but the package 'transparent.sty' is not loaded}%
    \renewcommand\transparent[1]{}%
  }%
  \providecommand\rotatebox[2]{#2}%
  \newcommand*\fsize{\dimexpr\f@size pt\relax}%
  \newcommand*\lineheight[1]{\fontsize{\fsize}{#1\fsize}\selectfont}%
  \ifx\svgwidth\undefined%
    \setlength{\unitlength}{51.37499903bp}%
    \ifx\svgscale\undefined%
      \relax%
    \else%
      \setlength{\unitlength}{\unitlength * \real{\svgscale}}%
    \fi%
  \else%
    \setlength{\unitlength}{\svgwidth}%
  \fi%
  \global\let\svgwidth\undefined%
  \global\let\svgscale\undefined%
  \makeatother%
  \begin{picture}(1,0.41608144)%
    \lineheight{1}%
    \setlength\tabcolsep{0pt}%
    \put(0,0){\includegraphics[width=\unitlength,page=1]{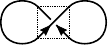}}%
  \end{picture}%
\endgroup%
}}, \quad \cL_-=\vcenter{\hbox{
\begingroup%
  \makeatletter%
  \providecommand\color[2][]{%
    \errmessage{(Inkscape) Color is used for the text in Inkscape, but the package 'color.sty' is not loaded}%
    \renewcommand\color[2][]{}%
  }%
  \providecommand\transparent[1]{%
    \errmessage{(Inkscape) Transparency is used (non-zero) for the text in Inkscape, but the package 'transparent.sty' is not loaded}%
    \renewcommand\transparent[1]{}%
  }%
  \providecommand\rotatebox[2]{#2}%
  \newcommand*\fsize{\dimexpr\f@size pt\relax}%
  \newcommand*\lineheight[1]{\fontsize{\fsize}{#1\fsize}\selectfont}%
  \ifx\svgwidth\undefined%
    \setlength{\unitlength}{51.37499903bp}%
    \ifx\svgscale\undefined%
      \relax%
    \else%
      \setlength{\unitlength}{\unitlength * \real{\svgscale}}%
    \fi%
  \else%
    \setlength{\unitlength}{\svgwidth}%
  \fi%
  \global\let\svgwidth\undefined%
  \global\let\svgscale\undefined%
  \makeatother%
  \begin{picture}(1,0.41608144)%
    \lineheight{1}%
    \setlength\tabcolsep{0pt}%
    \put(0,0){\includegraphics[width=\unitlength,page=1]{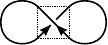}}%
  \end{picture}%
\endgroup%
}}, \quad \cL_0=\vcenter{\hbox{
\begingroup%
  \makeatletter%
  \providecommand\color[2][]{%
    \errmessage{(Inkscape) Color is used for the text in Inkscape, but the package 'color.sty' is not loaded}%
    \renewcommand\color[2][]{}%
  }%
  \providecommand\transparent[1]{%
    \errmessage{(Inkscape) Transparency is used (non-zero) for the text in Inkscape, but the package 'transparent.sty' is not loaded}%
    \renewcommand\transparent[1]{}%
  }%
  \providecommand\rotatebox[2]{#2}%
  \newcommand*\fsize{\dimexpr\f@size pt\relax}%
  \newcommand*\lineheight[1]{\fontsize{\fsize}{#1\fsize}\selectfont}%
  \ifx\svgwidth\undefined%
    \setlength{\unitlength}{51.37499903bp}%
    \ifx\svgscale\undefined%
      \relax%
    \else%
      \setlength{\unitlength}{\unitlength * \real{\svgscale}}%
    \fi%
  \else%
    \setlength{\unitlength}{\svgwidth}%
  \fi%
  \global\let\svgwidth\undefined%
  \global\let\svgscale\undefined%
  \makeatother%
  \begin{picture}(1,0.41608144)%
    \lineheight{1}%
    \setlength\tabcolsep{0pt}%
    \put(0,0){\includegraphics[width=\unitlength,page=1]{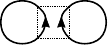}}%
  \end{picture}%
\endgroup%
}}.
\end{align*}
Now, the defining equation for the Alexander--Conway polynomial in Equation~\eqref{eq-AC} gives
\begin{align*}
z\nabla(\cL_0)= \nabla(\cL_+)- \nabla(\cL_-)=\nabla(\bigcirc)- \nabla(\bigcirc)=1-1=0,
\end{align*}
since both $\cL_+$ and $\cL_-$ are equivalent to the unknot, using the 1st Reidemeister move from Figure~\ref{figure-reidemeister}{\scshape (a)}. Recall that the value of the knot invariant $\nabla$ does not depend on the chosen orientation, which merely serves as a tool to apply the recursive formula. Thus, we conclude that
\begin{equation}\label{eq-AC-2circ}
\nabla\left(\bigcirc \; \bigcirc \right)=0.
\end{equation}
Similarly, we can compute the Jones polynomial of two copies of the unknot using Equation~\eqref{eq-Jones}. Namely, 
\begin{align*}
\rJ(\cL_0)=\frac{1}{q-q^{-1}}\left(q^{-2}\rJ(\cL_+)-q^2 \rJ(\cL_-)\right)=\frac{q^{-2}-q^2}{q-q^{-1}}=\frac{(q^{-1}-q)(q^{-1}+q)}{q-q^{-1}}.
\end{align*}
Thus, we conclude that 
\begin{equation}\label{eq-Jones-2circ}
\rJ\left(\bigcirc \; \bigcirc \right)=-q-q^{-1}.
\end{equation}

Now we can move on to compute the Alexander--Conway and Jones polynomial for the \emph{Hopf link} $\cL_H$ which is displayed in Equation~\eqref{eq-HL}. In this case, by fixing a single crossing, we obtain the following three links:
\begin{align*}
(\cL_H)_+=\vcenter{\hbox{
\begingroup%
  \makeatletter%
  \providecommand\color[2][]{%
    \errmessage{(Inkscape) Color is used for the text in Inkscape, but the package 'color.sty' is not loaded}%
    \renewcommand\color[2][]{}%
  }%
  \providecommand\transparent[1]{%
    \errmessage{(Inkscape) Transparency is used (non-zero) for the text in Inkscape, but the package 'transparent.sty' is not loaded}%
    \renewcommand\transparent[1]{}%
  }%
  \providecommand\rotatebox[2]{#2}%
  \newcommand*\fsize{\dimexpr\f@size pt\relax}%
  \newcommand*\lineheight[1]{\fontsize{\fsize}{#1\fsize}\selectfont}%
  \ifx\svgwidth\undefined%
    \setlength{\unitlength}{54.52825772bp}%
    \ifx\svgscale\undefined%
      \relax%
    \else%
      \setlength{\unitlength}{\unitlength * \real{\svgscale}}%
    \fi%
  \else%
    \setlength{\unitlength}{\svgwidth}%
  \fi%
  \global\let\svgwidth\undefined%
  \global\let\svgscale\undefined%
  \makeatother%
  \begin{picture}(1,0.68001086)%
    \lineheight{1}%
    \setlength\tabcolsep{0pt}%
    \put(0,0){\includegraphics[width=\unitlength,page=1]{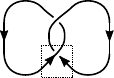}}%
  \end{picture}%
\endgroup%
}}, \quad (\cL_H)_-=\vcenter{\hbox{
\begingroup%
  \makeatletter%
  \providecommand\color[2][]{%
    \errmessage{(Inkscape) Color is used for the text in Inkscape, but the package 'color.sty' is not loaded}%
    \renewcommand\color[2][]{}%
  }%
  \providecommand\transparent[1]{%
    \errmessage{(Inkscape) Transparency is used (non-zero) for the text in Inkscape, but the package 'transparent.sty' is not loaded}%
    \renewcommand\transparent[1]{}%
  }%
  \providecommand\rotatebox[2]{#2}%
  \newcommand*\fsize{\dimexpr\f@size pt\relax}%
  \newcommand*\lineheight[1]{\fontsize{\fsize}{#1\fsize}\selectfont}%
  \ifx\svgwidth\undefined%
    \setlength{\unitlength}{121.93656537bp}%
    \ifx\svgscale\undefined%
      \relax%
    \else%
      \setlength{\unitlength}{\unitlength * \real{\svgscale}}%
    \fi%
  \else%
    \setlength{\unitlength}{\svgwidth}%
  \fi%
  \global\let\svgwidth\undefined%
  \global\let\svgscale\undefined%
  \makeatother%
  \begin{picture}(1,0.30409096)%
    \lineheight{1}%
    \setlength\tabcolsep{0pt}%
    \put(0,0){\includegraphics[width=\unitlength,page=1]{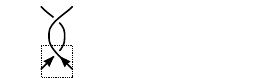}}%
    \put(0.4695006,0.13660346){\makebox(0,0)[lt]{\lineheight{1.25}\smash{\begin{tabular}[t]{l}$\to$\end{tabular}}}}%
    \put(0,0){\includegraphics[width=\unitlength,page=2]{HL-.pdf}}%
  \end{picture}%
\endgroup%
}}, \quad (\cL_H)_0=\vcenter{\hbox{
\begingroup%
  \makeatletter%
  \providecommand\color[2][]{%
    \errmessage{(Inkscape) Color is used for the text in Inkscape, but the package 'color.sty' is not loaded}%
    \renewcommand\color[2][]{}%
  }%
  \providecommand\transparent[1]{%
    \errmessage{(Inkscape) Transparency is used (non-zero) for the text in Inkscape, but the package 'transparent.sty' is not loaded}%
    \renewcommand\transparent[1]{}%
  }%
  \providecommand\rotatebox[2]{#2}%
  \newcommand*\fsize{\dimexpr\f@size pt\relax}%
  \newcommand*\lineheight[1]{\fontsize{\fsize}{#1\fsize}\selectfont}%
  \ifx\svgwidth\undefined%
    \setlength{\unitlength}{94.69390997bp}%
    \ifx\svgscale\undefined%
      \relax%
    \else%
      \setlength{\unitlength}{\unitlength * \real{\svgscale}}%
    \fi%
  \else%
    \setlength{\unitlength}{\svgwidth}%
  \fi%
  \global\let\svgwidth\undefined%
  \global\let\svgscale\undefined%
  \makeatother%
  \begin{picture}(1,0.39368343)%
    \lineheight{1}%
    \setlength\tabcolsep{0pt}%
    \put(0,0){\includegraphics[width=\unitlength,page=1]{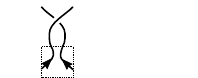}}%
    \put(0.60721437,0.17566106){\makebox(0,0)[lt]{\lineheight{1.25}\smash{\begin{tabular}[t]{l}$\to$\end{tabular}}}}%
    \put(0,0){\includegraphics[width=\unitlength,page=2]{HL0.pdf}}%
  \end{picture}%
\endgroup%
}}.
\end{align*}
Here, we use the 2nd Reidemeister move to transform $(\cL_H)_-$ into two circles, and the 1st Reidemeister move to transform $(\cL_H)_0$ into the unknot. From Equation~\eqref{eq-AC} we can now derive that 
$$\nabla(\cL_H)=\nabla(\bigcirc \; \bigcirc)+z \nabla(\bigcirc)=0+z,$$
under use of Equation~\eqref{eq-Jones-2circ} and hence conclude that
\begin{equation}\label{eq-AC-HL}
\nabla\left(\vcenter{\hbox{}} \right)=z.
\end{equation}
A similar computation for the Jones polynomial of Equation~\eqref{eq-Jones} gives
$$
\rJ(\cL_H)=q^4 \rJ(\bigcirc \; \bigcirc)+q^2(q-q^{-1})\rJ(\bigcirc)=-q^5-q^3+q^3-q,
$$
where we use Equation~\eqref{eq-Jones-2circ} for the value of $\rJ(\bigcirc \; \bigcirc)$.
This yields
\begin{equation}\label{eq-Jones-HL}
\rJ\left(\vcenter{\hbox{
\begingroup%
  \makeatletter%
  \providecommand\color[2][]{%
    \errmessage{(Inkscape) Color is used for the text in Inkscape, but the package 'color.sty' is not loaded}%
    \renewcommand\color[2][]{}%
  }%
  \providecommand\transparent[1]{%
    \errmessage{(Inkscape) Transparency is used (non-zero) for the text in Inkscape, but the package 'transparent.sty' is not loaded}%
    \renewcommand\transparent[1]{}%
  }%
  \providecommand\rotatebox[2]{#2}%
  \newcommand*\fsize{\dimexpr\f@size pt\relax}%
  \newcommand*\lineheight[1]{\fontsize{\fsize}{#1\fsize}\selectfont}%
  \ifx\svgwidth\undefined%
    \setlength{\unitlength}{51.36276612bp}%
    \ifx\svgscale\undefined%
      \relax%
    \else%
      \setlength{\unitlength}{\unitlength * \real{\svgscale}}%
    \fi%
  \else%
    \setlength{\unitlength}{\svgwidth}%
  \fi%
  \global\let\svgwidth\undefined%
  \global\let\svgscale\undefined%
  \makeatother%
  \begin{picture}(1,0.70952677)%
    \lineheight{1}%
    \setlength\tabcolsep{0pt}%
    \put(0,0){\includegraphics[width=\unitlength,page=1]{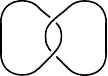}}%
  \end{picture}%
\endgroup%
}}\right)=-q^5-q.
\end{equation}
At this point, we have completed all calculations that were used to derive the values $\nabla(\cK_1)$ and $\rJ(\cK_1)$ in Section~\ref{section-invariants}, see Equations~\eqref{eq-AC-K1} and \eqref{eq-Jones-K1}.

We can now compute the invariants associated to the knot $\cK_2$. We choose a positive crossing at the bottom of the knot and, by replacing this crossing by a negative crossing and two parallel strands (to find $(\cK_2)_-$ and $(\cK_2)_0$) we obtain:
\begin{align*}
(\cK_2)_+=\vcenter{\hbox{
\begingroup%
  \makeatletter%
  \providecommand\color[2][]{%
    \errmessage{(Inkscape) Color is used for the text in Inkscape, but the package 'color.sty' is not loaded}%
    \renewcommand\color[2][]{}%
  }%
  \providecommand\transparent[1]{%
    \errmessage{(Inkscape) Transparency is used (non-zero) for the text in Inkscape, but the package 'transparent.sty' is not loaded}%
    \renewcommand\transparent[1]{}%
  }%
  \providecommand\rotatebox[2]{#2}%
  \newcommand*\fsize{\dimexpr\f@size pt\relax}%
  \newcommand*\lineheight[1]{\fontsize{\fsize}{#1\fsize}\selectfont}%
  \ifx\svgwidth\undefined%
    \setlength{\unitlength}{48.09253224bp}%
    \ifx\svgscale\undefined%
      \relax%
    \else%
      \setlength{\unitlength}{\unitlength * \real{\svgscale}}%
    \fi%
  \else%
    \setlength{\unitlength}{\svgwidth}%
  \fi%
  \global\let\svgwidth\undefined%
  \global\let\svgscale\undefined%
  \makeatother%
  \begin{picture}(1,1.00568712)%
    \lineheight{1}%
    \setlength\tabcolsep{0pt}%
    \put(0,0){\includegraphics[width=\unitlength,page=1]{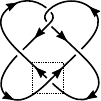}}%
  \end{picture}%
\endgroup%
}}, \quad (\cK_2)_-=\vcenter{\hbox{
\begingroup%
  \makeatletter%
  \providecommand\color[2][]{%
    \errmessage{(Inkscape) Color is used for the text in Inkscape, but the package 'color.sty' is not loaded}%
    \renewcommand\color[2][]{}%
  }%
  \providecommand\transparent[1]{%
    \errmessage{(Inkscape) Transparency is used (non-zero) for the text in Inkscape, but the package 'transparent.sty' is not loaded}%
    \renewcommand\transparent[1]{}%
  }%
  \providecommand\rotatebox[2]{#2}%
  \newcommand*\fsize{\dimexpr\f@size pt\relax}%
  \newcommand*\lineheight[1]{\fontsize{\fsize}{#1\fsize}\selectfont}%
  \ifx\svgwidth\undefined%
    \setlength{\unitlength}{117.31095813bp}%
    \ifx\svgscale\undefined%
      \relax%
    \else%
      \setlength{\unitlength}{\unitlength * \real{\svgscale}}%
    \fi%
  \else%
    \setlength{\unitlength}{\svgwidth}%
  \fi%
  \global\let\svgwidth\undefined%
  \global\let\svgscale\undefined%
  \makeatother%
  \begin{picture}(1,0.41228977)%
    \lineheight{1}%
    \setlength\tabcolsep{0pt}%
    \put(0,0){\includegraphics[width=\unitlength,page=1]{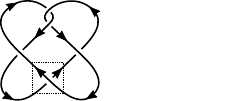}}%
    \put(0.41951788,0.19039976){\makebox(0,0)[lt]{\lineheight{1.25}\smash{\begin{tabular}[t]{l}$\to$\end{tabular}}}}%
    \put(0,0){\includegraphics[width=\unitlength,page=2]{K2-.pdf}}%
  \end{picture}%
\endgroup%
}}, \quad (\cK_2)_0=\vcenter{\hbox{
\begingroup%
  \makeatletter%
  \providecommand\color[2][]{%
    \errmessage{(Inkscape) Color is used for the text in Inkscape, but the package 'color.sty' is not loaded}%
    \renewcommand\color[2][]{}%
  }%
  \providecommand\transparent[1]{%
    \errmessage{(Inkscape) Transparency is used (non-zero) for the text in Inkscape, but the package 'transparent.sty' is not loaded}%
    \renewcommand\transparent[1]{}%
  }%
  \providecommand\rotatebox[2]{#2}%
  \newcommand*\fsize{\dimexpr\f@size pt\relax}%
  \newcommand*\lineheight[1]{\fontsize{\fsize}{#1\fsize}\selectfont}%
  \ifx\svgwidth\undefined%
    \setlength{\unitlength}{116.2209603bp}%
    \ifx\svgscale\undefined%
      \relax%
    \else%
      \setlength{\unitlength}{\unitlength * \real{\svgscale}}%
    \fi%
  \else%
    \setlength{\unitlength}{\svgwidth}%
  \fi%
  \global\let\svgwidth\undefined%
  \global\let\svgscale\undefined%
  \makeatother%
  \begin{picture}(1,0.41602761)%
    \lineheight{1}%
    \setlength\tabcolsep{0pt}%
    \put(0,0){\includegraphics[width=\unitlength,page=1]{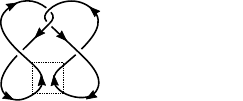}}%
    \put(0.42172406,0.17684548){\makebox(0,0)[lt]{\lineheight{1.25}\smash{\begin{tabular}[t]{l}$\to$\end{tabular}}}}%
    \put(0,0){\includegraphics[width=\unitlength,page=2]{K20.pdf}}%
  \end{picture}%
\endgroup%
}}.
\end{align*}
We see that $(\cK_2)_0$ is equivalent to a Hopf link $\cL_H$. One can also check, see Figure~\ref{figure-K1equivK3}, that $(\cK_2)-$ is equivalent to the trefoil knot $\cK_1$. This allows us to compute the invariants associated to the knot $\cK_2$. For the Alexander--Conway polynomial, we compute, again using Equation~\eqref{eq-AC}, that
\begin{align*}
\nabla(\cK_2)=\nabla((\cK_2)_+)&=\nabla((\cK_2)_-)+z \nabla((\cK_2)_0)\\
&=\nabla(\cK_1)+z \nabla(\cL_H)\\
&=1+z^2+z z,
\end{align*}
where we apply that $\nabla(\cK_1)=1+z^2$ by Equation~\eqref{eq-AC-K1} and $\nabla(\cL_H)=z$ by Equation~\eqref{eq-AC-HL}. Thus, we conclude that 
\begin{equation}\label{eq-AC-K2}
\nabla(\cK_2)=1+2z^2.
\end{equation}
Similarly, we compute the Jones polynomial using Equation~\eqref{eq-Jones}:
\begin{align*}
\rJ(\cK_2)=\rJ((\cK_2)_+)&=q^4\rJ((\cK_2)_-)+q^2(q-q^{-1}) \rJ((\cK_2)_0)\\
&=q^4\rJ(\cK_1)+(q^3-q) \rJ(\cL_H)\\
&=q^4(-q^8+q^6+q^2)+(q^3-q)(-q^5-q),
\end{align*}
where we apply that $\rJ(\cK_1)=-q^8+q^6+q^2$ by Equation~\eqref{eq-Jones-K1} and $\rJ(\cL_H)=-q^5-q$ by Equation~\eqref{eq-Jones-HL}. Summarizing and simplifying the powers of $q$ we conclude that 
\begin{equation}\label{eq-Jones-K2}
\rJ(\cK_2)=-q^{12}+q^{10}-q^8+2q^6-q^4+q^2.
\end{equation}
To complete the computation of knot invariants in Table~\ref{table-invariants}, we show that $\cK_1$ and $\cK_3$ are equivalent knots in Figure~\ref{figure-B6}. This way, $\nabla(\cK_1)=\nabla(\cK_3)$ and $\rJ(\cK_1)=\rJ(\cK_3)$. 

\begin{figure}
\centering
\begingroup%
  \makeatletter%
  \providecommand\color[2][]{%
    \errmessage{(Inkscape) Color is used for the text in Inkscape, but the package 'color.sty' is not loaded}%
    \renewcommand\color[2][]{}%
  }%
  \providecommand\transparent[1]{%
    \errmessage{(Inkscape) Transparency is used (non-zero) for the text in Inkscape, but the package 'transparent.sty' is not loaded}%
    \renewcommand\transparent[1]{}%
  }%
  \providecommand\rotatebox[2]{#2}%
  \newcommand*\fsize{\dimexpr\f@size pt\relax}%
  \newcommand*\lineheight[1]{\fontsize{\fsize}{#1\fsize}\selectfont}%
  \ifx\svgwidth\undefined%
    \setlength{\unitlength}{338.26539198bp}%
    \ifx\svgscale\undefined%
      \relax%
    \else%
      \setlength{\unitlength}{\unitlength * \real{\svgscale}}%
    \fi%
  \else%
    \setlength{\unitlength}{\svgwidth}%
  \fi%
  \global\let\svgwidth\undefined%
  \global\let\svgscale\undefined%
  \makeatother%
  \begin{picture}(1,0.17844428)%
    \lineheight{1}%
    \setlength\tabcolsep{0pt}%
    \put(0,0){\includegraphics[width=\unitlength,page=1]{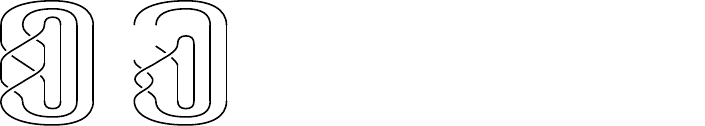}}%
    \put(0.14595759,0.0872639){\makebox(0,0)[lt]{\lineheight{1.25}\smash{\begin{tabular}[t]{l}$\sim$\end{tabular}}}}%
    \put(0,0){\includegraphics[width=\unitlength,page=2]{K1equivK3.pdf}}%
    \put(0.33228263,0.08752565){\makebox(0,0)[lt]{\lineheight{1.25}\smash{\begin{tabular}[t]{l}$\sim$\end{tabular}}}}%
    \put(0,0){\includegraphics[width=\unitlength,page=3]{K1equivK3.pdf}}%
    \put(0.4663436,0.08712977){\makebox(0,0)[lt]{\lineheight{1.25}\smash{\begin{tabular}[t]{l}$\sim$\end{tabular}}}}%
    \put(0,0){\includegraphics[width=\unitlength,page=4]{K1equivK3.pdf}}%
    \put(0.66550963,0.08817074){\makebox(0,0)[lt]{\lineheight{1.25}\smash{\begin{tabular}[t]{l}$\sim$\\\end{tabular}}}}%
    \put(0,0){\includegraphics[width=\unitlength,page=5]{K1equivK3.pdf}}%
    \put(0.82051523,0.08627037){\makebox(0,0)[lt]{\lineheight{1.25}\smash{\begin{tabular}[t]{l}$\sim$\\\end{tabular}}}}%
    \put(0,0){\includegraphics[width=\unitlength,page=6]{K1equivK3.pdf}}%
  \end{picture}%
\endgroup%

\caption{An equivalence of the two knots $\cK_1$ and $\cK_3$ from Figure \ref{figure-knots}}
\label{figure-K1equivK3}
\end{figure}



\bibliography{biblio}
\bibliographystyle{amsrefs}

\end{document}